\documentclass[12pt,leqno]{amsart}
\usepackage{amsmath,stmaryrd}
\usepackage{amssymb}
\usepackage{amsthm}
\usepackage{epsf}
\usepackage{graphicx}
\usepackage{esint} 
\usepackage{dsfont}
\usepackage[pagebackref]{hyperref} 
\hypersetup{pdfpagemode=FullScreen,  colorlinks=true} 
\usepackage{enumerate} 
\usepackage{xcolor,bbm,easybmat,bm}
\usepackage{tikz-cd}
\usepackage{centernot}

\numberwithin{equation}{section}

\newcommand\footnoteref[1]{\protected@xdef\@thefnmark{\ref{#1}}\@footnotemark}
\makeatother

\newcommand{\vp}{\varphi}
\newcommand{\dr}{\partial}
\DeclareMathOperator*{\osc}{osc}

\DeclareMathOperator{\diver}{div}

\DeclareMathOperator{\dist}{dist}
\DeclareMathOperator{\supp}{supp}
\DeclareMathOperator{\diam}{diam}

\DeclareMathOperator{\Tr}{Tr}

\def\div{\mathop{\operatorname{div}}}
\DeclareMathOperator{\loc}{loc}

\newcommand{\LforLp}{\rm{Loc}}
\newcommand{\lpformath}{\rm{(}}
\newcommand{\rpformath}{\rm{)}}
\newcommand{\PNformath}{\rm{PN}}
\newcommand{\Nformath}{\rm{N}}
\newcommand{\Dformath}{\rm{D}}
\newcommand{\Rformath}{\rm{R}}
\newcommand{\wPNformath}{\rm{wPN}}
\newcommand{\PRNformath}{\rm{PNR}}
\newcommand{\wPRNformath}{\rm{wPNR}}

\newcommand{\1}{{\mathds 1}}

\newcommand{\ms}{\medskip}

\newcommand{\R}{\mathbb R}
\newcommand{\bN}{\mathbb N}
\newcommand{\bZ}{\mathbb Z}

\newcommand{\bp}{\noindent {\em Proof: }}
\newcommand{\ep}{\hfill $\square$ \medskip}

\newcommand{\wt}{\widetilde}
\newcommand{\wh}{\widehat}

\newcommand{\A}{\mathcal A}
\newcommand{\N}{\mathcal N}
\newcommand{\C}{\mathcal C}

\newcommand{\Z}{\mathcal Z}

\newcommand{\NN}{\mathbb N}

\newcommand{\Rn}{\mathbb R^n}
\newcommand{\norm}[1]{\left\Vert#1\right\Vert}
\newcommand{\abs}[1]{\left\vert#1\right\vert}
\newcommand{\br}[1]{\left(#1\right)}
\newcommand{\set}[1]{\left\{#1\right\}}

\newcommand{\om}{\Omega}
\newcommand{\pom}{\partial\Omega}

\newcommand{\dint}{\int\!\!\!\!\!\int}
\def\Yint#1{\mathchoice
	{\YYint\displaystyle\textstyle{#1}}%
	{\YYint\textstyle\scriptstyle{#1}}%
	{\YYint\scriptstyle\scriptscriptstyle{#1}}%
	{\YYint\scriptscriptstyle\scriptscriptstyle{#1}}%
	\!\dint}
\def\YYint#1#2#3{{\setbox0=\hbox{$#1{#2#3}{\iint}$}
		\vcenter{\hbox{$#2#3$}}\kern-.51\wd0}}
\def\longdash{\mkern-1.5mu{-}\mkern-7.5mu{-}} 

\def\fiint{\Yint\longdash}

\theoremstyle{plain}
\newtheorem{theorem}[equation]{Theorem}
\newtheorem{lemma}[equation]{Lemma}
\newtheorem{corollary}[equation]{Corollary}
\newtheorem{proposition}[equation]{Proposition}

\newtheorem{definition}[equation]{Definition}

\theoremstyle{definition}

\theoremstyle{remark}

\newcommand{\Lp}{\hyperlink{Lp}{$\lpformath\LforLp_{p}\rpformath$}}
\newcommand{\Lwp}[1]{\hyperlink{Lp}{$\lpformath\LforLp_{#1}\rpformath$}}
\newcommand{\PNq}{\hyperlink{PNq}{$\lpformath\PNformath_{p'}\rpformath$}}

\newcommand{\Np}{\hyperlink{Np}{$\lpformath\Nformath_{p}\rpformath$}}
\newcommand{\Nwp}[1]{\hyperlink{Np}{$\lpformath\Nformath_{#1}\rpformath$}}
\newcommand{\wPNq}{\hyperlink{wPNq}{$\lpformath\wPNformath_{p'}\rpformath$}}
\newcommand{\wPNwq}[1]{\hyperlink{wPNq}{$\lpformath\wPNformath_{#1}\rpformath$}}
\newcommand{\PRNp}{\hyperlink{PRNp}{$\lpformath\PRNformath_{p}\rpformath$}}
\newcommand{\PRNwp}[1]{\hyperlink{PRNp}{$\lpformath\PRNformath_{#1}\rpformath$}}
\newcommand{\wPRNp}{\hyperlink{wPRNp}{$\lpformath\wPRNformath_{p}\rpformath$}}

\newcommand{\Dq}{\hyperlink{Dq}{$\lpformath\Dformath_{p'}\rpformath$}}

\newcommand{\Rp}{\hyperlink{Rp}{$\lpformath\Rformath_{p}\rpformath$}}
\newcommand{\Rwp}[1]{\hyperlink{Rp}{$\lpformath\Rformath_{#1}\rpformath$}}

\begin{document}

\title[The Poisson-Neumann problem]{The $L^p$ Poisson-Neumann problem and its relation to the Neumann problem}

\author[Feneuil]{Joseph Feneuil}
\address{Joseph Feneuil. Laboratoire de math\'ematiques d'Orsay, Universit\'e Paris-Saclay, France}
\email{joseph.feneuil@universite-paris-saclay.fr}

\author[Li]{Linhan Li}
\address{Linhan Li. School of Mathematics, The University of Edinburgh, Edinburgh, UK }
\email{linhan.li@ed.ac.uk}

\maketitle

\begin{abstract} 
We introduce the $L^p$ Poisson-Neumann problem for an uniformly elliptic operator $L=-\rm{div }A\nabla$ in divergence form in a bounded 1-sided Chord Arc Domain $\om$, which considers solutions to $Lu=h-\rm{div} \, \vec{F}$ in $\om$ with zero Neumann data on the boundary for $h$ and $\vec F$ in some tent spaces. We give different characterizations of solvability of the $L^p$ Poisson-Neumann problem and its weaker variants, and in particular, we show that solvability of the weak $L^p$ Poisson-Neumann probelm is equivalent to a weak reverse H\"older inequality. We show that the Poisson-Neumman problem is closely related to the $L^p$ Neumann problem, whose solvability is a long-standing open problem. We are able to improve the extrapolation of the $L^p$ Neumann problem from Kenig and Pipher \cite{KP93} by obtaining an extrapolation result on the Poisson-Neumann problem.
\end{abstract}

\ms\noindent{\bf Keywords: $L^p$ Neumann problem, Poisson-Neumann problem, Duality, Tent spaces}

\ms\noindent

\tableofcontents

\section{Introduction}

\subsection{State of the art}
Since the initial work of Dahlberg in \cite{Dah77}, there has been considerable interest in various boundary value problems for  elliptic operators in the form of $L=-\diver(A\nabla)$ with boundary data in $L^p$. The one that has seen the most achievements is the Dirichlet problem. In particular, there are many characterizations of solvability of the $L^p$ Dirichlet problem (see for instance \cite{JK82a,Ken94,KKPT00,DKP11,KKPT16,HL18,HLM19,CDMT22,BPTT??}), which leads to solvability of the Dirichlet problem for a large class of elliptic operators including those with non-smooth coefficients (\cite{JK81a, DJ90,KP01,DPP07,AAAHK11,HKMP15b,DP19,CHMT20,FP23, HLMP22}), and a clear picture of necessary and sufficient conditions on the domains (\cite{HMM16, GMT18,AHMMT20,HMMTZ21}) to ensure its solvability. 

In comparison to the Dirichlet problem, the Neumann problem with $L^p$ boundary data (see Definition \ref{def.Neu}) is much less understood.
Solvability of the $L^p$ Neumann problem has been established under much more restrictive conditions for both the operator and the domain compared to that of the Dirichlet problem (\cite{JK81b,KP93,KP95,MT03,KS08,KR08,AA11,AM14,DPR17,DHP23}). Roughly speaking, one cannot go too far from constant-coefficient operators on Lipschitz domains with small constant, or small perturbations from operators on  half space whose coefficient matrix is symmetric and independent of the transversal variable (``$t$-independent''), or one has to restrict to 2 dimension. In fact, it is a long standing open problem in the area to solve the $L^p$ Neumann problems in chord-arc and more general domains for the Laplacian, or for more general elliptic operators.   Lacking alternative formulations or characterizations of the Neumann problem is one of main obstacles in extending its solvability to a wider range of operators and domains. 

Sometimes viewed as a companion to the Neumann problem, the Regularity problem studies the relation between the gradient of the solution to the Dirichlet problem and the (tangential) derivatives of the boundary data (see Definition \ref{def.Reg}). In some special cases- when the dimension is 2, or when a ``Rellic identity'' holds
- the solvability of the $L^p$ Neumann problem can be deduced from the solvability of the $L^p$ Regularity problem. There have been big breakthroughs recently in solving the Regularity problem for elliptic operators which are (large) perturbations of constant-coefficient operators (\cite{MPT,DHP23}, see also \cite{Fen22}), building on earlier results of \cite{KP93, DPR17, MT??}. There are also results on the Regularity problem under different assumptions (\cite{HKMP15a,DK12,GMT??,DFMreg}).

Due to these advancements in the Dirichlet problem and the Regularity problem, it is reasonable to expect new developments in the Neumann problem. However, it is not clear at this stage how the Regularity problem (or the Dirichlet problem) would help solving the Neumann problem except for the special cases mentioned above. Moreover, the Neumann problem is kind of isolated from the Dirichlet problem and the Regularity problem, in the sense that there are counterexamples constructed in \cite{KP95} showing that solvability of the Regularity problem (or the Dirichlet problem) does NOT imply solvability of the Neumann problem, and vice versa. 

We decide to create a different world by introducing the $L^p$ Poisson-Neumann problem (see Definition \ref{def.PN}) and its (weaker) variants. As mentioned earlier, lacking equivalent characterizations of the Neumann problem and its isolation from  other boundary value problems are essentially the reasons that make solving the Neumann problem so hard. By considering the weak Poisson-Neumann problem for $L=-\diver(A\nabla)$ with bounded and measurable coefficients on 1-sided chord-arc domains (see Definition \ref{def.1CAD}), we are able to show that its solvability is equivalent to a (weak) reverse H\"older type inequality (Theorem \ref{MainTh}), which enjoys nice properties such as self-improvement, and has the potential to connect to many other conditions in rather general settings as we have seen in the $L^p$ Dirichlet problem.  
Not only the (weak) Poisson-Neumann problem is interesting by itself, our hope is that it can serve as a stepping stone to a better overall understanding of the Neumann problem. In fact, we can already improve some extrapolation result on the $L^p$ Neumann problem obtained by Kenig and Pipher (\cite{KP93}) with the help of the Poisson-Neumann problem (see Corollary \ref{cor.D+N=N}).

\medskip
Let us be more precise and start with giving the definitions of the various boundary value problems. We postpone the precise assumptions on the operator $L=-\diver A \nabla$ (which will however always have real coefficients in our article), the domain $\Omega$, as well as the definitions of the non-tangential maximal function $\N$, the square function $\A_1$, and the tent spaces $T^{p}_q(\Omega)$ to the preliminaries (Subsection \ref{SStent}). 

\begin{definition} \hypertarget{Dq}{}
If $p'\in (1,\infty)$, we say that the $L^{p'}$ Dirichlet problem \Dq\ - or \Dq$_L$ when we mention the operator - is solvable if there is $C>0$ such that for any $f\in C^\infty_0(\R^n)$, the solution $u_f$ to $Lu_f=0$ defined with the help of the elliptic measure by
\begin{equation} \label{defufhm}
u_f(X) := \int_{\partial \Omega} f(y) d\omega^X(y)
\end{equation}
satisfies
\[\|u_f\|_{T^{p'}_\infty(\Omega)} := \|\N(u_f)\|_{L^{p'}(\partial \Omega)} \leq C \|f\|_{L^{p'}(\partial \Omega)}.\]
\end{definition}

As mentioned at the beginning, the solvability of \Dq\ was studied in many scenarios, and it is not the purpose of this article to give a comprehensive survey on the literature.
Amongst the many results on the $L^p$ Dirichlet problem, we would like to bring out the following characterisations of \Dq\ by Mourgoglou, Poggi, and Tolsa (\cite{MPT}), which gives a well-rounded theory and generalizes some results from \cite{KP95}, and has inspired us into seeking their counterpart of the Neumann problem.

\begin{theorem}[\cite{MPT}] \label{ThMPT}
Let $\Omega \subset \R^n$ be 1-sided CAD\footnote{this is not optimal- see \cite{MPT} for a weaker assumption- but we shall limit ourselves to this setting for the present article. Moreover, Theorem \ref{ThMPT} is only a portion of the characterizations in \cite{MPT}.}, 
$L=-\diver A \nabla$ be a uniformly elliptic operator, and $p\in (1,\infty)$. The following are equivalent:
\begin{enumerate}
    \item \Dq$_{L^*}$ is solvable;
    \item the $L^{p'}$ Poisson-Dirichlet problem is solvable for $L^*$, meaning that there exists $C>0$ such that for any $\vec F \in L^\infty_c(\Omega,\R^n)$, the solution $u$ defined with the Green function as
    \[
        u(X) := \iint_\Omega \nabla_Y G(Y,X) \cdot \vec F(Y) \, dY 
    \]
    satisfies 
    \[ \|u\|_{\wt T^{p'}_\infty(\Omega)}:= \|\wt \N(u)\|_{L^{p'}(\partial \Omega)} \leq C \|\delta\vec F\|_{\wt T^{p'}_1(\Omega)}:= \|\wt \C_1(\delta \vec F)\|_{L^{p'}(\partial \Omega)},\]
    where $\delta(Z) = \dist(Z,\partial \Omega)$;
    \item the $L^p$ Poisson-regularity problem is solvable for $L$, meaning that there exists $C>0$ such that for any $h \in L^\infty_c(\Omega)$, the solution $u$ defined with the Green function as
    \[
        u(X) := \iint_\Omega G(X,Y) h(Y) \, dY 
    \]
    satisfies 
    \[\|\nabla u\|_{\wt T^{p}_\infty(\Omega)}:= \|\wt \N(\nabla u)\|_{L^{p}(\partial \Omega)} \leq C \|\delta h\|_{\wt T^{p}_1(\Omega)}:= \|\wt \C_1(\delta h)\|_{L^{p}(\partial \Omega)};\]
    \item There exists $C>0$ such that for any $Y\in \Omega$, the Green function satisfies the $L^p$ bound
    \begin{equation} \label{MPTbdG}
    \|\nabla G(.,Y) \1_{10B_Y \setminus B_Y}\|_{\wt T^{p}_\infty(\Omega)}:= \|\wt \N(\nabla G(.,Y) \1_{10B_Y \setminus B_Y})\|_{L^{p}(\partial \Omega)} \leq C \delta(Y)^{(1-n)/p'};
    \end{equation}
    where $B_Y := B(Y,\delta(Y)/4)$.
\end{enumerate}
\end{theorem}

\medskip
Next we present the definition of the $L^p$ Regularity problem.

\begin{definition}\label{def.Reg} \hypertarget{Rp}{}
Let $p\in (1,\infty)$. We say that the $L^p$ Regularity problem \Rp\ - or \Rp$_L$ - is solvable if there is $C>0$ such that for any $f\in C^\infty_0(\R^n)$, the solution $u_f$ to $Lu_f=0$ defined as in \eqref{defufhm} satisfies
\[\|\nabla u_f\|_{\wt T^{p}_\infty(\Omega)} := \|\wt\N(\nabla u_f)\|_{L^{p}(\partial \Omega)} \leq C \|\nabla f\|_{L^{p}(\partial \Omega)}.\]
\end{definition}

Note that $\|\nabla f\|_{L^{p}(\partial \Omega)}$ means that we have a notion of gradient on the boundary. We shall not spend to much time on it, because it is not the topic of this article. Briefly speaking, for rough boundaries (as in \cite{MT??,MPT}), we use the Haj\l asz upper gradient (see \cite{Haj96}), but when the boundary is Lipschitz, we can equivalently use the ``classical'' local tangential gradient. 
\medskip

There is some duality between the Dirichlet and the Regularity problem. One direction is given by the following result.

\begin{proposition}[\cite{MT??}, also \cite{KP93,DFMcarl}] \label{Rp=>Dq}
Let $\Omega \subset \R^n$ be a 1-sided CAD, $L=-\diver A \nabla$ be a uniformly elliptic operator, and $p\in (1,\infty)$. Then \Rp$_L$ $\implies$ \Dq$_{L^*}$.
\end{proposition}

The other direction is more delicate. It has been shown that \Dq$_{L^*}$ $\implies$ \Rp$_{L}$ in some cases such as for ``$t$-independent'' operators in Lipschitz graph domains (\cite{HKMP15a}), for operators that satisfy some Carleson measure condition 
(\cite{DPR17,DHP23,Fenreg}) in Lipschitz domains, and more generally, in corkscrew domains (\cite{MPT}). See also \cite{KS11reg, DK12,GMT??, AR12} for other settings.   

It is not yet known whether \Dq$_{L^*}$ $\implies$ \Rp$_{L}$ holds in general. However, there is a partial result by Shen (\cite{Shen07}), who showed that \Dq\  implies a dichotomy in the solvability of the Regularity problem:  
if  \Dq\  is solvable, then either \Rp\ is solvable or \Rwp{q} is not solvable for any $1 < q < \infty$.

\begin{proposition}[\cite{Shen07}] \label{Rq+Dp'=>Rp}
Let $\Omega \subset \R^n$ be a Lipschitz domain, and let $L=-\diver A \nabla$ be a uniformly elliptic operator\footnote{\cite{Shen07} is written for Lipschitz domains and symmetric coefficient matrix $A$. The authors guarantee the generalization to non-symmetric $A$ with the same proof; going beyond Lipschitz domain requires using Haj\l asz gradient and we are not aware whether it has been checked.}, and $p,q\in (1,\infty)$. Then 
\[\text{\Rwp{q}$_L$ $+$ \Dq$_{L^*}$ $\implies$ \Rp$_{L}$.}\]
\end{proposition}

\medskip

Let us now give the definition of the $L^p$ Neumann problem.

\begin{definition}\label{def.Neu} \hypertarget{Np}{}
If $p\in (1,\infty)$, we say that the $L^p$ Neumann problem \Np\ - or \Np$_L$ - is solvable if there is $C>0$ such that for any $f\in C^\infty_0(\R^n)$ satisfying $\int_{\partial \Omega} f\, d\sigma = 0$, the solution $u_f$ to $Lu_f=0$ defined with the help of the Neumann function as 
\begin{equation} \label{defufNp}
u_f(X) := \int_{\partial \Omega} N(X,y) f(y) \, d\sigma(y)
\end{equation}
satisfies
\[\|\nabla u_f\|_{\wt T^{p}_\infty(\Omega)} := \|\wt\N(\nabla u_f)\|_{L^{p}(\partial \Omega)} \leq C \|f\|_{L^{p}(\partial \Omega)}.\]
\end{definition}

It is not clear what the Neumann function is when the domain is too rough, as it was only constructed and studied in a ball (\cite{KP93}), but many of its properties are still valid in more general domains such as bounded 1-sided CAD. We shall give the definition and present the properties that we need in the preliminaries below. While writing this paper, we learned that the Neumann function in 1-sided CAD is constructed in a work in preparation \cite{Steve'student}), along with properties beyond what we need for this paper. 

Let us give a brief discussion on the history of the $L^p$ Neumann problem, and again, we are not trying to be exhaustive. The solvability of the $L^p$ Neumann problem  for the Laplacian in Lipschitz domains is established for $p\in (1,2+\epsilon)$ in \cite{DK87}, and it is optimal in the sense that for any $q>2$, we can find a Lipschitz domain for which \Np\ is not solvable.
Their result is based on the solvability of $L^2$ Neumann problem for the Laplacian in Lipschitz domains, which is proven in \cite{JK81b} via   establishing a ``Rellich identity'', that is, $\norm{\nabla_{\rm tan}u}_{L^2(\pom)}\approx\norm{\frac{\dr u}{\dr n}}_{L^2(\pom)}$. 
The ``Rellich identity" can be generalized to elliptic operators with real, symmetric coefficient matrix of ``$t$-independent" coefficients on half space, which is essentially proved in \cite{KP93} to obtain solvability of the $L^2$ Neumann problem in that setting. 
Solvability of the Neumann problem for non-symmetric $t$-independent operators in the half space is still an open problem, but many works pushed the limits beyond the symmetric $t$-independent condition (a Dini condition on the $t$-dependence \cite{MT03}, in the half plane \cite{KR08}, small complex $t$-independent perturbations \cite{AA11,AAAHK11}). Elliptic operators involving Carleson measures are studied in \cite{KP93,KP95,DPR17}, but it is still an open problem whether the $L^p$ Neumann problem is sovlable without the smallness assumption on the Carleson norm of the coefficients, except {\em a posteriori} in $\R^2$, where we can link the Neumann problem to the regularity problem with a trick from \cite{KR08} (see \cite{DHP23}).

\medskip
It is useful to know whether solving \Np\ from some $p\in(1,\infty)$ gives the solvability of the $L^q$ Neumann problem for any other $q$. In \cite[Theorem 6.3]{KP93}, it was shown that on balls,
\begin{equation} \label{Rp+Np=>Nq}
    \text{ \Rp$_L$ $+$ \Np$_L$ $\implies$ \Nwp{q}$_L$\ \ for $q\in (1,p+\epsilon)$.}
\end{equation}
for some $\epsilon>0$. When writing this paper, we learned that this result is generalized to CAD domains in \cite{Steve'student}. Since there are examples of \Np $\centernot\implies$ \Rp\  even for symmetric operators in nice domains (see \cite{KP93}), one cannot extend the range of solvabiltiy of the Neumann problem without knowing the solvability of the Regularity problem in general. One of our contributions is that we can improve this result by replacing \Rp$_L$ with \Dq$_{L^*}$, and our result holds for 1-sided CAD.

\subsection{Our results}
Throughout the article, $n\geq 2$ and $\Omega \subset \R^n$ is a bounded 1-sided Chord Arc Domain (see Subsection \ref{SSCAD}), $L=-\diver A\nabla$ is an uniformly elliptic (see \eqref{defelliptic}) operator with bounded and measurable coefficients on $\om$, and $L^*=-\diver A^T\nabla$ is the adjoint operator.  The limitation to bounded 1-sided Chord Arc Domain is only due to the lack of construction of a suitable elliptic theory around the Neumann problem in more general domains. Should the Neumann function and its properties be studied in more general settings (like in the mixed codimension setting \cite{DFMhighcod,DFMmixed} or in domains that lack connectivity \cite{MT??,MPT,CHM??}), our assumption would certainly be weaker.

\medskip

Let us introduce our additional boundary value problems. First is the Poisson-Neumann problem (we recall that the definitions of tent spaces, $\wt \N$ or $\wt \C_1$ are delayed to Subsection \ref{SStent}).

\begin{definition}\label{def.PN} \hypertarget{PNq}{}
Let $p\in (1,\infty)$. We say that the $L^{p'}$ Poisson-Neumann problem \PNq\ is solvable for $L^*$ if for any $\vec{F}\in L^\infty_c(\Omega,\R^n)$, the solution $u$ to $L^*u=-\diver \vec F$ with zero average - constructed for instance via the Neumann function as
\begin{equation} \label{eq.solwPN}
\bar u(X) := \iint_\Omega \nabla_Y N(Y,X) \cdot \vec F(Y) \, dY
\end{equation}
 and then $u := \bar u - u_{\Omega}$ - satisfies
\[\|u\|_{\wt T^{p'}_\infty}:= \|\wt \N(u)\|_{L^{p'}(\partial \Omega)} \leq C \|\delta \vec F\|_{\wt T^{p'}_1} := C  \|\wt\C_1(\delta\vec{F})\|_{L^{p'}(\partial \Omega)},\]
where $C$ is independent of $\vec F$. 

We say that the weak $L^{p'}$ Poisson-Neumann problem \wPNq\ is solvable for $L^*$ if for any $\vec{F}\in L^\infty_c(\Omega,\R^n)$, the solution $u$ constructed as in \eqref{eq.solwPN} satisfies
\[\|\wt \N(\delta \nabla u)\|_{L^{p'}(\partial \Omega)} \leq C  \|\wt\C_1(\delta\vec{F})\|_{L^{p'}(\partial \Omega)},\]
where $C$ is independent of $\vec F$. 

We write \PNq$_{L^*}$ or \wPNq$_{L^*}$ when we want to mention the operator.
\end{definition}

The next definition presents the ``Poisson-Neumann-Regularity problem'', that is the dual formulation of the Poisson-Neumann problem. We shall introduce the space $\widehat L^\infty_c(\Omega)$ of compactly supported and bounded measurable functions with zero average.

\begin{definition} \hypertarget{PRNp}{}
Let $p\in (1,\infty)$. We say that the $L^p$ Poisson-Neumann-Regularity \PRNp\ is solvable for $L$ if for any $h\in \widehat{L}^\infty_c(\Omega)$ and $\vec{F}\in L^\infty_c(\Omega,\R^n)$, the solution $u$ to $Lu=h-\diver \vec F$ constructed with the Neumann function as 
\begin{equation} \label{eq.solwPRN}
u(X):= \iint_\Omega \big[\nabla_Y N(X,Y) \cdot \vec F(Y) + N(X,Y) h(Y)\big] \, dY
\end{equation}
satisfies
\[\|\wt\N(\nabla u) \|_{L^{p}(\partial \Omega)} \leq C \|\wt\C_1(\delta |h| + |\vec{F}|)\|_{L^{p}(\partial \Omega)},\]
where $C$ is independent of $h$ and $\vec F$. 

We say that the property \wPRNp\ holds for $L$ if for any $\vec{F}\in L^\infty_c(\Omega,\R^n)$, the solution constructed in \eqref{eq.solwPRN} with $h\equiv 0$ satisfies
\[\|\wt\N(\nabla u) \|_{L^{p}(\partial \Omega)} \leq C  \|\wt\C_1(\vec{F})\|_{L^{p}(\partial \Omega)},\]
where $C$ is independent of $\vec F$.

We write \PRNp$_L$ or \wPRNp$_L$ when we want to mention the operator.
\end{definition}

Finally, our last property is the localization of solutions. It is meant to be a substitute to the property (4) of Theorem \ref{ThMPT}. Although it is stronger than the latter, the localization property is more flexible to use and its Dirichlet counterpart is a key component of the proof of Proposition \ref{Rq+Dp'=>Rp}.   

\begin{definition} \hypertarget{Lp}{} 
Let $\Omega \subset \R^n$ be a bounded 1-sided CAD, $L$ be a uniformly elliptic operator and $p \in [1,\infty)$. We say that the property \Lp\ holds for $L$ if for any ball $B=B(x,r)$ centered at the boundary, and any local solution $u \in W^{1,2}(2B\cap \Omega)$ to $Lu=0$ in $2B\cap\om$  with zero Neumann data on $2B\cap\pom$, we have
\begin{equation} \label{loc1}
\|\wt \N(|\nabla u|\1_B)\|_{L^p(\partial \Omega)} \leq C r^{(n-1)/p}\left ( \fiint_{2B\cap \om} |\nabla u|^2 dX \right)^\frac12.
\end{equation}
\end{definition}

A similar localization property was used in \cite{KS08}: it was proven there that for the Laplacian in bounded Lipschitz domains, \Np\ is equivalent to a weak reverse H\"older inequality for $p>2$, from which \Np\ is derived for the Laplacian on bounded convex domains for a range of $p>1$. 

\medskip

With all these properties in hand, we can start the presentation of our results. If one views the Poisson-Neumann and Poisson-Neumann-regularity problem as analogues of the Poisson-Dirichlet and Poisson-regularity problem, respectively, with zero Neumann data instead of zero Dirichlet data, then our first result is a natural analogue of  (2) $\iff$ (3) in Theorem \ref{ThMPT}.

\begin{proposition}
Let $\Omega \subset \R^n$ be a bounded 1-sided CAD, $L= - \diver A \nabla$ be a uniformly elliptic operator, and $p\in (1,\infty)$. The following are equivalent,
\begin{enumerate}[(i)]
\item The Poisson-Neumann problem \PNq$_{L^*}$ is solvable;
\item the Poisson-Neumann-Regularity \PRNp$_L$ is solvable;
\item the Poisson-Neumann-Regularity \PRNp$_L$ is solvable for $\vec F =0$.
\end{enumerate}
If one of the above condition is satisfied, we say that the strong $L^p$ Neumann problem is solvable for $L$. 
\end{proposition}

We use the name ``strong $L^p$ Neumann problem'' because the solvability of the $L^p$ Neumann problem \Np$_L$ is actually implied by \PNq$_{L^*}$ or \PRNp$_{L}$. We do not know whether the strong $L^p$ Neumann problem is equivalent to the $L^p$ Neumann problem, but we have a partial converse by assuming in addition \Dq$_{L^*}$. The following proposition is a summary of these observations.

\begin{proposition} \label{Np<=>PRNp} 
Let $\Omega \subset \R^n$ be a bounded 1-sided CAD, $L= - \diver A \nabla$ be a uniformly elliptic operator, and $p\in (1,\infty)$. Then 
\[ \text{\Np$_L$ $+$ \Dq$_{L^*}$ $\implies$ \PRNp$_L$ $\implies$ \Np$_L$.}\]
\end{proposition}

We came to the conclusion that the strong $L^p$ Neumann problem should be put in parallel to the Regularity problem \Rp$_L$, and that a bound on the Neumann function analoguous to \eqref{MPTbdG} would be a characterization of a ``missing'' problem should make the parrallel with the Dirichlet problem \Dq$_{L^*}$. We looked at how we should weaken \PNq$_{L^*}$ and \PRNp$_L$ so that it is more or less equivalent to a local $T^{p}_\infty(\Omega)$-bound on the (gradient of the) Neumann function, and we arrived at the following characterizations.

\begin{theorem} \label{MainTh}
    Let $\Omega \subset \R^n$ be a bounded 1-sided CAD, $L=-\div A \nabla$ be a uniformly elliptic operator on $\Omega$, and $p\in (1,\infty)$. The following are equivalent.
      \begin{enumerate}[(a)]
        \item \wPNq$_{L^*}$ is solvable;
        \item \wPRNp$_{L}$ is solvable;
        \item \Lp$_L$ holds.
    \end{enumerate}
    If one of the above conditions is true, we say that the weak $L^{p'}$ Poisson-Neumann problem is solvable for $L^*$, in reference to the property \wPNq$_{L^*}$. 
\end{theorem}

One can see now that the weak Poisson-Neumann is indeed weaker than the Poisson-Neumann problem, as \PRNp$_L$ implies \wPRNp$_L$ by definition. 

The localization property can be viewed as a weak reverse H\"older inequality, which is a nice property that we would like to have for the weak Poisson-Neumann problem. For one thing, we can show that the property \Lp$_L$ self-improves (see Lemma \ref{Lp=>Lp+e}). We remark that the localization property \Lp$_L$ is defined as it is instead of a bound on the Neumann function because we didn't succeed to self-improve the bound on the Neumann function from (local) $\wt T^{p}_\infty(\Omega)$ to (local) $\wt T^{p+\epsilon}_\infty(\Omega)$. 

Another advantage of having the localization property characterization of the weak Poisson-Neumann problem is that we can use it to prove a partial converse of \wPNq$_{L^*}$ $\implies$ \PRNp$_{L}$ in the spirit of Proposition \ref{Rq+Dp'=>Rp}.

\begin{theorem} \label{MainTh2}
    Let $\Omega \subset \R^n$ be a bounded 1-sided CAD, $L=-\div A \nabla$ be a uniformly elliptic operator on $\Omega$, and $p,q\in (1,\infty)$. Then 
    \[ \text{ \wPNq$_{L^*}$ $+$ \PRNwp{q}$_L$ $\implies$ \PRNp$_L$.}\]
\end{theorem}

As a consequence of all our results, we are able to show that the weak Poisson-Neumann problem and the strong Neumann problem are solvable on an open interval.

\begin{corollary}\label{cor.main}
    Let $\Omega \subset \R^n$ be a bounded 1-sided CAD, $L=-\div A \nabla$ be a uniformly elliptic operator on $\Omega$, and $p\in (1,\infty)$. 
    \begin{enumerate}
        \item If \wPNq$_{L^*}$ is solvable, then there exists $\epsilon>0$ such that \wPNwq{q'}$_{L^*}$ is solvable for any $q'\in (p'-\epsilon,\infty)$;
        \item If \PRNp$_{L}$ is solvable, then there exists $\epsilon>0$ such that \PRNwp{q}$_{L}$ is solvable for any $q\in (1,p+\epsilon)$.
    \end{enumerate}
\end{corollary}

Note that part (2) of Corollary \ref{cor.main} together with Proposition \ref{Np<=>PRNp} improves the earlier result of Kenig and Pipher \eqref{Rp+Np=>Nq}, as we now have 
\begin{corollary}\label{cor.D+N=N}
Let $\Omega \subset \R^n$ be a bounded 1-sided CAD, $L=-\div A \nabla$ be a uniformly elliptic operator on $\Omega$, and $p\in (1,\infty)$. 
    \begin{equation} \label{Dq+Nq=>Nq}
\text{\Dq$_{L^*}$ $+$ \Np$_L$ $\implies$ \Nwp{q}$_L$ \ for $q\in (1,p+\epsilon)$.}
\end{equation}
\end{corollary}

\medskip

\subsection{Plan of the article, and how to gather the results to get the main results.}

We summarize our results in the following diagram (dashed double line means that the implications are true under additional assumptions).

 
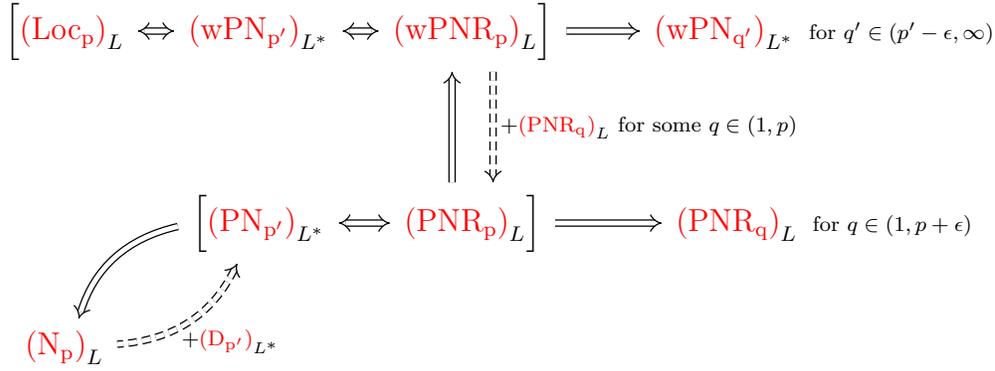
\begin{tikzcd}[column sep=small]
\Big[ \text{\Lp}_L \arrow[r,Leftrightarrow] & 
\text{\wPNq}_{L^*} \arrow[r,Leftrightarrow]  & 
\text{\wPRNp}_L \Big] \arrow[rr, Rightarrow] \arrow[dd,shift right=1.5ex,Leftarrow] \arrow[dd,Rightarrow,shift left=1.5ex, dashed,"+\text{\PRNwp{q}}_L \text{ for some $q\in (1,p)$}"] & & 
\text{\wPNwq{q'}}_{L^*} \, \text{ \tiny for $q'\in (p'-\epsilon,\infty)$} \\
\\
& \Big[ \text{\PNq}_{L^*} \arrow[r,Leftrightarrow] \arrow[dl,Leftarrow, dashed, bend left, "\hspace{1.1cm}+ \text{\Dq}_{L^*}" {below}]& 
\text{\PRNp}_L \Big] \arrow[rr,Rightarrow]& & 
\text{\PRNwp{q}}_L \, \text{ \tiny for $q\in (1,p+\epsilon)$} \\
\text{\Np}_L \arrow[ur,Leftarrow,bend left] 
\end{tikzcd}

\medskip

We let the reader check that those results are a consequence of the following implications\footnote{Single dashed arrow means that we also use the property at the tail for the implication at the tip}.
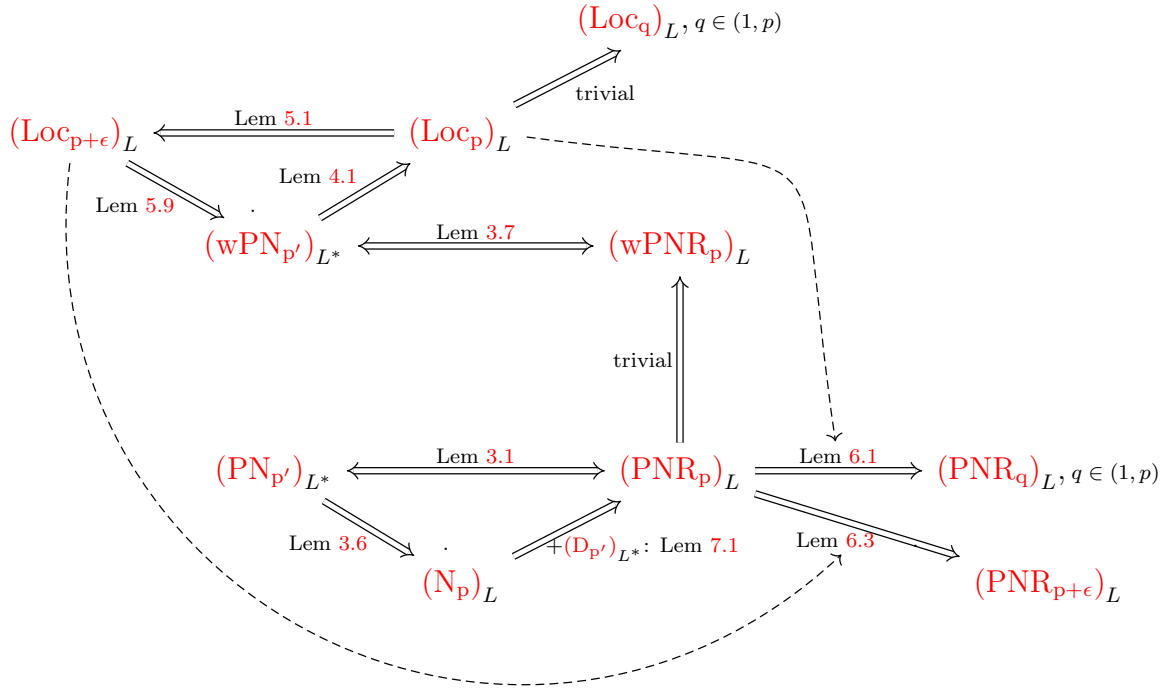
\smallskip

 \begin{tikzcd}[column sep=small]
 & & & & \text{\Lwp{q}}_{L}, \text{\tiny $q\in (1,p)$} \arrow[dl,Leftarrow, "\text{trivial}"]& \\
 &  \text{\Lwp{p+\epsilon}}_{L}  \arrow[dr,Rightarrow, "\text{Lem \ref{lemBp=>PNp'}} \hspace{1cm}." below] 
 & & \text{\Lp}_{L} \arrow[ll,Rightarrow, "\text{Lem \ref{Lp=>Lp+e}}" above]   &  &   \\
 & & \text{\wPNq}_{L^*} \arrow[rr,Leftrightarrow, "\text{Lem \ref{wPNq=>wPRNp}}"] \arrow[ur,Rightarrow, "\text{Lem \ref{wPNq=>Lp}}"]
 & & \text{\wPRNp}_{L}\\
 & & & & & & \\
 & & & & & & \\
 & & \text{\PNq}_{L^*} \arrow[rr,Leftrightarrow, "\text{Lem \ref{PNq<=>PRNp}}"] \arrow[dr,Rightarrow, "\text{Lem \ref{PNq=>Np}} \hspace{1cm}." below] 
 & & \text{\PRNp}_{L} \arrow[uuu,Rightarrow, "\text{trivial}"] \arrow[dl,Leftarrow, "\hspace{2cm}+\text{\Dq}_{L^*}: \ \text{Lem \ref{Np+Dq=>PRNp}}" below] 
 & \phantom{U} & \text{\PRNwp{q}}_{L}, \text{\tiny  $q\in (1,p)$} \arrow[ll,Leftarrow,"\text{Lem \ref{PRNp=>PRNq}}"{name=Y,above}]\\
 & & & \text{\Np}_L & & \phantom{U} & \text{\PRNwp{p+\epsilon}}_L \arrow[ull,Leftarrow,"\text{Lem \ref{lemShen}\hspace{0.5cm}.}"{name=Z,below}] \\
 \arrow[from=2-2, to=Z, dashed,bend right=70,looseness=1.3] 
 \arrow[from=2-4, to=Y, dashed,bend left=37,looseness=2]
\end{tikzcd}

Actually, we are missing one thing in the above diagram to complete the first diagram: our Lemma \ref{lemShen} also gives that for {\it any} $q\in (1,p)$,
\[\text{\Lwp{p+\epsilon}}_L + \text{\PRNwp{q}}_L \implies \text{\PRNp}_L.\]
In particular, when $q=p-\epsilon$, it implies  \PRNp$_L$ $\implies$ \PRNwp{p+\epsilon}$_L$.
\smallskip

The plan of our article is as follows. Section 2 is devoted to the preliminaries, where we first present our geometric and topological condition (Subsection \ref{SSCAD}), then we introduce the functional analysis (the  non-tangential function $\wt \N$, the square function $\wt\A_1$, the Carleson functional $\wt \C_1$ and their relation to tent space), and we finish by stating results related to the elliptic theory and the Neumann function. In Section \ref{S3}, we give results related to duality; Section \ref{S4} deals with the localization property; Section \ref{S5} shows how to go from the localization property to solvabiltiy of the weak Poisson Neumann; Section \ref{S6} is about the interpolation and extrapolation of the Poisson Neumann regularity problem. Finally, Section \ref{S7} establishes how to recover the Poisson Neumann problem from the Neumann problem. The Appendix proves an interpolation result on tents spaces.

\subsection{Acknowledgement} The authors want to thank Pascal Auscher for pointing out references on tents spaces, in particular \cite{Huang17}, and Guy David and Marco Michetti, for an early presentation of their work \cite{DDEMM??}.

\section{Preliminaries}
Throughout this article, we always assume that $\om\subset\Rn$ is an open and connected set (i.e. a domain), and for $X\in\om$, we use the notations $\delta(X):= \dist(X,\partial \Omega$, $B_X:=B(X,\delta(X)/4)$.
\subsection{Our conditions} \label{SSCAD} 

\begin{definition}
Let $E \subset \R^n$ be a closed set. We say that $E$ is $(n-1)$-Ahlfors regular if there exists a measure $\sigma$ - called Ahlfors regular measure - and $C=C_\sigma>0$ such that for any $x\in E$ and any $r \in (0,\diam E)$, we have
\[C {r}^{n-1} \leq \sigma(B(x,r)) \leq Cr^{n-1}.\]
Note that if $E$ is $(n-1)$-Ahlfors regular, $\mathcal H^{n-1}|_E$ is an Ahlfors regular measure on $E$. 
\end{definition}

Before we continue, let us introduce a Whitney decomposition of our open domain $\Omega$. Let $\mathbb D_{\R^n}:= \bigcup_{k\in \bZ} \mathbb D_{\R^n,k}:= \bigcup_{k\in \bZ} \{\prod_{i=1}^n [m_i2^{-k},(m_i+1)2^{-k})\}_{(m_1,\dots,m_n)\in \bZ^{n}}$ be the collection of dyadic cubes of $\R^n$. If $Q \in \mathbb D_{\R^n}$, we write $k(Q)$ for the value - called {\em generation} of $Q$ - such that $Q \in \mathbb D_{\R^n,k}$, $\ell(Q)$ for $2^{-k(Q)}$, $r(Q)$ for $\diam(Q)$, and $\lambda Q$ for $\{x\in \R^n, \, \dist(x,Q) \leq (\lambda-1)r(Q)\}$. Notice that the collection $\mathbb D_{\R^n}$ has a natural order (inclusion).

We construct $\mathcal W_\Omega$ - called Whitney decomposition - as the subcollection of the dyadic cubes $W$ that are maximal for the property $10W \subset \Omega$. We say that $W,W'$ are adjacent if 
\[\rm{Interior}(W)\cup \rm{Interior}(W')\subsetneqq{\rm Interior}\br{\overline{W\cup W'}}.\]
Observe that $\mathcal W_\Omega$ is a partition of $\Omega$, and that two adjacent Whitney cubes have generations $k(W),k(W')$ that differ from at most 1.

\begin{definition}[1-sided NTA domains and Harnack chains] \label{defHC}
    We say that $Y\in \Omega$ and $Z\in \Omega$ are linked by an {\bf $\kappa$-Harnack chain of Whitney cubes} if there exists a chain $W_1, \dots, W_N$ of Whitney cubes $W_i \in \mathcal W_\Omega$ such that 
    \begin{enumerate}
        \item $Y \in \overline{W_1}$ and $Z \in \overline{W_N}$;
        \item for each $i< N$, $W_i$ and $W_{i+1}$ are adjacent;
        \item for any $i\in \{1,\dots,N\}$ we have
        \[ r(W_i) \geq \kappa \min\left\{ \sum_{j=1}^i r(W_j) , \sum_{j=i}^N r(W_j) \right\}\]
        \item $\displaystyle |Y-Z| \geq \kappa \sum_{j=1}^N r(W_j)$.
    \end{enumerate}
    We use the notation
    \[H_\Omega(Y,Z) = {\rm Interior} \left( \bigcup_{i=1}^N \overline{W_i} \right)\]
    for the union of the Whitney cubes in the Harnack chain linking $Y$ and $Z$.

    We say that $\Omega$ is a {\bf John domain} if $\Omega$ is bounded and there exists a center $X_0$\footnote{By abuse of notation, we say that a ball $B\subset \Omega$ is the John center of $\Omega$ if $B=B_{X_0}=B(X_0,\delta(X_0)/4)$} and a constant $\kappa>0$ such that $\delta(X_0)  \geq \kappa \diam(\Omega)$ and for each $Y\in \Omega$, $Y$ and $X_0$ are linked by a $\kappa$-Harnack chain of Whitney cubes.

   Finally, we say that $\Omega$ is a {\bf 1-sided NTA domain} if there exists $\kappa>0$ such that for any $Y,Z\in \Omega$, $Y$ and $Z$ are linked by a $\kappa$-Harnack chain of Whitney cubes.

\end{definition}

The above definition is a unclassical, as it relies on the Whitney cubes, but we let the reader check as an easy exercise that our definition is equivalent to the usual definitions of John domains and 1-sided NTA domains that are familiar to the reader. Another observation that we want to raise up is the fact that, if $\Omega$ is 1-sided NTA and $E \subset \Omega$, then
\begin{equation} \label{constructionS} 
S := \bigcup_{Y,Z \in E} H_\Omega(Y,Z)
\end{equation}
is a John domain, which allows a simple construction of $1$-sided NTA subdomains.

\medskip

John domains are the domains that are well adapted to Poincar\'e inequalities. 

\begin{proposition}[Poincar\'e inequality]  \label{prop.Poincare}
If $S$ is a John domain, then for any $w\in W^{1,2}(S)$, we have
    \[\left( \iint_S |w-w_S|^2 \, dX \right)^\frac12  \leq C \diam(S) \left( \iint_S |\nabla w|^2 \, dX \right)^\frac12,\]
    where $w_S$ is either the average $\fiint_S w$ or the average on the John center $\fiint_{B_{X_0}}$, and where the constant $C>0$ depends only on the John constant.

    If particular, if $\Omega$ is 1-sided NTA, then there exists $C,K\geq 1$ (depending only on the 1-sided NTA constant of $\Omega$) such that for any $B=B(x,r)$ centered on $\partial \Omega$ and any $w\in W^{1,2}(KB\cap \Omega)$,
    \[ \left( \iint_{B \cap \Omega} |w-w_B|^2 \, dX \right)^\frac12  \leq C r \left( \iint_{KB \cap \Omega} |\nabla w|^2 \, dX \right)^\frac12.\]
\end{proposition}

\bp
The first part is a consequence of \cite[Theorem 1]{HK95}, since $L^2$-Poincar\'e inequality on balls are known to be true. For the second part, construct the John domain $B \cap \Omega \subset S$ as in \eqref{constructionS}, and notice that $S \subset KB \cap \Omega$ for some $K$ independent of $B$, then use the first part. We can even have construct $S$ to be 1-sided NTA, see \cite[Appendix A]{HM14}.  
\ep

We conclude our paragraph with the definition of our main geometric setting.

\begin{definition}\label{def.1CAD}
    We say that $\Omega \subset \Omega$ is a 1-sided Chord Arc Domain (1-sided CAD for short) if $\partial \Omega$ is $(n-1)$-Ahlfors regular and $\Omega$ is bounded whenever $\partial \Omega$ is bounded and 1-sided NTA. Moreover, $\sigma$ will always denote the Ahlfors regular measure on the previously introduced 1-sided CAD. 
\end{definition}

\subsection{Tents spaces} \label{SStent}
Before we introduce all our properties, we require to (recall or) give a bit more definitions. Let $\Omega \subset \R^n$ be a 1-sided CAD, and $\sigma$ be its Ahlfors regular measure. For $X\in\om$, recall that $\delta(X):=\dist(X,\pom)$ and $B_X := B(X,\delta(X)/4)$.

For $a>1$, define the cones $\gamma_a(x)$ with vertex in $x$ and aperture $a$ as
\[\gamma_{a}(x):=\set{X\in\om: \abs{X-x}<a\dist(X,\pom)}.\]
This is used to define the modified non-tangential maximal function and the area function of $F\in L_{\loc}^2(\om)$. If $c\in(0,1)$, we have
\[
\wt\N^{a,c}(u)(x):= \sup_{X\in \gamma_a(x)}\br{\fiint_{B(X,c\delta(X))}\abs{u(Y)}^2dY}^{1/2}, \quad x\in\pom,
\]
and
\begin{equation}\label{defeq.wtArea}
\wt\A_1^{a,c}(F)(x):=\iint_{\gamma_a(x)} \br{\fiint_{B(X,c\delta(X))} |F|^2}^{1/2} \frac{dX}{\delta(X)^n},  
\end{equation}
while the averaged Carleson functional of a function $F\in L_{\loc}^2(\om)$ is
\begin{equation}\label{defeq.CavgF}
    \wt\C_1^{c}(F)(x) := \sup_{r>0} \frac{1}{\sigma(B(x,r))}\iint_{B(x,r) \cap \om} \br{\fiint_{B(X,c\delta(X))} |F|^2}^{1/2} \, \frac{dX}{\delta(X)}, \quad \text{for }x\in \partial \Omega.
\end{equation}
Moreover, we set
\[\gamma(x):= \gamma_2(x), \qquad \wt \N := \wt \N^{2,1/4}, \qquad \wt \A_1 := \wt \A^{2,1/4}_1  \quad \text{ and } \quad  \wt \C_1 := \wt \C^{1/4}_1 \]
In the above definitions, we use $\sim$ over the functionals $\N$, $\A_1$ and $\C_1$ to indicate that we are taking an average on Whithney balls. We keep the subscript $1$ in the definitions of $\wt \C$ and $\wt \A$, because we want to emphasize that we are employing the functionals related to the tent spaces $T^{p}_1(\Omega)$ (the most common occurence of the area functional $\A$ in PDE is the ``square functional'' - i.e. the $L^2$-version of $\A_1$) and not $T^p_2(\Omega)$.

\begin{lemma} \label{lemNCA}
Let $\Omega \subset \R^n$ be 1-sided CAD, $a >1$, $c\in (0,1)$, and $p\in (0,\infty]$. Then 
\[\|\wt \N^{a,c}(u)\|_{L^p(\partial \Omega,\sigma)} \approx \|\wt \N(u)\|_{L^p(\partial \Omega,\sigma)},\]
\[\|\wt \A_1^{a,c}(u)\|_{L^p(\partial \Omega,\sigma)} \approx \|\wt \A_1(u)\|_{L^p(\partial \Omega,\sigma)},\]
and
\[\|\wt \C_1^{c}(u)\|_{L^p(\partial \Omega,\sigma)} \approx \|\wt \C_1(u)\|_{L^p(\partial \Omega,\sigma)},\]
where the constants involved depend on $\Omega$, $a$, $c$ and $p$. Moreover, if $p\in (1,\infty)$, then
\[\|\wt \C_1(u)\|_{L^p(\partial \Omega,\sigma)} \approx \|\wt \A_1(u)\|_{L^p(\partial \Omega,\sigma)},\]
where the constant depend on $\Omega$ and $p$.
\end{lemma}

\bp
The fact that we can change the aperture is a well known fact. The proof is written in \cite[Chapter II, 2.5.1]{Stein93} and \cite[Proposition 4]{CMS85} when $\Omega = \R^n_+$ but the proof can be immediately adapted to 1-sided CAD.

The fact that we can change $c$ can be proven by hand, since increasing $c$ can be compensated by increasing $a$ instead (see for instance \cite{MPT} Lemma 2.1 when $p>1$, but the proof result easily extend to $p\in (0,1]$).

Finally, the last equivalence between the square functional $\wt \A_1$ and the Carleson functional $\wt \C_1$ is again common knowledge to people working on tent spaces (see for instance \cite[Theorem 3]{CMS85}).
\ep

With those in hand, the space of functions $u$ for which $\|\wt \N^{a,c}(u)\|_{p}$ or $\|\wt \A_1^{c}(u)\|_{p}$ are finite doesn't depend on $a$ or $c$. We call them ($L^2$ averaged) tent spaces $\wt T^p_\infty(\Omega)$ and $\wt T^p_1(\Omega)$ whose norm are
\[\|u\|_{\wt T^p_\infty(\Omega)} := \|\wt \N^{a,c}(u)\|_{L^p(\partial \Omega,\sigma)},\]
\[\|u\|_{\wt T^p_1(\Omega)} := \|\wt \A_1^{a,c}(u)\|_{L^p(\partial \Omega,\sigma)}\qquad\text{for }p\in(0,\infty),\]
\[\|u\|_{\wt T^\infty_1(\Omega)} := \|\wt \C_1^{a,c}(u)\|_{L^\infty(\partial \Omega,\sigma)}.\]
Those averaged tent spaces are a common occurence in PDE: the non-tangential maximal function and the square functional are basic tools to estimates the solutions of boundary value problems with $L^p$ data. 

In our article, we shall use the notation $\wt \N^*$, $\wt \C^*_1$, or $\wt \A^*_1$ when we want to enlarge the aperture $a$ or the radius constant $c$ (instead of explicitly write $a$ and $c$). And we can come back to the $\wt \N$, $\wt \C_1$, or $\wt \A_1$ by invoking the following lemma.

\medskip

We shall frequently use the following duality between the Carleson functional $\wt\C_1$ and the non-tangential maximal function $\wt\N$.

\begin{lemma} If $\Omega$ is 1-sided CAD and $p\in [1,\infty)$, then for any $u\in \wt T^p_\infty(\Omega)$ and $F \in \wt T^{p'}_1(\Omega)$,
\begin{equation}\label{eq.CNdual}
    \abs{\iint_{\om}uF\,dX}\lesssim \norm{\wt\C_1(\delta F)}_{L^{p'}(\pom)}\norm{\wt\N(u)}_{L^p(\pom)}.
\end{equation}
Moreover
\begin{equation}\label{eq.NpC}
    \norm{\wt\N(u)}_{L^p(\pom)}\lesssim \sup_{F:\norm{\wt\C_1(\delta F)}_{L^{p'}(\pom)}=1}\abs{\iint_{\om}uF\,dX},
\end{equation}
and 
\begin{equation}\label{eq.NpC2}
    \norm{\wt\C_1(u)}_{L^{p'}(\pom)}\lesssim \sup_{F:\norm{\wt\N(\delta F)}_{L^{p}(\pom)}=1}\abs{\iint_{\om}uF\,dX}.
\end{equation}
\end{lemma}

\bp
See Lemmas A.18 and A.22 in \cite{MPT}, whose proof are based on a work from \cite{HR13}.
\ep

\begin{lemma} \label{lemTp1density}
Let $\Omega \subset \R^n$ be a 1-sided CAD, and $p\in [1,\infty)$. Then $L^\infty_c(\Omega)$ is dense in $\wt T^p_1(\Omega)$.
\end{lemma}

\bp
Let $u\in \wt T^p_1(\Omega)$. Without loss of generality, we can assume $u\geq 0$. Define $u_N(X) = \min\{N,u(X)\} \1_{\delta(X) \leq 1/N}$. We have $\wt \A_1(u-u_N)(x) \leq \wt \A_1(u)(x)$ and $\|\wt \A_1(u-u_N)\|_{p} \leq \|\wt \A_1(u)\|_{p}$, so by the Lebesgue domination theorem, we have
\[ \wt \A_1(u-u_N)(x) \to 0 \qquad \text{ whenever }  \wt \A_1(u)(x) < +\infty,\]
and by reusing the Lebesgue domination theorem, 
\[\|\wt \A_1(u-u_N)\|_{p} \to 0.\]
The lemma follows.
\ep

Let us turn to the more subtle properties of the tent spaces. We want an interpolation between tent spaces that is not included in existing literature to the best of our knowledge. We present here a version adapted to what we need.

\begin{theorem} \label{Thinterpolation}
Let $\Omega \subset \R^n$ be a 1-sided CAD and $p\in (1,\infty)$. Assume that a linear operator $\mathcal Z$ is bounded from $\wh{\wt T^p_{1}}(\Omega) := \{f\in \wt T^p_{1}(\Omega), \, [f/\delta]_\Omega:= \iint_\Omega f/\delta = 0\}$ to $\wt T^p_\infty(\Omega)$ and satisfies, for any ball $B=B(x,r)$ centered at the boundary, any integer $j\geq 5$, and any $f \in L^\infty_c(\Omega \cap B)$ such that $[f/\delta]_\Omega = 0$, we have
\begin{equation} \label{interpolation1}
\left( (2^jr)^{1-n} \int_{\pom} \Big|\wt \N[\mathcal Z(f)\1_{2^{j+1}B\setminus 2^jB}]\Big|^p \, d\sigma \right)^\frac1p \leq C g(j) r^{1-n} \int_{\pom} \Big|\wt \A_1(f)\Big| \, d\sigma,
\end{equation}
where $g(j)$ is such that $\sum_{j\geq 3} g(j) 2^{j(n-1)} <+\infty$ and $C$ is independent of $B$, $j$ and $f$. Then, for any $q\in (1,p)$, the operator $\mathcal Z$ is  bounded from $\wh{\wt T^q_{1}}(\Omega)$ to $\wt T^q_\infty(\Omega)$, that is
\[\|\mathcal Z(f)\|_{\wt T^q_\infty(\Omega)} := \|\wt \N[\mathcal Z(f)]\|_{L^q(\partial \Omega)} \leq C_q \|f\|_{\wt T^q_1(\Omega)} := C_q \|\wt \A_1 (f)\|_{L^q(\partial \Omega)}.\]
for any $f \in L^\infty_c(\Omega)$ such that $[f/\delta]_\Omega = 0$.
\end{theorem}

\bp See Appendix \ref{SAppendix}. \ep

\subsection{Elliptic theory: Introduction}

Let $\Omega \subset \R^n$ be an open domain. We say that an operator $L:=-\diver A \nabla$ is uniformly elliptic (on $\Omega$) if the matrix $A$ is measurable on $\Omega$ and if there exists $C>0$ such that 
\begin{equation} \label{defelliptic}
    \begin{array}{ll}
        |A(X)\xi \cdot \zeta| \leq C|\xi||\zeta| \qquad \text{ for } X\in \Omega, \, \xi,\zeta\in \R^n, \\
        A(X)\xi \cdot \xi \geq C^{-1} |\xi|^2  \qquad \text{ for } X\in \Omega, \, \xi\in \R^n.
    \end{array}
\end{equation}

\begin{definition}
Let $L = - \diver A \nabla $ be a uniformly elliptic operator and $h\in L^2(B\cap \Omega), \vec F \in L^2(B\cap \Omega,\R^n)$. We say that $u\in W^{1,2}(B\cap \Omega)$ is a (weak) local solution to $Lu=h - \diver \vec F$ if
\[ \iint_{B \cap \Omega} A \nabla u \cdot \nabla \varphi \, dX = \iint_{B \cap \Omega} \Big(h \varphi + \vec F \cdot \nabla \varphi\Big) \, dX \qquad \text{ for any } \varphi \in C^\infty_0(B \cap \Omega).\]

We say that $u\in W^{1,2}(B\cap \Omega)$ is a (weak) local solution to $Lu=h- \diver \vec F$ with zero Neumann data if 
\[ \iint_{\Omega} A \nabla u \cdot \nabla \varphi \, dX = \iint_\Omega \Big(h \varphi + \vec F \cdot \nabla \varphi\Big) \, dX \qquad \text{ for any } \varphi \in C^\infty_0(B).\]
\end{definition}

\begin{theorem} \label{ThLM}
Let $\Omega \subset \R^n$ be a bounded 1-sided CAD domain, and $L=-\diver A \nabla$ be a uniformly elliptic operator. Let $W_0^{1,2}(\Omega)$ be the subspace of $W^{1,2}(\Omega)$ with zero trace on $\partial \Omega$ (see for instance \cite{DFMmixed} for the construction and the property of traces) and $\widehat{W}^{1,2}(\Omega)$ be the subspace of $W^{1,2}(\Omega)$ of functions whose trace has zero average, that is,
\[
\widehat W^{1,2}(\om):=\set{f\in W^{1,2}(\om), \, \int_{\pom}\Tr f\,d\sigma=0}.
\]

Then both $W_0^{1,2}(\Omega)$ and $\widehat W^{1,2}(\Omega)$ are Hilbert space equipped with the inner product $\left<u, v\right> = \iint_\Omega \nabla u \cdot \nabla v\, dX$. So by the Lax-Milgram theorem we have
\begin{itemize}
    \item For any $h\in L^2(\Omega)$, $\vec F \in L^2(\Omega,\R^n)$, there exists a unique $u_{D} \in W_0^{1,2}(\Omega)$ such that 
    \begin{equation} \label{LMDirichlet} 
    \iint_{\Omega} A \nabla u_{D} \cdot \nabla \varphi \, dX = \iint_\Omega \Big( h \varphi  + \vec F \cdot \nabla \varphi\Big) \, dX \qquad \text{ for any } \varphi \in  W_0^{1,2}(\Omega).
    \end{equation}
    In particular, $u_D\in W_0^{1,2}(\om)$ is a solution to $-\diver A\nabla u_D=h-\diver\vec F$ in $\om$. 
    \item For any $h\in L^2(\Omega)$, $\vec F \in L^2(\Omega,\R^n)$, there exists a unique $u_{N} \in \widehat{W}^{1,2}(\Omega)$ such that 
    \begin{equation} \label{LMNeumann} 
    \iint_{\Omega} A \nabla u_{N} \cdot \nabla \varphi \, dX = \iint_\Omega \Big( h \varphi  + \vec F \cdot \nabla \varphi\Big) \, dX \qquad \text{ for any } \varphi \in  \widehat{W}^{1,2}(\Omega).
    \end{equation}
    If $\iint_\Omega h\, dX = 0$, such $u_N$ satisfies
    \begin{equation} \label{LMNeumann2} 
    \iint_{\Omega} A \nabla u_{N} \cdot \nabla \varphi \, dX = \iint_\Omega \Big(h \varphi  + \vec F \cdot \nabla \varphi \Big) \, dX \qquad \text{ for any } \varphi \in  W^{1,2}(\Omega),
    \end{equation}
    in particular, $u_N \in \widehat{W}^{1,2}(\Omega)$ is a weak solution to $Lu=h-\diver \vec F$ with zero Neumann data.
\end{itemize}
\end{theorem}

\bp
See \cite{DFMmixed} for the construction of the trace and the Poincar\'e inequalities needed to prove the fact that $(\widehat W^{1,2}(\Omega), \|\nabla .\|_2)$ is a complete space. The rest is the Lax-Milgram theorem.
\ep

Solutions given by the Lax-Milgram theorem can be represented by integrals of the Green function and the Neumann function. We start with the representation in terms of the Green function.
\begin{theorem} \label{ThGreen}
Let $\Omega \subset \R^n$ be a bounded 1-sided CAD domain, and $L=-\diver A \nabla$ be a uniformly elliptic operator. There exists a unique function $G$ defined on $\Omega \times \Omega$ such that $G(X,.)$ is continuous on $\Omega\setminus \{X\}$, locally integrable in $\Omega$, and such that, for any $h\in L^\infty_c(\Omega)$, the solution $u_{D,h}\in W_0^{1,2}(\om)$ defined in \eqref{LMDirichlet} can be represented as
\[
u_{D,h}(X) = \iint_\Omega G(X,Y) h(Y) \, dY.
\]
Moreover
\begin{enumerate}[(i)]
    \item for any $X,Y\in \Omega$, $G^T(Y,X) = G(X,Y)$, where $G^T(X,Y)$ is the Green function for the adjoint operator  $L^*=-\diver(A^T\nabla)$;
    \item for any $Y \in \Omega$ and any $r\in(0,\delta(Y)/2)$,
    \begin{equation} \label{pwtbdonG}
    \iint_{\Omega \setminus B(Y,r)} |\nabla G(X,Y)|^2 \, dX \leq Cr^{2-n},
    \end{equation}
    which implies that for any ball $B=B(x,r)$ centered on the boundary and verifying $Y \in 100B \setminus 1.1B$, we have
    \begin{equation} \label{pwtbdonG2}
        \left(\fiint_{B\cap \Omega} |\nabla_X G(X,Y)|^2 \, dX \right)^\frac12 \leq C r^{1-n},
    \end{equation}
    with a constant that depends only on $\Omega$ and $L$.
\end{enumerate}
\end{theorem}

\bp
See for instance \cite[Section 14]{DFMmixed}, that deals with a more general case, and \cite{GW82} for the original proof from which \cite{DFMmixed} is inspired. 
\ep

Let us make a step on the side to define the elliptic measure that we do not really need, but was introduced in the definition of \Dq.

\begin{theorem} \label{Thhm}
    Let $\Omega \subset \R^n$ be bounded 1-sided CAD domain, and $L=-\diver A \nabla$ be a uniformly elliptic operator. There exists a unique collection of probability measure $\{\omega^X\}_{X\in \Omega}$ on $\partial \Omega$ - called {\em elliptic measure} - such that, for any $f\in C^\infty(\R^n)$, the function constructed as
    \[u_f(X) := \int_{\partial \Omega} f(y) \, d\omega^X(y)\]
    is a weak solution $u\in W^{1,2}(\Omega)$ to $Lu_f=0$ in $\Omega$ (hence is continuous in $\Omega$), and can be extended by continuity on $\partial \Omega$ by $u_f|_{\partial \Omega} = f|_{\partial \Omega}$.
\end{theorem}

\bp
The proof can again be found in \cite[Section 14]{DFMmixed}, but in a more general context. The existence of $\omega^X$ is based on the maximum principle, while the continuity of $u_f$ is classical elliptic theory. See \cite{Ken94} for an earlier reference. 
\ep

\subsection{Estimates on solution with zero Neumann data}

Our next goal is to prove the analogue of Theorem \ref{ThGreen} for the solution given by \eqref{LMNeumann}. We start with the basic estimates.

\begin{lemma}[Cacciopoli inequality] \label{Caccio}
    Let $\Omega \subset \R^n$ be a 1-sided CAD, $L=-\diver A \nabla$ be a uniformly elliptic operator. Then there exists $C>0$ such that:
    \begin{itemize}
    \item for any ball $B=B(X,r)$ such that $2B \subset \Omega$, and any weak solution $u\in W^{1,2}(2B)$ to $Lu=0$ in $2B$, we have
    \[ \left( \fiint_B |\nabla u|^2 \, dX \right)^\frac12  \leq C  \frac1r \left( \fiint_{2B} |u|^2 \, dX \right)^\frac12;\]
    \item for any ball $B=B(x,r)$ centered on the boundary ($r<\diam \Omega$) and any weak local solution $u\in W^{1,2}(2B\cap \Omega)$ to $Lu=0$ with zero Neumann data,
    \[ \left( \fiint_{B\cap \Omega} |\nabla u|^2 \, dX \right)^\frac12  \leq C  \frac1r \left( \fiint_{2B \cap \Omega} |u|^2 \, dX \right)^\frac12.\]
    \end{itemize}
\end{lemma}

\begin{proposition}[Moser estimates] \label{Moser}
    Let $\Omega \subset \R^n$ be a 1-sided CAD, $L=-\diver A \nabla$ be a uniformly elliptic operator. Then there exists $C>0$ such that:
    \begin{itemize}
    \item for any ball $B=B(X,r)$ such that $2B \subset \Omega$, and any weak solution $u\in W^{1,2}(2B)$ to $Lu=0$ in $2B$, we have
    \[ \osc_B u \, dX \leq C \fiint_{2B} |u - u_B| \, dX;\]
    \item for any ball $B=B(x,r)$ centered on the boundary ($r<\diam \Omega$) and any weak local solution $u\in W^{1,2}(2B\cap \Omega)$ to $Lu=0$ with zero Neumann data,
    \[ \osc_{B\cap \Omega} u \, dX \leq C \fiint_{2B \cap \Omega} |u - u_{B\cap \Omega}| \, dX.\]
    \end{itemize}
\end{proposition}

\begin{proposition}[H\"older continuity] \label{Holder}
    Let $\Omega \subset \R^n$ be a 1-sided CAD, $L=-\diver A \nabla$ be a uniformly elliptic operator. Then there exist $C>0$ and $\alpha \in (0,1]$ such that:
    \begin{itemize}
    \item for any ball $B=B(X,r) \subset \Omega$, any weak solution $u\in W^{1,2}(B)$ to $Lu=0$ in $B$, and any $\epsilon \in (0,1)$, we have
    \[ \osc_{\epsilon B} u \, dX \leq C \epsilon^\alpha \osc_B u;\]
    \item for any ball $B=B(x,r)$ centered on the boundary ($r<\diam \Omega$) and any weak solution $u\in W^{1,2}(B\cap \Omega)$ to $Lu=0$ with zero Neumann data,
    \[ \osc_{\epsilon B\cap \Omega} u \, dX \leq C \epsilon^\alpha \osc_{B\cap \Omega}.\]
    \end{itemize}
\end{proposition}

\noindent {\em Proof of Lemma \ref{Caccio} and Propositions \ref{Moser} and \ref{Holder}:} The first part (that is, interior estimates) of each result is classical. The proof of the second part can be found in \cite{DDEMM??} or \cite{Steve'student}.
\ep

\begin{theorem} \label{ThNeumann}
Let $\Omega \subset \R^n$ be bounded 1-sided CAD domain, and $L=-\diver A \nabla$ be a uniformly elliptic operator. There exists a unique function $N$ - called Neumann function - defined on $\Omega \times \Omega$ such that $N(X,\cdot)$ is continuous on $\overline{\Omega}\setminus \{X\}$, locally integrable in $\Omega$, $\int_{\pom}N(X,y)d\sigma(y)=0$, and that the 
 following representation formula holds: for any $Y\in\overline{\om}$, 
any $\Phi\in W^{1,2}(\om)$ and $\phi=\Tr(\Phi)$, 
\begin{equation}\label{eq.RieszNmn'}
     \iint_{\om}A(X)\nabla_X N(X,Y)\cdot \nabla \Phi(X)dX=\Phi(Y)- \fint_{\pom}\phi(y)d\sigma(y).
\end{equation}
Moreover,
\begin{enumerate}[(i)]
    \item for any $X,Y\in \Omega$, $N^T(Y,X) = N(X,Y)$, where $N^T(X,Y)$ is the Neumann function for the adjoint operator  $L^*=-\diver(A^T\nabla)$;
    \item for any ball $B=B(x,r)$ centered on the boundary and verifying $Y \in \Omega \setminus 1.1B$, we have
    \begin{equation} \label{ptwbdonN}
        \left(\fiint_{B\cap \Omega} |\nabla_X N(X,Y)|^2 \, dX \right)^\frac12 \leq C r^{1-n},
    \end{equation}
    or equivalently
    \begin{equation} \label{ptwbdonN2}
       \osc_{X\in B\cap \Omega} N(X,Y) \leq C r^{2-n},
    \end{equation}
    with a constant that depends only on $\Omega$ and $L$. 
\end{enumerate}
\end{theorem}

\bp
When $n\geq 3$ and the domain is a ball, the proof is written in \cite{KP93}, but their proof doesn't trivially extend to $n=2$ and 1-sided CAD (because they use `reflection' across the boundary). However, their construction/existence and the property (i) of the Neumann function is valid for all open domains (as long as we have a boundary Poincar\'e inequality, as in \cite{DFMmixed} or \cite{HMTbook14}). In particular, one can construct the Neumaman function as \[N(X,Y):=G(X,Y)+v(X,Y),\]
where $v(\cdot,Y)\in \widehat W^{1,2}(\om)$ is a solution given by the Lax-Milgram theorem (for each $Y\in\om$) to   \begin{equation}\label{eq.defv}
    \iint_\Omega A(X) \nabla_X v(X,Y)\cdot \nabla \varphi (X) \, dX = \int_{\partial \Omega} \varphi(x) \, d\omega_T^Y(x) \qquad \text{ for } \varphi \in \widehat W^{1,2}(\Omega),
\end{equation}
where $\omega_T$ is the the elliptic measure associated to the adjoint operator $L^*$. 

The proof of \eqref{ptwbdonN} when $n\geq 3$ is in \cite{Steve'student}, while \eqref{ptwbdonN2} is a simple consequence of the Moser estimate (Proposition \ref{Moser}) and the Poincar\'e inequality (Proposition \ref{prop.Poincare}). 

\medskip

The bound \eqref{ptwbdonN} when $n=2$ is easy. We just need to observe that the norm of the functional $\varphi \to \int_{\partial \Omega} \varphi(x) \, d\omega_T^Y(x)$ in $[\widehat W^{1,2}(\Omega)]^*$ is bounded by $C\delta(Y)^{1-n/2} = C$; see the proof of Lemma 2.3 in \cite{KP93}, and the equality is due to the fact that $n=2$. So if $B$ is the ball in \eqref{ptwbdonN}, and $v$ is as in \eqref{eq.defv}, by letting $\vp=v(\cdot,Y)$ in \eqref{eq.defv}, we obtain that
\[\left(\fiint_{B} |\nabla_X v(X,Y)|^2 \, dX\right)^\frac12 \lesssim \frac1r \left(\iint_{B} |\nabla_X v(X,Y)|^2 \, dX\right)^\frac12 \lesssim \frac1r\]
since $n=2$. We conclude by using the bound \eqref{pwtbdonG2} on the Green function.
\ep

By taking transpose of \eqref{eq.RieszNmn'} and the property $(i)$ in Theorem \ref{ThNeumann}, we get the following representation formula for the solution $u_N$.
\begin{proposition}
   For any $h\in L^\infty_c(\Omega)$, $\vec F \in L^\infty_c(\Omega,\R^n)$, the solution $u_{N}$ defined in \eqref{LMNeumann} can be represented as
\[
u_{N}(X) = \iint_\Omega \Big( N(X,Y) h(Y) + \nabla_Y N(X,Y) \cdot \vec F(Y)\Big) \, dY.
\]
\end{proposition}

\section{Duality results} \label{S3}

In this section, we shall present why the Poisson-Neumann problem and the Poisson-Neumann-regularity are equivalent.
Like we say in the introduction, the result is a simple application of the duality given by \eqref{eq.CNdual} and \eqref{eq.NpC}.

\begin{lemma} \label{PNq<=>PRNp}
Let $\Omega \subset \R^n$ be a bounded 1-sided CAD, $L=-\diver A \nabla$ be a uniformly elliptic operator and $p\in (1,\infty)$. The following are equivalent,
\begin{enumerate}[(i)]
\item The Poisson-Neumann problem \PNq$_{L^*}$ is solvable;
\item the Poisson-Neumann-Regularity \PRNp$_L$ is solvable;
\item the Poisson-Neumann-Regularity \PRNp$_L$ is solvable for $\vec F =0$.
\end{enumerate}
\end{lemma}

\bp We just need to prove $(i) \implies (ii)$ and $(iii) \implies (i)$, as $(ii) \implies (iii)$ is trivial.

The proof is done by duality and is almost identical to the proof of Theorem 1.22, (c)$\implies$(d) and (e) $\implies$ (c) in \cite{MPT}. 

\medskip

We shall prove first $(iii) \implies (i)$. Let $\vec F \in L^\infty_c(\Omega)$ and let $u$ is the weak solution to $L^*u= - \div \vec F$ with zero Neumann data. We want to prove that 
\[\|\wt \N(u)\|_{L^{p'}(\partial \Omega)} \lesssim \|\wt \C_1(\delta |\vec F|)\|_{L^{p'}(\partial \Omega)}.\]

Let $K$ any compact subset of $\Omega$. In this case $\|\wt \N(u\1_K)\|_{p'} < +\infty$, and so by the duality \eqref{eq.NpC}, and the density of $L^\infty_c(\Omega)$ in the space $\wt T^p_1(\Omega)$, there exists $g = g_K \in L^\infty_c(\Omega)$ such that 
\begin{equation} \label{loc6b}
    \|\wt \C_1(\delta g)\|_{L^{p'}(\partial \Omega)} \lesssim 1
\end{equation}
and
\[I := \|\wt \N(u\1_K)\|_{L^{p'}(\partial \Omega)} \leq \iint_\Omega u g \, dX.\]
Since $u$ has zero average on $\Omega$, the right-hand side above do not change if we replace $g$ by $g-\fiint_\Omega g$. So without loss of generality, we can assume that $g\in \widehat{L^\infty_c}(\Omega)$.

\medskip

We let $v$ be a weak solution to $L v= g$ in $\Omega$ with zero Neumann data, and the bound on $I$ becomes
\begin{equation} \label{loc8b}
    I \leq \iint_\Omega u\,  Lv \, dX = \int_\Omega A^* \nabla u\,  \nabla v \, dX = \iint_\Omega \vec F \cdot \nabla v \, dX
\end{equation}
because $u$ is a weak solution to $L^* u = -\div \vec F$ with zero Neumann data. By the Carleson inequality \eqref{eq.CNdual}, $(iii)$ applied to $v$, and \eqref{loc6b}, we have that
\begin{multline*} 
\|\wt \N(u\1_K)\|_{L^{p'}(\partial \Omega)} = I \lesssim \|\wt \C_1(\delta \vec F)\|_{L^{p'}(\partial \Omega)} \|\wt \N(\nabla v)\|_{L^{p}} \\
\lesssim \|\wt \C_1(\delta \vec F)\|_{L^{p'}(\partial \Omega)} 
\|\wt \C_1(\delta g)\|_{L^{p'}(\partial \Omega)} \lesssim \|\wt \C_1(\delta \vec F)\|_{L^{p'}(\partial \Omega)}.
\end{multline*}
Since the bound is independent of the compact $K$, we take $K \uparrow \Omega$ and we have the desired result.

\medskip

We turn to the proof of $(i) \implies (ii)$. Let $h\in \widehat{L^\infty_c(\Omega)}$, $\vec F \in L^\infty_c(\Omega)$ and let $u$ be the weak solution to $L^*u= - \div \vec F$ with zero Neumann data. We want to prove that 
\[\|\wt \N(\nabla u)\|_{L^{p}(\partial \Omega)} \lesssim \|\wt \C_1(\vec F)\|_{L^{p}(\partial \Omega)} + \|\wt \C_1(\delta h)\|_{L^{p}(\partial \Omega)}.\]

We take again a compact $K\subset \Omega$. By the duality \eqref{eq.NpC}, and by the density of $L^\infty_c(\Omega)$ in the space $T^{p'}_1(\Omega)$, there exists $\vec G = \vec G_K \in L^\infty_c(\Omega)$ such that 
\begin{equation} \label{loc6c}
    \|\wt \C_1(\delta \vec G)\|_{L^{p}(\partial \Omega)} \lesssim 1
\end{equation}
and
\[II:= \|\wt \N(\nabla u \1_K)\|_{L^{p}(\partial \Omega)} \leq \iint_\Omega \nabla u \cdot \vec G \, dX = -\iint_\Omega u \, (\diver \vec G)\, dX.\]

\medskip

We write $v$ for a weak solution to $L^* v= -\diver \vec G$ in $\Omega$ with zero Neumann data, and the bound on $II$ becomes
\begin{equation} 
    II \leq \iint_\Omega u\,  Lv \, dX = \iint_\Omega A \nabla u \cdot  \nabla v \, dX = \iint_\Omega \big[ hv + \vec F \cdot \nabla v\big] \, dX
\end{equation}
because $u$ is a weak solution to $L u = h -\div \vec F$ with zero Neumann data. By the Carleson inequality \eqref{eq.CNdual}, we have that
\[ II \lesssim \|\wt \C_1(\delta h)\|_{L^{p}(\partial \Omega)} \|\wt \N(v)\|_{L^{p'}} +  \|\wt \C_1(\vec F)\|_{L^{p}(\partial \Omega)} \|\wt \N(\delta \nabla v)\|_{L^{p'}}\]
But due to the Cacciopolli inequality (Lemma \ref{Caccio}), for any $\xi \in \partial \Omega$ 
\[\wt \N(\delta \nabla v)(\xi)\lesssim \wt \N^*(\vec G)(\xi) + \wt \N^*(v)(\xi),\]
where $\wt \N^*$ is a non-tangential maximal function defined with cones of larger aperture. Since $\wt \N$ and $\wt \N^*$ have equivalent $L^{p'}$ norms (Lemma \ref{lemNCA}), we deduce
\[
\|\wt \N(\delta \nabla v)\|_{L^{p'}} + \|\wt \N(v)\|_{L^{p'}} \lesssim \|\wt \N(\vec G)\|_{L^{p'}} + \|\wt \N(v)\|_{L^{p'}} \\ \lesssim \|\wt \C_1(\vec G)\|_{L^{p'}} \lesssim 1.
\]
where, for the second inequality, we apply (i) to $v$ and we use the fact that $\wt \N(\vec G)(\xi) \lesssim \wt \C_1(\vec G)(\xi)$, and the last inequality is \eqref{loc6c}. Going back to $II$, we have 
\[\|\wt \N(\nabla u \1_K)\|_{L^{p}(\partial \Omega)} = II \lesssim \|\wt \C_1(\delta h)\|_{L^{p}(\partial \Omega)} +  \|\wt \C_1(\vec F)\|_{L^{p}(\partial \Omega)}.\]
The lemma follows by taking $K\uparrow \Omega$.
\ep

\begin{lemma} \label{PNq=>Np}
 Let $\Omega \subset \R^n$ be a 1-sided CAD, $L = - \div A \nabla$ be a uniformly elliptic operator, and $p\in (1,\infty)$. Then
 \[ \text{\PNq$_{L^*}$ $\implies$ \Np$_L$}.\]
\end{lemma}

\bp The proof is a variant of $(i) \implies (ii)$ from Lemma \ref{PNq<=>PRNp}, to which we refer for details. Let $h\in L^p(\partial \Omega)$ and let $u$ be the weak solution to $Lu=0$ with Neumann data $h$. 

Take a compact $K\subset \Omega$, and then construct vector $\vec G \in L^\infty_c(\Omega)$, and a weak solution $v \in W^{1,2}(\Omega)$ to $L^*v = -\diver \vec G$ with zero Neumann data and zero average such that $\|\wt \C_1(\delta \vec G)\|_{p'} \lesssim 1$ and 
\[\|\wt \N(\nabla u\1_K)\|_{p} \leq \iint_\Omega \nabla u \cdot \vec G \, dX = \iint_\Omega A \nabla u \cdot \nabla v \, dX = \int_{\partial \Omega} h \Tr(v)\, d\sigma.\]
We can use the H\"older inequality to further bound 
\[
\|\wt \N(\nabla u\1_K)\|_{p} \leq \|h\|_{p} \|\Tr(v)\|_{p'} \leq \|h\|_{p} \|\wt \N(v)\|_{p'} \lesssim \|h\|_p \|\wt \C_1(\delta \vec G)\|_{p'} \lesssim \|h\|_p 
\]
by \PNq$_{L^*}$ applied to $v$ and the fact that $\|\wt \C_1(\delta \vec G)\|_{p'} \lesssim 1$. The lemma follows.
\ep

\begin{lemma} \label{wPNq=>wPRNp}
Let $\Omega \subset \R^n$ be a 1-sided CAD, $L=-\diver A \nabla$ be a uniformly elliptic operator and $p\in (1,\infty)$. Then \wPNq$_{L^*}$ is equivalent to \wPRNp$_L$.
\end{lemma}

\bp The proof is done again by duality, it is analogous to the proof of Theorem 1.22, (c)$\implies$(d) and (e) $\implies$ (c) in \cite{MPT}. 
Let us demonstrate only \wPRNp$_{L}$ $\implies$ \wPNq$_{L^*}$, since the converse is almost identical. 

\medskip

Let $\vec F \in L^\infty_c(\Omega)$, and
\begin{equation} \label{loc7}
u(X):= \iint_\Omega \nabla_Y N(Y,X) \cdot \vec F(Y) \, dY,
\end{equation}
which is a solution to $L^*u= - \div \vec F$. We want to prove that 
\[I:= \|\wt \N(\delta \nabla u)\|_{L^{p'}(\partial \Omega)} \lesssim \|\wt \C_1(\delta \vec F)\|_{L^{p'}(\partial \Omega)}.\]

By the duality \eqref{eq.NpC}, and by the density of $L^\infty_c(\Omega)$ in the space $\mathbf{C}_p$, there exists $\vec G \in L^\infty_c(\Omega,\R^n)$ such that 
\begin{equation} \label{loc6}
    \|\wt \C_1(\delta \vec G)\|_{L^{p}(\partial \Omega)} \lesssim 1
\end{equation}
and
\[I \leq \iint_\Omega \delta \nabla u \cdot \vec G \, dX = -\iint_\Omega u \diver[\delta \vec G] \, dX.\]
We let $v$ be the weak solution to $L v= - \diver[\delta G]$ in $\Omega$ with zero Neumann data, that is
\[v(X) := \iint_\Omega \nabla_Y N(X,Y) \cdot [\delta \vec G](Y) \, dY,\]
and the bound on $I$ becomes
\begin{equation} \label{loc8}
    I \leq \iint_\Omega u\,  Lv \, dX = \iint_\Omega A^* \nabla u\,  \nabla v \, dX = \iint_\Omega \vec F \cdot \nabla v \, dX
\end{equation}
because $u$ is a weak solution to $L^*u = -\div \vec F$ with zero Neumann data. Using successively the Carleson inequality \eqref{eq.CNdual}, \wPRNp$_{L}$ to $v$, and \eqref{loc6} entails then that
\[ I \lesssim \|\wt \C_1(\delta \vec F)\|_{L^{p'}(\partial \Omega)} \|\wt \N(\nabla v)\|_{L^{p}} \lesssim \|\wt \C_1(\delta \vec F)\|_{L^{p'}(\partial \Omega)} \|\wt \C_1(\delta \vec G)\|_{L^{p}(\partial \Omega)} \lesssim \|\wt \C_1(\delta \vec F)\|_{L^p(\partial \Omega)}.\]
The lemma follows.
\ep

\section{Localization} \label{S4}

In this section, we want to prove the localization of the weak Poisson-Neumann problem and of the strong Neumann problem. 

We start with the following lemma, which proves a stronger statement than Theorem \ref{MainTh} $(a) \implies (c)$, and will yield localization of the weak Poisson-Neumann problem (Corollary \ref{wPNq=>Lp2}).

\begin{lemma} \label{wPNq=>Lp}
Let $\Omega \subset \R^n$ be a bounded 1-sided CAD, $L=-\diver A \nabla$ be a uniformly elliptic operator and $p\in (1,\infty)$. If \wPNq$_{L^*}$ holds, then for any nonnegative $\alpha,\beta,\gamma$ that satisfy $\max\{1,\alpha\} < \beta \leq \gamma \leq 2$, there exists $C=C(\beta-\alpha)>0$ such that for any ball $B=B(x,r)$ centered on the boundary, and any local solution $u \in W^{1,2}(\gamma B\cap \Omega)$ to $Lu=0$ in $\gamma B\cap\om$ with zero Neumann data, we have
\begin{equation} \label{locCase1}
\|\wt \N(|\nabla u|\1_B)\|_{L^p(\partial \Omega)} \leq C r^{(n-1)/p} \fiint_{\beta B\setminus \alpha B} |\nabla u| \, dX .
\end{equation}
Similarly, if $u\in W^{1,2}(\Omega \setminus \frac1{\gamma}B)$ is a local solution to $Lu=0$ with zero Neumann data, then 
\begin{equation} \label{locCase2}
\|\wt \N(|\nabla u|\1_{\Omega \setminus B})\|_{L^p(\partial \Omega)} \leq C r^{(n-1)/p} \fiint_{\frac1{\alpha}B \setminus \frac{1}{\beta}B} |\nabla u| \, dX .
\end{equation}
\end{lemma}

\bp
We shall only prove \eqref{locCase1} when $\alpha=1$ and $\beta = \gamma =2$ (which trivially implies \Lp$_L$) and we let the reader check that the proof for \eqref{locCase2} and other values of $\alpha$, $\beta$ and $\gamma$ can easily be adapted.

\medskip

\noindent {\bf Step 1} We claim that to show that  \eqref{locCase1} holds for $\alpha = 1$, $\beta = \gamma = 2$, it suffices to prove that there exists $\Lambda \geq 1$ such that for any ball $B' = B(x',r')$ centered at the boundary and any local solution $u \in W^{1,2}(\Lambda B'\cap \Omega)$ to $Lu=0$ in $\Lambda B\cap \om$ with zero Neumann data, there holds
\begin{equation} \label{loc99}
\|\wt \N(|\nabla u|\1_{B'})\|_{L^p(\partial \Omega)} \lesssim r'^{(n-1)/p}\Big( \fiint_{\Lambda B' \cap \Omega} |\nabla u|^2 dX \Big)^\frac12.
\end{equation}

\medskip

In fact, if \eqref{loc99} holds, then for a given ball $B= B(x,r)$ centered on the boundary, we cover $B$ by a finite collection of balls $B'_i$ of radius $r_i'\approx r$ such that 
\begin{enumerate}[(i)]
    \item either $B'_i$ is centered at the boundary and $\Lambda B'_i \subset \frac98B$, 
    \item or $2B'_i\subset \frac98B \cap \Omega$.
\end{enumerate}
As there are finite many $B_i'$ and the number depends only on the dimension and the CAD constants,  we have that
\[\|\wt \N(|\nabla u|\1_{B})\|_{L^p(\partial \Omega)} \lesssim \sup_{i} \|\wt \N(|\nabla u|\1_{B_i'})\|_{L^p(\partial \Omega)}.\]
If $B'_i$ is as in (i), then we use \eqref{loc99} to get that 
\[
\|\wt \N(|\nabla u|\1_{B_i'})\|_p\lesssim (r')^{(n-1)/p}\br{\fiint_{\Lambda B_i'\cap\om}\abs{\nabla u}^2}^{1/2}\lesssim r^{(n-1)/p}\br{\fiint_{\frac98 B\cap\om}\abs{\nabla u}^2}^{1/2}.
\]
If $B_i'$ is as in (ii), simply observe that $\wt \N(|\nabla u|\1_{B_i'})$ is supported on a large shadow if $B_i'$ on $\pom$, which is of measure $C(r')^{n-1}$. Therefore,
\[ \|\wt \N(|\nabla u|\1_{B_i'})\|_{L^p(\partial \Omega)} \approx (r')^{(n-1)/p} \left(\fiint_{B_i'} |\nabla u|^2 \, dX \right)^\frac12 \approx r^{(n-1)/p} \left(\fiint_{\frac98B\cap\om} |\nabla u|^2 \right)^\frac12.\]
Ultimately, we have showed that
\begin{equation} \label{loc98}
\|\wt \N(|\nabla u|\1_{B})\|_{L^p(\partial \Omega)} \lesssim r^{(n-1)/p}\Big(\fiint_{\frac98 B \cap \Omega} |\nabla u|^2 dX \Big)^\frac12.
\end{equation}

So it remains to show that we can change the right-hand side of \eqref{loc98} from a $L^2$ average to a $L^1$ average, and to show that we can integrate over $2B\setminus B$ instead of $\frac98B$. 
We want to first use the Caccioppoli inequality and then use the Poincar\'e inequality to control the right-hand side of \eqref{loc98}, but since the set $\frac{9}{8}B\cap\om$ might not be a John domain (it might not be connected), we need to be careful when applying the Poncar\'e inequality. Therefore, we construct the set $S$ as
\[ S := {\rm Interior}\left( \bigcup_{\begin{subarray}{c} Y,Z\in \br{\frac74 B \setminus \frac65 B} \cap \Omega \\ \delta(Y)+ \delta(Z)+|Y-Z|< r/K \end{subarray}} \overline{H_\Omega(Y,Z)}\right),\]
where $H_\Omega(Y,Z)$ is the Harnack chain of Whitney cubes between $Y$ and $Z$ (see Definition \ref{defHC}). If $K$ is sufficiently large (depending on the constant in the Harnack chain condition), $S$ is contained in $\frac{15}{8}B \setminus \frac{7}{6}B$. Moreover, $S$ is a union of disjoint 1-sided John domains that contains $\br{\frac74B\setminus\frac65 B}\cap \Omega$; in particular, each connected component $S'$ of $S$ is a John domain and hence satisfies the Poincar\'e inequality
\begin{equation} \label{S'Poincare}
\iint_{S'} |w-w_{S'}| \, dZ \lesssim r \iint_{S'} |\nabla w| \, dZ 
\end{equation}
whenever $w\in W^{1,2}(S')$ and $w_{S'}$ is the average of $w$ on the John center of $S'$ (see Definition \ref{defHC}). So if $w\in W^{1,2}(S)$ and 
\begin{equation} \label{defwS}
w_S(Z) := w_{S'} \quad \text{when $Z \in S'$},
\end{equation}
then we have
\begin{equation} \label{SPoincare}
\fiint_{S} |w-w_S| \, dZ  \lesssim r  \fiint_{S} |\nabla w| \, dZ .
\end{equation}

Let us come back to $u$. Since $u-u_S$ is a solution to $Lu=0$ on each connected component of $\frac87B$, the Caccioppoli inequality and the Moser estimate entail that 
\begin{multline*}
    \left( \fiint_{\frac98 B \cap \Omega} |\nabla u|^2 dX\right)^\frac12  \lesssim \frac1{r}  \left( \fiint_{(\frac32B \setminus \frac54B) \cap \Omega} |u-u_S|^2 dX \right)^\frac12
    \lesssim r^{-1} \sup_{(\frac32B \setminus \frac54B)\cap\om}\abs{u-u_S}\\
    \lesssim  r^{-1}\fiint_{(\frac74 B \setminus \frac65B) \cap \Omega} |u-u_S| dX,\end{multline*}
so \eqref{loc98} becomes
\begin{multline*}
    \|\wt \N(|\nabla u|\1_{B})\|_{L^p(\partial \Omega)} \lesssim r^{\frac{n-1}{p} - 1} \fiint_{S} |u-u_S| dX  \lesssim r^{(n-1)/p} \fiint_{S} |\nabla u| dX\\
    \lesssim r^{(n-1)/p} \fiint_{2B\setminus B} |\nabla u| dX .
\end{multline*}
by \eqref{SPoincare} and by definition of $S$. Step 1 follows.

\medskip

\noindent {\bf Step 2:} Let $B=B(x,r)$ be a ball centered at the boundary. We assume that $u \in W^{1,2}(\Lambda B\cap \Omega)$ is a weak solution to $Lu=0$ with zero Neumann data for some $\Lambda\geq 1$ large (independent of $B$ and $u$) to be determined during the proof, and we want to prove \eqref{loc99}, that is
\begin{equation} \label{loc97}
\|\wt \N(|\nabla u|\1_{B})\|_{L^p(\partial \Omega)} \lesssim r^{(n-1)/p}\Big( \fiint_{\Lambda B \cap \Omega} |\nabla u|^2 dX \Big)^\frac12.
\end{equation} 
We construct then $\varphi\in C^\infty_0(\R^n)$ such that $\varphi \equiv 1$ on $\frac32 B$, $\varphi \equiv 0$ outside $2B$, and $|\nabla \varphi| \leq \frac4r$. 
By duality, i.e. \eqref{eq.NpC}, there exists $\vec G \in L^\infty_c(\Omega)$ such that
\begin{equation} \label{loc2}
    \|\wt \C_1(\delta \vec G)\|_{L^{p'}(\partial \Omega)} \leq 1
\end{equation}
and
\begin{equation} \label{loc3}
    \|\wt \N(\nabla u \1_B)\|_{L^p(\partial \Omega)} \leq C \iint_\Omega \1_B\nabla u \cdot \vec G\, dX
\end{equation}
Without loss of generality, we can assume that $\vec G = \vec G \1_B$. Let $c_u$ be a constant to be determined later, we have that
\begin{multline*}
    \|\wt \N(\nabla u \1_B)\|_{L^p(\partial \Omega)} \lesssim \iint_\Omega \nabla u \cdot \vec G\, \varphi \, dX = \iint_\Omega \nabla(u-c_u) \cdot \vec G\, \varphi \, dX\\
    = -\iint_\Omega (u-c_u) \varphi \diver \vec G \, dX - \iint_{\Omega} (u-c_u) \nabla \varphi \cdot \vec G \, dX =: I_0 + I_1
\end{multline*}
We define $v$ to be the weak solution to $L^*v= - \diver \vec G$ in $\Omega$ with zero Neumann data, that is
\[v(X) := \iint_\Omega \nabla_Y N(Y,X) \cdot \vec G(Y) \, dY.\]
One has 
\begin{multline*}
I_0 =  \iint_\Omega (u-c_u) \varphi  L^* v \, dX = \iint_\Omega A \nabla [(u-c_u)\varphi] \cdot \nabla v \, dX \\
= \iint_\Omega A \nabla u \cdot \nabla v \, \varphi \, dX + \iint_\Omega A \nabla \varphi \cdot \nabla v \, (u-c_u) \, dX =: I_{00} + I_2.
\end{multline*}
Using the fact that $Lu = 0$ in $\Omega \cap B$, we further obtain,
\begin{multline} \label{loc95}
I_{00} = \iint_\Omega A \nabla u \cdot \nabla [(v-c_v) \varphi] \, dX  - \iint_\Omega A \nabla u \cdot \nabla \varphi \, (v-c_v) \, dX
\\ 
=  - \iint_\Omega A \nabla u \cdot \nabla \varphi \, (v-c_v) \, dX
=: I_{3},
\end{multline}
where $c_v$ is another constant to be determined later.

\medskip

It remains to bound $I_1$, $I_2$ and $I_3$. Let us treat $I_1$ first. Choosing $c_u=\fiint_{4B\cap\om}u$ and $K$ large enough, one obtains that
\begin{multline*}
I_1 \lesssim r^{n-1} \left( \fiint_{2B \cap \Omega} |\vec G| \, dX\right) \left( \sup_{2B \cap \Omega} |u-c_u| \right)  \lesssim r^{n-1} \left( \inf_{x\in 2B \cap \partial \Omega} \wt \C_1(\delta G)(x) \right) \left( \fiint_{4KB \cap \Omega} |\nabla u|^2 \right)^\frac12 
\end{multline*}
by the definition of $\wt \C_1$, the Moser estimate for $u-c_u$, and the Poincar\'e inequality (Proposition \ref{prop.Poincare}, where we have fixed $c_u = \fiint_{4B\cap\om} u$). Since
\[\inf_{x\in 2B \cap \partial \Omega} \wt \C_1(\delta G)(x) \lesssim r^{(1-n)/p'}\|\wt \C_1(\delta G) \|_{L^{p'}(\partial \Omega)} \lesssim r^{(1-n)/p'} \]
by \eqref{loc2}, and that $1-\frac1{p'} = \frac1p$, we deduce that
\[I_1 \lesssim r^{(n-1)/p} \left( \fiint_{4KB \cap \Omega} |\nabla u|^2 \right)^\frac12\]
as desired.

\medskip

Before we deal with $I_2$ and $I_3$, let's take a break to prove the following estimate: for some $c_v$ that will be specified shortly, and $K$ large enough, there holds
\begin{equation} \label{loc.claim}
\fiint_{4B\cap\om} |v-c_v| \, dX \lesssim r^{1-n} \|\wt \N(\delta|\nabla v|\1_{4KB})\|_{L^{1}(\partial \Omega)}.
\end{equation}
Observe that
\begin{multline}\label{loc.claim1}  
\fiint_{4B\cap\om} |v-c_v| \, dX 
\lesssim r^{-n}\int_{y\in\pom}\iint_{\gamma(y)}\abs{v(X)-c_v}\1_{4B}(X)\frac{dX}{\delta(X)^{n-1}}d\sigma(y)\\
\lesssim r^{-n}  \int_{y \in\partial \Omega} \sum_{W \in \mathcal W \atop y \in 100W} \ell(W) \left( \fiint_{W} |v-c_v|^{2} \1_{4B} \, dZ \right)^\frac12 d\sigma(y).
\end{multline}
We take $c_v$ to be the average of $v$ over $W_{B}$, where $W_B \in \mathcal W$ is a Whitney cube associated to $B$ (that is $W_B \subset B$, $\delta(W_B) := \dist(W_B,\partial \Omega) \gtrsim r$). For a given $W$ such that $W\cap 4B \neq \emptyset$, we construct the Harnack chain of Whitney cubes $W=W_1,\dots,W_{N}=W_B$. Such Harnack chain satisfies 
\begin{enumerate}
    \item $N \lesssim 1 + \ln(r/\delta(W))$ by definition;
    \item for any $1\leq i \leq N$ we have $W_i \in \gamma^*(y)$ for $y\in 100W \cap \partial \Omega$ and a cone $\gamma^*(y)$ with large enough aperture since $y\in 100W$;
    \item $W_i \subset 4KB$ for $K$ large enough as $W \cap 4B \neq \emptyset$. 
\end{enumerate}
Properties on the Harnack chains and the Poincar\'e inequalities
allow us to show that 
\begin{equation}\label{eq.v-cv}
    \left(\fiint_{W} |v-c_v|^{2} \1_{4B} \, dZ\right)^\frac12 \lesssim \sum_{i=1}^N \left(\fiint_{W_i} |\delta \nabla v|^{2} \1_{4KB} \, dZ\right)^\frac12.
    \end{equation} 
In fact, recalling $c_v=\fiint_{W_N}v\,dZ$, we can write
\[
\left(\fiint_{W} |v-c_v|^{2} \1_{4B} \, dZ\right)^\frac12 \le \br{\fiint_{W}\abs{v-(v)_{W_1}}^2\1_{4B}}^{1/2}+\sum_{i=2}^N\abs{(v)_{W_{i-1}}-(v)_{W_i}},
\]
where $(v)_{W_i}:=\fiint_{W_i}v\,dZ$. Since $W_i$ and $W_{i-1}$ are adjacent, we can find a ball $S_i\subset W_{i-1}\cup W_i$ such that both $S_i\cap W_{i-1}$ and $S_i\cap W_i$ have measures comparable to $\ell(W_i)^n$. By applying the triangle inequality and the Poincar\'e inequality several times, we have 
\begin{multline*}
    \abs{(v)_{W_{i-1}}-(v)_{W_i}}\\
    \le \abs{(v)_{W_{i-1}}-(v)_{S_i\cap W_{i-1}}}+ \abs{(v)_{S_i\cap W_{i-1}} (v)_{S_i}}
    +\abs{(v)_{S_i}-(v)_{S_i\cap W_{i}}}+\abs{(v)_{S_i\cap W_i}-(v)_{W_i}}\\
    \lesssim \fiint_{W_{i-1}}\abs{\delta\nabla v}dZ+\fiint_{S_i}\abs{\delta\nabla v}dZ+\fiint_{W_i}\abs{\delta\nabla v}dZ\\
    \lesssim \fiint_{W_{i-1}}\abs{\delta\nabla v}dZ+\fiint_{W_i}\abs{\delta\nabla v}dZ,
\end{multline*}
which gives \eqref{eq.v-cv} by the H\"older inequality and property (3) of the Whitney cubes $W_i$.  
Note that property (2) ensures that $\left(\fiint_{W_i} |\delta \nabla v|^{2} \1_{4KB} \, dZ\right)^\frac12\le \wt\N^*(\delta |\nabla v|\1_{4KB})(y)$ for all $1\le i\le N$, where $\wt \N^*$ denotes the modified non-tangential maximal function defined with the cones $\gamma^*$. Now by property (1) of the Harnack chain and \eqref{eq.v-cv}, one sees that 
\[\left(\fiint_{W} |v-c_v|^{2} \1_{4B} \, dZ\right)^\frac12 \lesssim \ln\Big(\frac{r}{\delta(W)} \Big)\wt \N^*(\delta |\nabla v|\1_{4KB})(y).\]
 Plugging this estimate into \eqref{loc.claim1}, we obtain that
\begin{multline*}
    \fiint_{4B\cap\om} |v-c_v| \, dX \lesssim r^{-n}  \int_{\partial \Omega} \wt \N^*(\delta |\nabla v| \1_{4KB} )(y) \sum_{\begin{subarray}{c} W \in \mathcal W \\ y \in 100W \\ W \cap 4B \neq \emptyset\end{subarray}} \ell(W) \ln\Big(\frac{r}{\delta(W)}\Big)   \, d\sigma(y) \\
    \lesssim r^{1-n}  \int_{\partial \Omega} \wt \N^*(\delta |\nabla v| \1_{4KB})(y) \, d\sigma(y).
\end{multline*}
To see why the last inequality holds, we check that for a given $y\in \partial \Omega$ and $k\geq -4$, there is a uniformly finite number of $W \in \mathcal W$ satisfying $100W \ni y$ and $\ell(W) =2^{-k}r$ for each $k\in\mathbb{Z}$, and therefore,
\[
\sum_{\begin{subarray}{c} W \in \mathcal W \\ y \in 100W \\ W \cap 4B \neq \emptyset\end{subarray}} \ell(W) \ln\Big(\frac{r}{\delta(W)}\Big)  
\lesssim \sum_{k= -4}^\infty k2^{-k}r  \lesssim r.
\]
Yet the equivalence of $L^{1}$-norms of non-tangential maximal function with different aperture (Lemma \ref{lemNCA}) yields that 
\[ \|\wt \N^*(\delta |\nabla v| \1_{4KB} )\|_{L^1} \approx \|\wt \N(\delta |\nabla v| \1_{4KB})\|_{L^1},\] and hence 
the estimate \eqref{loc.claim} follows. 

We further control the right-hand side of \eqref{loc.claim} using \wPNq$_{L^*}$, and \eqref{loc2} as follows. Since $\wt \N(\delta |\nabla v| \1_{4KB} )$ is supported in $20KB\cap \partial \Omega$), by the H\"older inequality, we get that 
\begin{multline*}
    r^{1-n}\norm{\wt\N(\delta\abs{\nabla v}\1_{4KB})}_{L^1(\pom)}
\lesssim\br{\fint_{20KB\cap\pom}\wt\N(\delta\abs{\nabla v})(y)^{p'}d\sigma(y)}^{1/{p'}}\\
    \lesssim r^{(1-n)/p'}\norm{\wt\N(\delta\abs{\nabla v})}_{L^{p'}(\pom)}\lesssim r^{(1-n)/p'} \|\wt \C_1(\delta \vec G)\|_{L^{p'}(\partial \Omega)} \lesssim r^{(1-n)/p'}
    \end{multline*}
by \wPNq$_{L^*}$ and \eqref{loc2}. This and \eqref{loc.claim} yield that 
\begin{equation} \label{loc5}
    \fiint_{4B\cap\om} |v-c_v| \, dX \lesssim r^{(1-n)/p'}.
\end{equation}

\medskip
We return to $I_2$ and $I_3$.
For $I_3$, using first the Cauchy-Schwarz inequality, we get that
\begin{multline*}
    |I_3| \lesssim r^{n-1} \left( \fiint_{2B\cap \Omega} |\nabla u|^2 \, dX \right)^\frac12 \left( \fiint_{(2B \setminus \frac32B)\cap \Omega} |v-c_v|^2 \, dX \right)^\frac12\\
\le  r^{n-1} \left( \fiint_{2B\cap \Omega} |\nabla u|^2 \, dX \right)^\frac12 \sup_{(2B\setminus\frac32B)\cap \om}\abs{v-c_v}.
\end{multline*}
Since $\vec G=\vec G\1_B$, $v-c_v$ is a solution to $L^*v=0$ in $(4B\setminus B)\cap\om$, with zero Neumann data. By the Moser estimate and \eqref{loc5}, we have that 
\begin{multline*}
\abs{I_3}\lesssim r^{n-1} \left( \fiint_{2B \cap \Omega} |\nabla u|^2 \, dX \right)^\frac12 \left( \fiint_{(4B \setminus B)\cap \Omega} |v-c_v| \, dX \right)\\
\lesssim r^{(n-1)/p} \left( \fiint_{(2B \setminus B)\cap \Omega} |\nabla u|^2 \, dX \right)^\frac12.
\end{multline*}

It remains to estimate $I_2$, which is similar to $I_3$, but we have the gradient on $v$ instead of $u$. We only have to use the Cacciopoli inequality (on $v$) and the Poincar\'e inequality (Proposition \ref{prop.Poincare}) to come back to the same situation as $I_3$:
\begin{multline*}
|I_2| \lesssim r^{n-1} \left( \fiint_{2B \cap \Omega} |u-c_u|^2 \, dX \right)^\frac12 \left( \fiint_{(2B \setminus \frac32B)\cap \Omega} |\nabla v|^2 \, dX \right)^\frac12 \\
\lesssim r^{n-1} \left( \fiint_{2KB \cap \Omega} |\nabla u|^2 \, dX \right)^\frac12 \left( \fiint_{(4B \setminus B)\cap \Omega} |v-c_v| \, dX \right) \lesssim r^{(n-1)/p} \left( \fiint_{2KB \cap \Omega} |\nabla u|^2 \, dX \right)^\frac12.
\end{multline*}
This completes the proof of the lemma.
\ep

\begin{corollary}[Localization of the weak Poisson-Neumann problem] \label{wPNq=>Lp2}
Let $\Omega \subset \R^n$ be a bounded 1-sided CAD, $L=-\diver A \nabla$ be a uniformly elliptic operator and $p\in (1,\infty)$. If \wPNq$_{L^*}$ holds, then there exists $C>0$ such that for any ball $B=B(x,r)$ centered on the boundary, any $\vec F \in L^\infty_c(\Omega)$, and any local solution $u \in W^{1,2}(2B\cap \Omega)$ to $Lu = - \diver \vec F$ in $2B\cap \Omega$ with zero Neumann data on $2B\cap\pom$, we have
\[
\|\wt \N(|\nabla u|\1_B)\|_{L^p(\partial \Omega)} \leq C \left( \|\wt \C_1(\vec F \1_{2B})\|_{L^p(\partial \Omega)} +  r^{(n-1)/p} \fiint_{(2B\setminus B) \cap \Omega} |\nabla u| dX \right).
\]
\end{corollary}

\bp
Let $u_0$ be the solution to $Lu_0 = -\diver (\vec F \1_{2B})$ in $\Omega$, with zero Neumann data and $\fiint_{\om}u_0\,dX=0$. Since $L(u-u_0)=0$ in $2B\cap\om$ and that $u-u_0$ has zero Neumann data on $2B\cap\pom$, we can apply
Lemma \ref{wPNq=>Lp} to obtain that 
\[\|\wt \N(|\nabla(u-u_0)|\1_B)\|_{L^p(\partial \Omega)} \lesssim r^{(n-1)/p} \fiint_{(2B\setminus B) \cap \Omega} |\nabla (u-u_0)| dX ,\]
which implies that
\[
    \|\wt \N(|\nabla u|\1_B)\|_{L^p(\partial \Omega)} \lesssim r^{(n-1)/p} \br{\fiint_{(2B\setminus B) \cap \Omega} |\nabla u| dX + \fiint_{2B \cap \Omega} |\nabla u_0| dX }
    + \|\wt \N(\nabla u_0)\|_{L^p(\partial \Omega)}  .
\]
However, 
\begin{multline*}
    \fiint_{2B \cap \Omega} |\nabla u_0| dX \lesssim r^{-n} \int_{\partial \Omega}  \iint_{\gamma(x)} \left(\fiint_{B_Y} |\nabla u_0(Z)|^2 \1_{2B}(Z)\, dZ \right)^\frac12  \, \frac{dY}{\delta(Y)^{n-1}} \, d\sigma(x)  \\
    \lesssim r^{-n} \int_{\partial \Omega} \wt \N(|\nabla u_0| \1_{2B})(x)\iint_{\gamma(x)}\1_{3B}(Y)\frac{dY}{\delta(Y)^{n-1}} d\sigma(x).
\end{multline*}
A simple computation shows that 
\(\iint_{\gamma(x)} \1_{3B}(Y) \frac{dY}{\delta(Y)^{n-1}} \lesssim r\), and thus
\begin{multline*}
    \fiint_{\frac32B \cap \Omega} |\nabla u_0| dX \lesssim r^{1-n}\int_{\partial \Omega} \wt \N(|\nabla u_0| \1_{2B})(x)d\sigma(x)\lesssim \fint_{\Lambda B\cap\pom}\wt \N(|\nabla u_0| \1_{2B})(x)d\sigma(x)\\
    \le \br{\fint_{\Lambda B\cap\pom}\wt \N(|\nabla u_0| \1_{2B})(x)^pd\sigma(x)}^{1/p}
    \lesssim r^{-(n-1)/p}\|\wt \N(\nabla u_0 \1_{2B})\|_{L^p(\partial \Omega)},
\end{multline*}
 where the second inequality is due to the fact that $\wt \N(|\nabla u_0| \1_{2B})$ is supported on $\Lambda B \cap \partial \Omega$ for some $\Lambda\ge2$ depending only on  the aperture.  
Altogether, we have proved that
\begin{multline*}
    \|\wt \N(|\nabla u|\1_B)\|_{L^p(\partial \Omega)} \lesssim r^{(n-1)/p} \fiint_{(2B\setminus B) \cap \Omega} |\nabla u| dX + \|\wt \N(\nabla u_0)\|_{L^p(\partial \Omega)}  \\
    \lesssim r^{(n-1)/p} \fiint_{(2B\setminus B) \cap \Omega} |\nabla u| dX + \|\wt \C_1(\vec F \1_{2B})\|_{L^p(\partial \Omega)}
\end{multline*} 
by \wPRNp$_L$ (which is equivalent to our assumption \wPNq$_{L^*}$ by Lemma \ref{wPNq=>wPRNp}). The corollary follows. \ep

We can use a variant of the previous result to localize the strong Neumann problem.

\begin{corollary}[Localization of the strong Neumann problem] \label{PNq=>Lp}
Let $\Omega \subset \R^n$ be a 1-sided CAD, $L=-\diver A \nabla$ be a uniformly elliptic operator and $p\in (1,\infty)$. If \PNq$_{L^*}$ holds, then there exists $C>0$ for any ball $B=B(x,r)$ centered on the boundary, any $h,\vec F \in L^\infty_c(2B \cap \Omega)$, and any local solution $u \in W^{1,2}(2B\cap \Omega)$ to $Lu = h -\diver \vec F$ in $2B\cap\om$ with zero Neumann data on $2B\cap\pom$, we have
\[
\|\wt \N(|\nabla u|\1_B)\|_{L^p(\partial \Omega)} \leq C\left( \|\wt \C_1([|\delta h| + |\vec F|] \1_{2B})\|_{L^p(\partial \Omega)} + r^{(n-1)/p} \fiint_{(2B\setminus B)\cap\om} |\nabla u|\, dX \right).
\]
\end{corollary}

\bp We define 
\[
\widehat h :=
\begin{cases}
    h \qquad \text{in }2B\cap\om,\\
    -\frac{1}{\abs{W_B}}\iint_{2B\cap\om}h(Z)dZ\qquad\text{in }W_B,\\
    0 \qquad\text{elsewhere},
\end{cases}
\]
where $W_B\subset \gamma(x)\setminus2B$ is a Corkscrew ball associated to $B\cap\pom$. It is easy to check that the function $\widehat h$ satisfies the following properties:
 \[\widehat h \equiv h \text{ in } 2B \cap \Omega;\]
 \[\fiint_\Omega \widehat h(Z) \, dZ = 0;\]
\[\|\wt \C_1(\delta \widehat h)\|_{p} \lesssim \|\wt \C_1(\delta h\1_{2B})\|_{p}.\]
Take $u_0$ to be the weak solution to \(Lu_0 = \widehat h -\diver(\vec F \1_{2B})\) in $\om$ 
with 0 Neumann data and $\fiint_{\om}u_0\,dZ=0$. Note that $L(u-u_0)=0$ in $2B\cap\om$ and that $u-u_0$ has 0 Neumman data on $2B\cap\pom$. 

\medskip

Recall that \PNq$_{L^*}$ implies \wPNq$_{L^*}$. So we can mimic the proof of Corollary \ref{wPNq=>Lp2} until we reach 
\[\|\wt \N(|\nabla u|\1_B)\|_{L^p(\partial \Omega)} \lesssim r^{(n-1)/p} \fiint_{(2B\setminus B) \cap \Omega} |\nabla u|\,  dX + \|\wt \N(\nabla u_0)\|_{L^p(\partial \Omega)}\]
and we use \PRNp$_L$ (which is equivalent to \PNq$_{L^*}$ by Lemma \ref{PNq<=>PRNp}) to say that 
\[\|\wt \N(\nabla u_0)\|_{L^p(\partial \Omega)} \lesssim \|\wt \C_1(\delta \widehat h)\|_{p} + \|\wt \C_1(\vec F\1_{2B})\|_{p} \lesssim \|\wt \C_1([\delta |h| + |\vec F|]\1_{2B})\|_{p},\]
which completes the proof.
\ep

\section{The localization property implies Poisson Neumann problem} \label{S5}
In this section, we study the localization property \Lp. 
We start with showing that the localization property self improves.
\begin{lemma} \label{Lp=>Lp+e}
    Let $\Omega \subset \R^n$ be a bounded 1-sided CAD, $L=-\diver A \nabla$ be a uniformly elliptic operator. If \Lp$_L$ holds for some $p>1$, then there exists $\epsilon>0$ such that \Lwp{q}$_L$ holds for any $q\in [1,p+\epsilon)$.
\end{lemma}

\bp 
First, we observe that $\wt \N(|\nabla u|\1_B)$ is supported in $5B \cap \partial \Omega$, so for any $q \in [1,p]$ the H\"older inequality implies that
\[ \left( \fint_{5B \cap \partial \Omega} |\wt \N(|\nabla u|\1_B)|^q \, dX \right)^\frac1q \leq \left( \fint_{5B \cap \partial \Omega} |\wt \N(|\nabla u|\1_B)|^p \, dZ \right)^\frac1p \lesssim \left( \fiint_{\Omega \cap 2B} |\nabla u|^2 \, dX \right)^\frac12\]
by \eqref{loc1}, that is \Lwp{q} holds for $q\in[1,p]$.

\medskip

It remains to prove that \Lp\ implies \Lwp{q} for some $q>p$. First, similar to Step 1 in the proof of Lemma \ref{wPNq=>Lp}, we note that to prove for a ball $B_0=B(x_0,r_0)$ centered at $\pom$, that 
\begin{equation}\label{eq.Locp+eps}
    \norm{\wt\N(\abs{\nabla u}\1_{B_0})}_{L^{p+\epsilon}}\lesssim r_0^{\frac{n-1}{p+\epsilon}}\br{\fiint_{2B_0\cap\om}\abs{\nabla u}^2dX}^{1/2}
\end{equation}
holds for any local solution $u\in W^{1,2}(2B_0\cap\om)$ with zero Neumann data on $2B_0\cap\pom$, it suffices to show that there exists $\Lambda\ge 1$ such that for any ball $B=B(x,r)$ centered at the boundary and any local solution $u\in W^{1,2}(3\Lambda B\cap \Omega)$ to $Lu=0$ in $3\Lambda B\cap\om$ with 
zero Neumann data on $3\Lambda B\cap\pom$, we have
\begin{equation} \label{loc1b}
\|\wt \N(|\nabla u|\1_{\frac1{5\Lambda}B})\|_{L^{p+\epsilon}(\partial \Omega)} \leq C r^{(n-1)/(p+\epsilon)}\left ( \fiint_{2\Lambda B} |\nabla u|^2 dX \right)^\frac12. 
\end{equation}
Indeed, we can go from \eqref{loc1b} to \eqref{eq.Locp+eps} by a standard argument that consists of covering $B_0$ by a finite collection of balls $B''$ such that $10\Lambda^2 B'' \subset 2B_0$.

Therefore, it reduces to proving \eqref{loc1b} for a fixed ball $B=B(x,r)$ centered at the boundary. 
We take $\Lambda \ge 40 K$, for some $K\ge1$ to be determined later, and define 
\[
g(y):=\wt \N(|\nabla u|\1_{\Lambda B})(y)\quad \text{for }y\in\pom.
\]
Our goal is to show that for any ball $B'$ centered at $\pom$ that satisfies $20KB'\subset 2B$, there exists a constant $C$ independent of $B$, $B'$ and $u$, such that 
\begin{equation} \label{loc10}
    \fint_{B'\cap \partial \Omega} g^p \, d\sigma \leq C \left(\fint_{10KB'\cap \partial \Omega} g \, d\sigma \right)^p.
\end{equation}
Once we have \eqref{loc10}, we can invoke \cite[Proposition 1.1, p.122]{Giaquinta} to find some $\epsilon>0$ such that 
\begin{equation}\label{eq.gp+eps}
    \left( \fint_{\frac1{10K}B\cap \partial \Omega} g^{p+\epsilon} \, d\sigma \right)^\frac1{p+\epsilon}\lesssim \left(\fint_{2B\cap \partial \Omega} g^p \, d\sigma \right)^{p}.
\end{equation}
Note that although \cite[Proposition 1.1, p.122]{Giaquinta} is stated for cubes in the Euclidean space, the proof can be carried over to any doubling space with only minor modifications. Therefore it applies to our setting as $\pom$ is $(n-1)$-Ahlfors regular. 

We claim that \eqref{eq.gp+eps} gives the desired estimate \eqref{loc1b}. In fact, as we have chosen $\Lambda\ge 40K$, \eqref{eq.gp+eps} implies that
\begin{multline*}
    r^{-\frac{n-1}{p+\epsilon}}\norm{\wt\N(\abs{\nabla u}\1_{\frac{1}{5\Lambda}B})}_{L^{p+\epsilon}(\pom)}\lesssim
    \br{\fint_{\frac1{\Lambda}B\cap\pom}\wt\N(\abs{\nabla u}\1_{\Lambda B})^{p+\epsilon}d\sigma}^{1/(p+\epsilon)}\\
    \lesssim\br{\fint_{2B\cap\pom}\wt\N(\abs{\nabla u}\1_{\Lambda B})^pd\sigma}^{1/p}
    \lesssim r^{-(n-1)/p}\norm{\wt\N(\abs{\nabla u}\1_{\Lambda B})}_{L^p(\pom)}\\
    \lesssim \br{\fiint_{2\Lambda B}\abs{\nabla u}^2dX}^{1/2}
\end{multline*}
by \Lp.  From this \eqref{loc1b} follows.

\medskip
Now we return to prove \eqref{loc10}. Fix a $B'=B(\xi,s)$ with $\xi\in\pom$ and $20KB'\subset 2B$. 
For $y\in B'\cap \partial \Omega$ we decompose
\[ |g(y)| \leq \wt \N(|\nabla u| \1_{2B'})(y) + \wt \N(|\nabla u| \1_{\Lambda B\setminus 2B'})(y) =: g_1(y) + g_2(y).\]
For $g_1$, we use \Lp, Caccioppoli's inequality, and Moser estimate (Proposition \ref{Moser}) in that order to obtain that  \begin{multline*}
    \norm{g_1}_{L^p(\pom)}=\norm{\wt \N(|\nabla u| \1_{2B'})}_{L^p(\pom)}
    \lesssim s^{\frac{n-1}{p}}\br{\fiint_{4B'\cap\om}\abs{\nabla u}^2dX}^{/12}\\ \lesssim s^{\frac{n-1}{p}-1}\br{\fiint_{\frac92B'\cap\om}\abs{u-c_u}^2dX}^{1/2}
    \lesssim s^{\frac{n-1}{p}-1}\fiint_{5B'\cap\om}\abs{u-c_u}dX,
\end{multline*}
where $c_u$ is any constant. Take $c_u=\fiint_{5B'\cap\om}u$. Applying Poincar\'e inequality (Proposition \ref{prop.Poincare}) to the last integral above entails that there exists $K\ge 1$ such that 
\begin{equation}\label{eq.LocK}
    \norm{g_1}_{L^p(\pom)}\lesssim s^{\frac{n-1}{p}-1}\fiint_{5KB'\cap\om}\abs{\nabla u}dX.
\end{equation}
Therefore, 
\begin{equation}\label{eq.Locg1}
    \br{\fint_{B'\cap\pom}g_1^pd\sigma}^{1/p}\le s^{-\frac{n-1}{p}}\norm{g_1}_{L^p(\pom)}\lesssim \fiint_{5KB'\cap\om}\abs{\nabla u}dX.
\end{equation}
We write 
\begin{multline*}
    \fiint_{5KB'\cap\om}\abs{\nabla u}dX\lesssim (Ks)^{-n}\int_{\pom}\iint_{\gamma(y)}\abs{\nabla u(Y)}\1_{5KB'}(Y)\frac{dY}{\delta(Y)^{n-1}}d\sigma(y)\\
    \approx (Ks)^{-1}\fint_{10KB'\cap\pom}\iint_{\gamma(y)}\abs{\nabla u(Y)}\1_{5KB'}(Y)\frac{dY}{\delta(Y)^{n-1}}d\sigma(y)\\
    \lesssim \fint_{10KB'\cap\pom}\wt\N(\abs{\nabla u}\1_{5KB'})(y)\,d\sigma(y)
\end{multline*}
because the integral $\iint_{\gamma(y)}\1_{5KB'}(Y)\delta(Y)^{1-n}dY\lesssim Ks$. Plugging this estimate into the right-hand side of \eqref{eq.Locg1}, we get that 
\begin{equation}\label{eq.Locg1_p}
    \br{\fint_{B'\cap\pom}g_1^pd\sigma}^{1/p}\lesssim \fint_{10KB'\cap\pom}\wt\N(\abs{\nabla u}\1_{5KB'})(y)\,d\sigma(y)\le \fint_{10KB'\cap\pom}g(y)d\sigma(y)
\end{equation}
since $\Lambda \ge 40 K$.

For $g_2$, we see from its definition that  
\[
g_2(y)=\sup_{Y\in\gamma(y), B_Y\cap\br{\Lambda B'\setminus 2B'}\neq \emptyset}\br{\fiint_{B_Y}\abs{\nabla u}^2\1_{\Lambda B'\setminus 2B'}}^{1/2},
\]
and thus for any $y\in B'\cap\pom$, the value $g_2(y)$ is attained necessarily at a point $Y \in \gamma(y)$ such that $\delta(Y)\geq \frac{s}{2a}$, where $a>1$ is the aperture of cones in the definition of the nontangential maximal function. In fact, for any $y\in B'\cap\pom$, such $Y$ has to satisfy $\abs{Y-\xi}>3s/2$, and so $\abs{Y-y}\ge 3/2$, which entails that $\delta(Y)>\frac{s}{2a}$ as $Y\in\gamma(y)$. Then for any $z\in \frac32 B'\cap\pom$, $\abs{Y-z}\le \abs{Y-y}+\abs{y-z}\le a\,\delta(Y)+\frac52s<a^*\delta(Y)$ for some $a^*>a$. This means that $Y\in\gamma^*(z)$ for any $z\in \frac32 B'\cap\pom$, where $\gamma^*$ are cones with aperture $a^*$. Therefore, 
\[
    g_2(y) \lesssim \fint_{\frac32B' \cap \partial \Omega} \wt\N_*(\abs{\nabla u}\1_{\Lambda B})(z)\, d\sigma(z) 
    \qquad \text{ for } y \in B' \cap \partial \Omega,\]
where $\wt\N_*$ is the nontangential maximal function with aperture $a^*$. Hence, there is some $K'\ge1$ so that 
\[
g_2(y)\lesssim \fint_{2K'B' \cap \partial \Omega} \wt\N(\abs{\nabla u}\1_{\Lambda B})(z)\, d\sigma(z) 
    \qquad \text{ for } y \in B' \cap \partial \Omega.
\]
If $K'$ is greater than the constant $K$ from \eqref{eq.LocK}, we take $K=K'$. So we can take average in $y\in B'\cap\pom$ and get that 
\begin{equation*}
    \fint_{B' \cap \partial \Omega} g_2(y)^p \, d\sigma(y) \lesssim \left( \fint_{2KB' \cap \partial \Omega} g(z)\, dz \right)^p.
\end{equation*}
The desired estimate \eqref{loc10} follows from this estimate and \eqref{eq.Locg1_p}. 
\ep

\begin{lemma} \label{lemBp=>PNp'}
Let $\Omega \subset \R^n$ be a bounded 1-sided CAD, $L$ be an uniformly elliptic operator and $p\in (1,\infty)$. The property \Lp$_{L^*}$ implies \wPNq$_L$.
\end{lemma}

\bp This proof is inspired from the proof of Theorem 1.11 in \cite{MPT} and Lemma 2.10 in \cite{KP95}. Let $\vec F \in L^\infty_c(\Omega,\R^n)$, and let $u$ be the solution to $Lu=-\div \vec F$ with zero Neumann data, which is defined for $X\in \Omega$ as 
\[u(X)=\iint_{\om}\vec F(Y)\cdot\nabla_Y N(X,Y)dY.\] 

We will show the following pointwise estimate on $\wt \N(u)(\xi)$: 
\begin{equation}\label{eq.Nu_ptw}
    \wt \N(\delta \nabla u)(\xi)\lesssim \mathcal M\br{\wt \C_1(\delta \vec F)^{p_1'}}(\xi)^{1/{p_1'}}+ \wt \A^*_1(\delta\vec F)(\xi) \qquad\text{for }\xi\in\pom,
\end{equation}
where $p_1>p$ is close enough of $p$ so that \Lwp{p_1}$_{L^*}$ is true (see Lemma \ref{Lp=>Lp+e}), and $\wt \A^*_1$ is defined as $\wt \A_1$ but with a cone of larger aperture. The claim \eqref{eq.Nu_ptw} yields the desired estimate 
\[
\br{\int_{\pom}\wt \N(u)^{p'}d\sigma}^{1/{p'}}\lesssim \norm{\wt \C_1(\delta \vec F)}_{L^{p'}(\pom)}
\]
because, when $p_1>p>1$, one has the $L^{p'/p'_1}$-boundedness of the Hardy-Littlewood maximal function $\mathcal M$ and the equivalence $\norm{\wt \A_1^*(\delta\vec F)}_{L^{p'}(\pom)} \approx \norm{\wt \A_1(\delta\vec F)}_{L^{p'}(\pom)} \lesssim \norm{\wt \C_1(\delta \vec F)}_{L^{p'}(\pom)}$, see Lemma \ref{lemNCA}. 

\medskip

To prove \eqref{eq.Nu_ptw}, we fix any $\xi \in \partial \Omega$ and then $X \in \gamma(\xi)$. For $Z\in \Omega$, we split 
\[u(Z) = u_1(Z) + u_2(Z) + u_3(Z),\]
where 
\[
u_1(Z):= \iint_{2B_X}\vec F(Y)\cdot\nabla_Y N(Z,Y)dY,
\]
\[
u_2(Z):= \iint_{\gamma^*_X(\xi)\setminus 2B_X}\vec F(Y)\cdot\nabla_Y N(Z,Y)dY,
\]
where 
\[\gamma^*_X(\xi) = \{Z\in \Omega, \, \delta(X)/C^* < |Z-\xi| < C^*\delta(Z)\}
\] 
is a cone with a large aperture $C^*$ be a determined later (the same as the one used to define $\wt \A_1$ in \eqref{eq.Nu_ptw}) removing a small ball centered at $\xi$, and
\[
u_3(Z):= \iint_{\Omega \setminus \gamma^*_X(\xi)}\vec F(Y)\cdot\nabla_Y N(Z,Y)dY.
\]

We start with $u_1$, which is the solution to the Poisson-Neumann problem $-\diver(A\nabla w)=-\diver(\vec F\1_{2B_X})$. By the Lax-Milgram theorem, we have
\begin{multline*}
\left( \fiint_{B_X/4} |\delta \nabla u_1|^2 \, dZ\right)^\frac12  \lesssim \delta(X)^{1-\frac n2} \left(\iint_{\Omega} |\nabla u_1|^2 \, dZ\right)^\frac12 \lesssim \delta(X)^{1-\frac n2} \left(\iint_{\Omega} |\vec F \1_{2B_X}|^2 \, dZ\right)^\frac12 \\\ \lesssim \left(\fiint_{2B_X} |\delta \vec F|^2 \, dZ\right)^\frac12.
\end{multline*}
Note that $2B_X \subset \frac52B_Y$ for any $Y\in B_X/8$, so we have
\begin{multline*}
\left(\fiint_{2B_X} |\delta \vec F|^2 \, dZ\right)^\frac12 
=\fiint_{B_X/8}\br{\fiint_{2B_X}\abs{\delta\vec F}^2dZ}^{\frac12}dY
\lesssim \fiint_{B_X/8} \left( \fiint_{\frac52B_Y} |\delta \vec F|^2 dZ \right)^{\frac12} \, dY \\
\lesssim \delta(X)^{1-n} \iint_{B(\xi,2\delta(X)) \cap \Omega} \left( \fiint_{\frac52B_Y} |\delta \vec F|^2 dZ \right)^{\frac12} \, \frac{dY}{\delta(Y)} \lesssim \wt \C_1(\delta \vec F)(\xi).
\end{multline*}
Altogether, we have
\begin{equation} \label{eq.cclu1}
    \left( \fiint_{B_X/4} |\delta \nabla u_1|^2 \, dZ\right)^\frac12 \lesssim \wt \C_1(\delta \vec F)(\xi)
\end{equation}
as desired for \eqref{eq.Nu_ptw}.

\medskip

Let us turn to the bound on $u_2$. Since $L u_2 = 0$ in $2B_X$, we can invoke the Cacciopoli inequality to obtain that
\[\left( \fiint_{B_X/4} |\delta \nabla u_2|^2 \, dZ\right)^\frac12 \lesssim \left( \fiint_{B_X/2} |u_2|^2 \, dZ\right)^\frac12 \lesssim \sup_{Z \in B_X/2} |u_2(Z)|.\]
Using the expression of $u_2$, we continue as
\begin{multline*}
    \sup_{Z \in B_X/2} |u_2(Z)| = \sup_{Z \in B_{X}/2} \left| \iint_{\gamma^*_X(\xi) \setminus 2B_X} \vec F(Y) \cdot \nabla_Y N(Z,Y) \, dY \right| \\
 \leq \sup_{Z \in B_{X}/2} \iint_{\gamma^*_X(\xi)\setminus 2B_X}\fiint_{\frac14B_Y}\abs{\vec F(\Theta)\cdot \nabla N_{\Theta}(Z,\Theta)}d\Theta\, dY \\
 \leq \sup_{Z \in B_{X}/2} \iint_{\gamma^*_X(\xi)\setminus 2B_X}\br{\fiint_{\frac14B_Y}\abs{\vec F}^2}^{1/2}\br{\fiint_{\frac14 B_Y}\abs{\nabla_{\Theta} N(Z,\Theta)}^2d\Theta}^{1/2} dY.
\end{multline*}
We claim that for $Y\in \gamma^*_X(\xi) \setminus 2B_X$, $B_Z\subset 2a(C^*)^2B_Y\setminus \frac12 B_Y$ for all $Z\in\frac12 B_X$. This ensures that we can apply \eqref{ptwbdonN} to obtain that 
\[
    \br{\fiint_{\frac14 B_Y}\abs{\nabla_{\Theta} N(Z,\Theta)}^2d\Theta}^{1/2}   \lesssim \delta(Y)^{1-n} \qquad \text{for }Y\in \gamma^*_X(\xi) \setminus 2B_X, \, Z\in B_X/2.
\]
To verify the claim on the geometric relation, it is enough to check that $B_X\subset 2a(C^*)^2B_Y\setminus \frac12 B_Y$ for $Y\in\gamma_X^*(\xi)$. We have $\abs{X-Y}\le\abs{Y-\xi}+\abs{X-\xi}\le C^*\delta(Y)+a\delta(X)$, which is obviously less than $(C^*+2a)\delta(Y)$ if $\delta(X)<2\delta(Y)$, and is less than $(C^*+a(C^*)^2)\delta(Y)$ if $\delta(X)\ge 2\delta(Y)$, because we still have $\abs{Y-\xi}\ge \delta(X)/C^*$ and so $\delta(X)\le C^*\abs{Y-\xi}\le (C^*)^2\delta(Y)$. This implies that $B_X\subset 2a(C^*)^2B_Y$. To see that $B_X\cap \frac12 B_Y=\emptyset$, we observe that for any $Y'\in \frac12 B_Y$, $\abs{Y'-X}\ge \abs{Y-X}-\delta(Y)/8$, which is greater than $\delta(X)/4$ if $\delta(Y)<2\delta(X)$, because $Y\notin 2B_X$. If $\delta(Y)\ge 2\delta(X)$, then we use $\abs{Y-X}\ge\abs{Y-\xi}-\abs{X-\xi}\ge \delta(Y)-a\delta(X)$ to get that $\abs{Y'-X}\ge \frac78\delta(Y)-a\delta(X)\ge \br{\frac74-a}\delta(X)$, which is greater than $\delta(X)/4$ if we require the aperture $a<3/2$.

Altogether, we get 
\begin{equation} \label{eq.cclu2}
    \left( \fiint_{B_X/4} |\delta \nabla u_2|^2 \, dZ\right)^\frac12 \lesssim \iint_{\gamma^*(\xi)\setminus B_X}\br{\fiint_{\frac14B_Y}\abs{\delta \vec F}^2}^{1/2} \frac{dY}{\delta(Y)^n} \lesssim \wt \A^*_1(\delta \vec F)(\xi)
\end{equation}
as desired for \eqref{eq.Nu_ptw}.

It remains to treat the most complicated term $u_3$. We construct a Whitney type decomposition of $\Omega \setminus \gamma^*_X(\xi)$. Indeed, if the aperture $C^*$ of $\gamma^*_X(\xi)$ is large enough ($C^* \geq 301$ works), the collection of balls $\{B(y,|y-X|/100)\}_{y\in \partial \Omega}$ covers $\Omega \setminus \gamma^*_X(\xi)$. By Vitali's covering lemma, there exists a non overlapping subcollection $\{B(y_m,|y_m-X|/100)\}_{m\in I}$ such that $\{B_m:= B(y_m,|y_m-X|/20)\}_{m\in I}$ is a finitely overlapping cover of $\Omega \setminus \gamma^*_X(\xi)$. 

\noindent Note that the radius of $B_m$ is bounded from below by $\delta(X)/20$, and there can be only a finite number of balls with same radius, so we can order $B_m$ from small radius to large radius (we identify $I$ to a subset of $\bN$), and we define 
\[D_m := \big(B_m \cap \Omega\big) \setminus \left( \gamma^*_X(\xi) \cup \bigcup_{m'<m} B_{m'} \right),\]
so that $D_m$ forms a partition of $\Omega \setminus \gamma^*_X(\xi)$. Moreover, since the number of balls in $\{B_m\}_{m\in I}$ of radius between $2^{k}\delta(X)/20$ and $2^{k+1}\delta(X)/20$ is uniformly bounded in $k\in\NN$, we have
\begin{equation} \label{sumrm}
\sum_{m\in I} \left( \frac{\delta(X)}{|y_m-X|}\right)^\beta  \lesssim 1 \qquad \text{for any }\beta>0. 
\end{equation}
Note also that 
\begin{equation}\label{eq.BmBX}
    10 B_m\cap B_X=\emptyset \qquad\text{for }m\in I
\end{equation}
as we can quickly check that for any $Z\in10 B_m$ and any $X'\in B_X$, $\abs{Z-X'}\ge \delta(X)/4>0$.
\noindent With all this preparation, we have the bound
\begin{equation} \label{u3byu3m}
    \left( \fiint_{B_X/4} |\delta \nabla u_3|^2 \, dZ \right)^\frac12 \lesssim \sum_{m\in I} \left( \fiint_{B_X/4} |\delta \nabla u_{3,m}|^2 \, dZ \right)^\frac12,
\end{equation}
where 
\[u_{3,m}(Z):= \iint_{D_m} \vec F(Y) \cdot \nabla_Y N(Z,Y) \, dY.\]

Now, we need to deal with each $u_{3,m}$, which is a solution to $Lw=-\diver\br{\vec F \1_{D_m}}$ with zero Neumann data. Using the Cacciopoli inequality and then H\"older continuity of the solution $u_{3,m}$ (see Proposition \ref{Holder}), we have
\begin{equation*}
    \left( \fiint_{B_X/4} |\delta \nabla u_{3,m}|^2 \, dZ \right)^\frac12 \lesssim \osc_{B_X/2} u_{3,m} 
    \lesssim \left( \frac{\delta(X)}{|y_m-X|} \right)^\alpha \osc_{B'_m\cap \om} u_{3,m}
    \leq \left( \frac{\delta(X)}{|y_m-X|} \right)^\alpha \sup_{B'_m\cap \om} |u_{3,m}|
\end{equation*}
where $B'_m \supset B_X$ is the largest ball centered $X$ that doesn't intersect $4B_m$. Note that $B'_m$ is well-defined because of \eqref{eq.BmBX}. Using the expression of $u_{3,m}$ and the Carleson inequality \eqref{eq.CNdual}, we have
\begin{multline*}
\sup_{B'_m\cap\om} |u_{3,m}| \leq \sup_{Z\in B'_m\cap\om} \iint_{D_m} |\vec F| |\nabla_Y N(Z,Y)| \, dY\\
\lesssim \sup_{Z\in B'_m\cap\om} \|\wt \C_1(\delta \vec F \1_{B_m})\|_{L^{p'_1}(\pom)} \|\wt \N(\nabla N(Z,.)\1_{B_m})\|_{L^{p_1}(\pom)}.
\end{multline*}
On one hand, since $p_1$ is chosen so that \Lwp{p_1}$_{L^*}$ is true, and since $Y \mapsto N(Z,Y)$ is a solution to $L^*u = 0$ in $3B_m$ (recall that $Z\in B_m'\cap\om$ and that $B_m'\cap 4B_m=\emptyset$), we have
\begin{multline*} \|\wt \N(\nabla N(Z,.)\1_{B_m})\|_{p_1} \lesssim |y_m-X|^{\frac{n-1}{p_1}} \left( \fiint_{2B_m\cap\om }|\nabla_Y N(Z,Y)|^2 \, dY \right)^\frac12 \lesssim |y_m-X|^{-\frac{n-1}{p_1'}}
\end{multline*}
by \eqref{ptwbdonN}. On the other hand,
\[
     \|\wt \C_1(\delta \vec F \1_{B_m})\|_{p'_1}^{p_1'}
    \lesssim   \int_{2B_m \cap \partial \Omega} \wt \C_1(\delta \vec F)^{p_1'} \, d\sigma 
    +\sum_{k=1}^\infty\int_{\br{2^{k+1}B_m\setminus2^kB_m}\cap\pom}\wt \C_1(\delta \vec F\1_{B_m})^{p_1'} \, d\sigma  
\]
For any $y\in \br{2^{k+1}B_m\setminus2^kB_m}\cap\pom$, observe that
\begin{multline*}
    \wt \C_1(\delta \vec F\1_{B_m})(y)
    \lesssim\sup_{r>2^k\abs{y_m-X}}r^{1-n}\iint_{B(y,r)\cap\om}\br{\fiint_{B_Y}\abs{\vec F}^2\1_{B_m}}^{1/2}dY\\
    =2^{-k(n-1)}\abs{y_m-X}^{1-n}\iint_{2B_m\cap\om}\br{\fiint_{B_Y}\abs{\vec F}^2\1_{B_m}}^{1/2}dY\\
    \lesssim 2^{-k(n-1)}\wt\C_1(\delta\vec F\1_{B_m})(y') \qquad \text{for all }y'\in B_m\cap\pom,
\end{multline*}
and thus,
\begin{multline*}
\sum_{k=1}^\infty\int_{\br{2^{k+1}B_m\setminus2^kB_m}\cap\pom}\wt \C_1(\delta \vec F\1_{B_m})^{p_1'} \, d\sigma  \lesssim\sum_{k=1}^\infty 2^{-k(n-1)(p_1'-1)}\br{\fint_{B_m\cap\pom}\wt\C_1(\delta\vec F\1_{B_m})d\sigma}^{p_1'}\\
\lesssim \int_{B_m\cap\pom}\wt\C_1(\delta\vec F)^{p_1'}d\sigma \qquad\text{since }p_1'>1.
\end{multline*}
Then 
\begin{multline*}
    |y_m-X|^{-\frac{n-1}{p_1'}}\|\wt \C_1(\delta \vec F \1_{B_m})\|_{p'_1}\lesssim\br{\fint_{2B_m\cap\pom}\wt\C_1(\delta\vec F)^{p_1'}d\sigma}^{1/{p_1'}}\\
    \lesssim \br{\fint_{B(\xi,10\abs{y_m-X})\cap\pom}\wt\C_1(\delta\vec F)^{p_1'}d\sigma}^{1/{p_1'}}
    \lesssim \Big|\mathcal M\Big(\wt \C_1(\delta \vec F)^{p'_1}\big)(\xi)\Big|^{\frac1{p'_1}}.
\end{multline*}

Altogether, we have for some $\alpha \in (0,1)$
\begin{equation} \label{eq.cclu3m}
    \left( \fiint_{B_X/4} |\delta \nabla u_{3,m}|^2 \, dZ\right)^\frac12 \lesssim \left(\frac{\delta(X)}{|y_m-X|}\right)^\alpha\Big|\mathcal M\Big(\wt \C_1(\delta \vec F)^{p'_1}\big)(\xi)\Big|^{\frac1{p'_1}} \qquad\text{for }m\in I,
\end{equation}
that is
\begin{equation} \label{eq.cclu3}
    \left( \fiint_{B_X/4} |\delta \nabla u_{3}|^2 \, dZ\right)^\frac12 \lesssim \Big|\mathcal M\Big(\wt \C_1(\delta \vec F)^{p'_1}\big)(\xi)\Big|^{\frac1{p'_1}},
\end{equation}
by \eqref{u3byu3m} and \eqref{sumrm}.

\medskip

To conclude,  we have chosen $\xi \in \partial \Omega$, $X\in \gamma(\xi)$, decomposed $u = u_1+u_2+u_3$, and then the bounds \eqref{eq.cclu1}, \eqref{eq.cclu2} and \eqref{eq.cclu3} on $u_i$, $i=1,2,3$ respectively give that
\[\left(\fiint_{B_X/4} |\delta \nabla u|^2 \, dZ \right)^\frac12 \lesssim \mathcal M\br{\wt \C_1(\delta \vec F)^{p_1'}}(\xi)^{1/{p_1'}}+ \wt \A^*_1(\delta\vec F)(\xi),\]
which in turn gives the desired bound \eqref{eq.Nu_ptw}. The theorem follows.
\ep

\section{The strong Neumann problem is solvable in $L^p$ for $p$ in an open interval} \label{S6}

The goal of this section is twofold. First, we want to prove that \PRNp $\implies$ \PRNwp{q} for $q\in (1,p)$, and - when $p>q$ - that \wPNq$_{L^*}$ $+$ \PRNwp{q}$_L$ $\implies$ \PRNp$_L$. A combination of the two (with results from the previous sections) gives the fact that \PRNp $\implies$ \PRNwp{q} for $q\in (1,p+\epsilon)$.

\begin{lemma} \label{PRNp=>PRNq}
Let $\Omega \subset \R^n$ be a bounded 1-sided CAD domain, $L = - \diver A \nabla$ be a uniformly elliptic operator, and $p\in (1,\infty)$.
\[\text{\PRNp $\implies$ \PRNwp{q} for $q\in (1,p)$}.\]
\end{lemma}

\bp Our goal is to use the interpolation result given as Theorem \ref{Thinterpolation}. We construct, for $f\in L^\infty_c(\Omega)$, $[f/\delta]_\Omega = 0$, the map
\[\mathcal Z(f) := \nabla u_f\]
where $u_f$ is the function constructed as
\[u_f(X):=\iint_\Omega N(X,Y) \frac{f(Y)}{\delta(Y)}\, dY,\] and we want to prove \eqref{interpolation1}.

\medskip

Take $B=B(x,r)$, $j$ and $f$ as in Theorem \ref{Thinterpolation}, in particular $h:= f/\delta$ has zero average. We write $C_j(B)$ for $2^{j+1}B \setminus 2^j B$ and $C^*_j(B)$ for $2^{j+2}B \setminus 2^{j-1} B$. We first claim that 
\begin{equation} \label{interpolation2}
    (2^jr)^{(1-n)/p} \|\wt \N[\nabla u_f \1_{C_j(B)}]\|_p \lesssim \osc_{C_j^*(B) \cap \Omega} u_f \qquad\text{for }j\ge5.
\end{equation}
This is pretty straightforward. Fix any integer $j\ge5$, like in Step 1 of the proof of Lemma \ref{wPNq=>Lp}, we cover $C_j(B)$  
by a finite (the number is independent of $B$ and $j$) collection of ball $\set{B'_i}_{i=1}^N$ with radius $r' \approx 2^jr$ such that 
\begin{enumerate}[(i)]
    \item either $B_i'$ is centered at the boundary and $3B'_i \subset C_j^*B$,
    \item or $3B'_i \subset C_j^*(B) \cap \Omega$.
\end{enumerate}
If the ball $B'_i$ satisfy (i), Lemma \ref{wPNq=>Lp} (which is satisfied since \PRNp$_L$ implies \wPNq$_{L^*}$) entails that
\[ (2^jr)^{(1-n)/p} \|\wt \N[\nabla u_f \1_{B'_i}]\|_p \lesssim \fiint_{2B'_i\cap\om} |\nabla u_f| \, dX 
\lesssim (2^jr)^{-1} \osc_{3B'_i\cap\om} u_f  
\lesssim  (2^jr)^{-1} \osc_{C_j^*(B) \cap \Omega} u_f,\]
where the second inequality is due to H\"older's inequality  and Caccioppoli inequality. If the ball $B'_i$ satisfies (ii), we directly have
\[ (2^jr)^{(1-n)/p} \|\wt \N[\nabla u_f \1_{B'_i}]\|_p \lesssim \left(\fiint_{B'_i} |\nabla u_f|^2 \, dX\right)^\frac12 \lesssim (2^jr)^{-1} \osc_{3B'_i} u_f  \lesssim  (2^jr)^{-1} \osc_{C_j^*(B) \cap \Omega} u_f.\]
Altogether, we have
\[r^{(1-n)/p}\|\wt \N[\nabla u_f \1_{C_j(B)}]\|_p \lesssim \sup_{i} \|\wt \N[\nabla u_f \1_{B'_i}]\|_p \lesssim (2^jr)^{-1} \osc_{C_j^*(B) \cap \Omega} u_f \]
as desired for the claim \eqref{interpolation2}.

\medskip

Since $h=f/\delta$ has zero average and is supported in $B\cap \Omega$, for any $X \in C^*_j(B)$ 
\begin{multline*}
|u_f(X)| =\abs{\iint_{B\cap\om} (N(X,Y)-N(X,x))\frac{f(Y)}{\delta(Y)}dY }
\leq \Big(\osc_{Y\in B \cap \Omega} N(X,Y) \Big) \left(\iint_\Omega \Big|\frac{f}\delta\Big| \, dY\right) \\ \lesssim 2^{-j\alpha} \Big(\osc_{Y\in 2^{j-4}B \cap \Omega} N(X,Y) \Big) \left(\int_{\partial \Omega} \wt \A_1(f) \, d\sigma\right)
\lesssim 2^{-j\alpha} (2^{j}r)^{2-n} \left(\int_{\partial \Omega} \wt \A_1(f) \, d\sigma\right)
\end{multline*}
for some $\alpha>0$ by H\"older continuity and \eqref{ptwbdonN2}. 

Altogether, we proved that for any integer $j\ge 5$,
\[(2^jr)^{(1-n)/p}\|\wt \N[\nabla u_f \1_{B'_i}]\|_p \lesssim 2^{j(1-n-\alpha)} r^{1-n} \int_{\partial \Omega} \wt \A_1(f) \, d\sigma,\]
which is the bound \eqref{interpolation1} with $g(j) := 2^{j(1-n-\alpha)}$. For any $q\in (1,p)$, Theorem \ref{Thinterpolation} and Lemma \ref{lemNCA} give that for any $f \in L^\infty_c(\Omega)$ with $[f/\delta]_\Omega = 0$, there holds
\[\|\wt \N[\mathcal Z(f)]\|_{L^q(\partial \Omega)} \lesssim \|\wt \A_1 (f)\|_{L^q(\partial \Omega)}\approx \norm{\wt\C_1(f)}_{L^q(\pom)},\]
 or equivalently,
\[\|\wt \N(\nabla u)\|_{L^q(\partial \Omega)} \lesssim \|\wt \C_1 (\delta h)\|_{L^q(\partial \Omega)}\]
whenever $h\in \widehat{L}^\infty_c(\Omega)$ and $u(X) := \iint_\Omega N(X,Y) h(Y) \, dY$. We conclude that \PRNwp{q}$_L$ is solvable for any $q\in(1,p)$, as desired.
\ep

\begin{lemma} \label{lemShen}
Let $\Omega \subset \R^n$ be a bounded 1-sided CAD, $L=-\diver A \nabla$ be a uniformly elliptic operator, and $1<q<p<r < \infty$. Then
\[ \text{\PRNwp{q}$_L$ $+$ \Lwp{r}$_{L}$ $\implies$ \PRNp$_L$.}\]
\end{lemma}

\bp The proof is close to the one of Shen in \cite{Shen07}, and follows the now classical strategy using a good-lambda argument. 

\medskip

\noindent {\bf Step 1: Good lambda argument.} 
Let $h \in \widehat{L}^\infty_c(\Omega)$, $u$ be a solution to $Lu=h$ with zero Neumann data. 

We want to estimate the set 
\[E(\lambda) := \{x\in \partial \Omega, \, \mathcal M(|\wt \N(\nabla u)|^q)(x) > \lambda\},\]
which has finite measure since  the Hardy-Littlewood maximal function $\mathcal M$ is weak $(1,1)$ and $\N(\nabla u) \in L^q(\partial \Omega, \sigma)$. 
More precisely, we want to show that there exists $C>0$ and $r>p$  (independent of $h$) such that for any $\eta,\gamma \in (0,1)$ we have 
\begin{equation} \label{Shenclaim}
|E(\lambda) \cap \{x\in \partial \Omega, \, \mathcal M (|\wt \C_1(\delta h)|^q)(x) \le \gamma \lambda\}| \leq C(\gamma + \eta^{r/q})  |E(\eta \lambda)|.
\end{equation}
\medskip

We define a Whitney decomposition\footnote{If $E(\lambda) = \partial \Omega$, it means that the Whitney decomposition is $\{\partial \Omega\}$} as $E(\eta \lambda)$ as follows: take the collection $\{\Delta(x,r_x)\}_{x\in E(\eta \lambda)}$ - where $r_x = \dist(x,\partial \Omega \setminus E(\eta\lambda))/50$ - of boundary balls that covers $E(\eta\lambda)$ and, by the Vitali lemma, take a covering subcollection $\{\Delta_k:=\Delta(x_k,r_{x_k})\}$ such that $\{\Delta_i/5\}$ is non-overlapping. 

The claim \eqref{Shenclaim} will be proved if we can show that, for any $k\geq 1$
\begin{equation} \label{Shen1}
|E(\lambda) \cap \Delta_k| \lesssim (\gamma + \eta^{r/q})  |\Delta_k|.
\end{equation}
whenever we have a $z\in \Delta_k$ such that
\begin{equation} \label{Shen2}
\{x\in \Delta_k, \, \mathcal M (|\wt \C_1(\delta h)|^q)(x) \le \gamma \lambda\} \supset \{z\} \neq \emptyset.
\end{equation}

The fact that $\Delta_k$ are Whitney balls for the set $E(\eta \lambda) \supset E(\lambda)$ means that 
\begin{equation}\label{eq.Shen_y}
    \text{there exists }y\in 100\Delta_k \quad\text{such that }\mathcal M(\wt \N(\nabla u)|^q)(y) \leq \eta \lambda,
\end{equation}
which in turn means that $\abs{\nabla u}$ cannot take large value too far from $\Delta_k$, whence, for $x\in \Delta_k \cap E(\lambda)$, 
\begin{equation} \label{Shen3}
\lambda < \mathcal M(|\wt \N(\nabla u)|^q)(x) = \mathcal M(|\wt \N(\nabla u \1_{2B_k})|^q)(x),
\end{equation}
where $B_k:= B(x_k,r_{x_k})$ and as long as $\eta$ is small enough.

Now, we construct $u_0$ to be the solution to $Lu_0 = h\1_{6B_k}$ in $\Omega$ that has zero Neumann data on $\pom$, and we have
\begin{multline*}
    |E(\lambda) \cap \Delta_k| = |\{x\in \Delta_k, \, \mathcal M(|\wt \N(\nabla u \1_{2B_k})|^q)(x) > \lambda\}| \\
    \leq |\{x\in \Delta_k, \, \mathcal M(|\wt \N(\nabla[u-u_0]\1_{2B_k})|^q)(x) > \lambda\}| + |\{x\in \Delta_k, \, \mathcal M(|\wt \N(\nabla u_0 \1_{2B_k})|^q)(x) > \lambda\}| \\ =:S_1+S_2.
\end{multline*}
We estimate $S_2$ using the weak $(1,1)$ boundedness of $\mathcal M$ and \PRNwp{q}$_L$, and we get
\begin{multline*}
    S_2 \lesssim \frac1\lambda \int_{\partial \Omega} |\wt \N(\nabla u_0 \1_{2B_k})(x)|^q d\sigma(x)
    \lesssim \frac1\lambda \|\wt \C_1(\delta h \1_{6B_k})\|_{L^q}^q 
    \approx\frac{1}{\lambda}\int_{8\Delta_k}\wt\C_1(\delta h\1_{6B_k})^qd\sigma    
    \\ \lesssim \frac{|\Delta_k|}{\lambda} \mathcal M(|\wt \C_1(\delta h \1_{6B_k})|^q)(z) \lesssim \gamma |\Delta_k|
\end{multline*}
by \eqref{Shen2} and the fact that $\|\wt \C_1(\delta h \1_{6B_k})\|_{L^q} \approx \|\wt \C_1(\delta h \1_{6B_k})\|_{L^q(8\Delta_k)}$ (see the computation in the proof of Lemma \ref{lemBp=>PNp'}). 

We turn to the bound on $S_1$. Since \Lwp{r}$_L$ holds, by the weak $(\frac rq,\frac rq)$ bound on $\mathcal M$ (we also use the improvement given in Lemma \ref{wPNq=>Lp} for the solution $u-u_0$ to $L(u-u_0)=0$ in $6B_k\cap\om$), we have
\[
     S_1 \lesssim \lambda^{-r/q} \int_{\partial \Omega} |\wt \N(\nabla (u-u_0) \1_{2B_k})(x)|^r d\sigma(x) 
       \lesssim \lambda^{-r/q} |\Delta_k| \left( \fiint_{4B_k\cap\om} |\nabla[u-u_0]| \, dX \right)^{r}.
\]
We can bound the solid integral on $|\nabla (u-u_0)|$ by an integral of the nontangential maximal function, that is,
\begin{multline*}
    \fiint_{4B_k} |\nabla(u-u_0)| \, dX
    \lesssim\int_{y\in\pom}\iint_{\gamma(y)}\abs{\nabla(u-u_0)(X)}\1_{4B_k}(X)\frac{dX}{\delta(X)^{n-1}}d\sigma(y)\\
    \le \int_{10\Delta_k}\wt\N(\nabla(u-u_0)\1_{4B_k})(y)\iint_{\gamma(y)\cap 4B_k}\frac{dX}{\delta(X)^{n-1}}d\sigma(y)
    \lesssim \fint_{10\Delta_k}\wt\N(\nabla(u-u_0)\1_{4B_k})\,d\sigma\\
    \lesssim\fint_{100\Delta_k}\wt\N(\nabla u)\,d\sigma
    +\fint_{10\Delta_k}\wt\N(\nabla u_0)\,d\sigma.
\end{multline*}
By H\"older's inequality, \PRNwp{q}, \eqref{Shen2} and \eqref{eq.Shen_y}, we obtain that  
\begin{multline*}
    \fiint_{4B_k} |\nabla(u-u_0)| \, dX
    \lesssim
    \left(\fint_{100\Delta_k} |\wt \N(\nabla u)|^q \, d\sigma\right)^\frac1q 
    + \abs{\Delta_k}^{-1/q}\norm{\wt\C_1(\delta h\1_{6B_k})}_{L^q(\pom)} \\
    \lesssim
    \br{\mathcal M(|\wt \N(\nabla u)|^q)(y)}^\frac1q  + \left(\fint_{12\Delta_k} |\wt \C_1(\delta h)|^q \, d\sigma\right)^\frac1q \\
    \lesssim (\eta \lambda)^\frac1q + \br{\mathcal M(|\wt \C_1(\delta h)|^q)(z)}^\frac1q \lesssim \lambda^\frac1q \br{\eta^{\frac1q}+\gamma^{\frac1q}}.
\end{multline*}
Altogether, 
\[S_1\lesssim (\eta^{r/q} + \gamma^{r/q}) |\Delta_k| \leq (\eta^{r/q} + \gamma) |\Delta_k|\]
and then
\[|E(\lambda)\cap \Delta_k| \lesssim S_1 + S_2 \lesssim (\eta^{r/q} + \gamma) |\Delta_k|.\]
The claims \eqref{Shen1} and then \eqref{Shenclaim} follow.

\medskip

\noindent {\bf Step 2: Conclusion.} From Step 1, we have for any $\lambda >0$
\[ |E(\lambda)| \leq C(\gamma + \eta^{r/q}) |E(\eta \lambda)| + |\{x\in \partial \Omega,\, \mathcal M(|\wt C_1(\delta h)|^q)> \gamma \lambda\}|\]
If we multiply $E(\lambda)$ by $\lambda^{p/q-1}$ and integrate, for any $\Lambda>0$, we have that 
\begin{multline} \label{Shen4}
    \int_0^\Lambda |E(\lambda)| \lambda^{\frac pq-1} d\lambda \leq C(\gamma + \eta^{\frac rq}) \int_0^\Lambda |E(\eta \lambda)| \lambda^{\frac pq-1} d\lambda \\ 
    + \int_0^\Lambda |\{x\in \partial \Omega,\, \mathcal M(|\wt C_1(\delta h)|^q)> \gamma \lambda\}| \lambda^{\frac pq-1} d\lambda \\
    \leq C(\gamma + \eta^{\frac rq})\eta^{-\frac pq} \int_0^{\eta\Lambda} |E(\lambda)| \lambda^{\frac pq-1} d\lambda + C_\gamma \int_{\partial \Omega} |\wt C_1(\delta h)|^p \, d\sigma
\end{multline}
by the $L^{p/q}$-boundedness of the Hardy-Littlewood maximal function $\mathcal M$. Note that the left-hand side above is finite, as 
\[
\int_0^\Lambda |E(\lambda)|\lambda^{\frac{p}{q}-1}d\lambda\le \Lambda^{p/q-1} \abs{\pom} < +\infty, 
\]
Now in particular, we can hide the first term of the right-hand side of \eqref{Shen4} to the left-hand side, provided that $C(\gamma + \eta^{r/q})\eta^{-p/q} \leq \frac12$. It is possible because $r>p$, and so we can first fix $\eta>0$ such that $C \eta^{(r-p)/q} \leq \frac14$, and then we fix $\gamma>0$ such that $C \gamma \eta^{-p/q} \leq \frac 14$. With those choices of $\eta$ and $\gamma$, we have
\[\int_0^\Lambda |E(\lambda)| \lambda^{\frac pq-1} d\lambda \leq C_\gamma \int_{\partial \Omega} |\wt C_1(\delta h)|^p \, d\sigma,\]
and taking $\Lambda \to \infty$ finally gives
\[\|\wt \N(|\nabla u)|)\|_{p}^p\le \int_{\pom}\mathcal{M}\br{\wt\N(\nabla u)^q}^{p/q}d\sigma=  \frac{p}{q}\int_0^\infty |E(\lambda)| \lambda^{\frac pq-1} d\lambda \lesssim \|\wt C_1(\delta h)\|_{p}^p,\]
which proves \PRNp\ and hence the lemma.
\ep

\section{Links with Neumann problem} \label{S7}

The purpose of this section is to prove the last implcation that we haven't proven yet, that is the first implcation of Proposition \ref{Np<=>PRNp}.

\begin{lemma} \label{Np+Dq=>PRNp}
 Let $\Omega \subset \R^n$ be a 1-sided CAD, $L = - \div A \nabla$ be a uniformly elliptic operator, and $p\in (1,\infty)$. Then
 \[ \text{\Np$_L$ $+$ \Dq$_{L^*}$ $\implies$ \PRNp$_L$}.\]
\end{lemma}

\bp
The proof of this implication is non-trivial, but is actually easy if we utilize the Poisson-Regularity estimates from \cite{MPT}.

So let $h \in \widehat{L^\infty_c}(\Omega)$, and let $u$ be the weak solution to $Lu=h$ in $\Omega$ with zero Neumann data. We want to prove that 
\[\|\wt \N(\nabla u)\|_{L^p(\partial \Omega)} \leq C \|\wt \C_1(\delta h)\|_{L^p(\partial \Omega)}\]
with a constant $C>0$ independent of $h$.

Take $u_D \in W^{1,2}(\Omega)$ to be the weak solution of $Lu=h$ in $\Omega$ with zero Dirichlet data instead. Morally, since $u-u_D$ is a solution to $Lw=0$ in $\Omega$, the solvability of the Neumann problem \Np$_L$ implies
\begin{equation} \label{Np1}
\|\wt \N(\nabla[u-u_D])\|_{L^p(\partial \Omega)} \lesssim \|(A\nabla u_D)_n\|_{L^p(\partial \Omega)} \lesssim \|\wt \N(\nabla u_D)\|_{L^p(\partial \Omega)} \lesssim \|\wt \C_1(\delta h)\|_{L^p(\partial \Omega)}, 
\end{equation}
where $(A\nabla u_D)_n$ is the Neumann boundary data of $u_D$, and the last inequality is by the assumption \Dq$_{L^*}$ holds and $(1)\implies (3)$ of Theorem \ref{ThMPT} (i.e. \cite[Theorem 1.22]{MPT}).

\medskip

We need to prove \eqref{Np1} rigorously by justifying that 
\begin{equation}\label{eq.u-uD=intN}
    u(X)-u_D(X)=\int_{\pom}N(X,y)g(y)\,d\sigma(y) +C
\end{equation}
for some $g\in L^p(\pom)$ and an arbitrary constant $C$. To that purpose, 
we define for $\phi\in Lip(\pom)$,
\[\ell(\phi): = \iint_\Omega A \nabla u_D\cdot \nabla \Phi\, dX - \iint_{\Omega} h \, \Phi\, dX,\]
where $\Phi$ is any extension of $\phi$ in $W^{1,2}(\Omega) \cap Lip(\Omega)$. Note that this is well-defined because  if $\Phi_1$ and $\Phi_2$ are two extensions of $\phi$ in $W^{1,2}(\Omega)$, then 
\[\iint_\Omega A \nabla u_D\cdot \nabla (\Phi_1- \Phi_2)\, dX - \iint_{\Omega} h \, (\Phi_1-\Phi_2)\, dX = 0\]
because $u_D \in W^{1,2}(\Omega)$ is a weak solution to $Lu_D = h$ and $\Phi_1 - \Phi_2 \in W_0^{1,2}(\Omega)$. Since $\phi \in H^{1/2}(\partial \Omega,\sigma) \cap L^{p'}(\partial \Omega,\sigma)$, we can conveniently use the Varopoulos extension $\Phi$ constructed in \cite[Theorem 1.4]{MZ??} (see also \cite{HR18}), which belongs to $W^{1,2}(\Omega) \cap Lip(\Omega)$ and satisfies 
\begin{equation} \label{Np2}
\|\wt \N(\Phi)\|_{p'} + \|\wt \C_1(\delta \nabla \Phi)\|_{p'} \leq C \|\phi\|_{p'}
\end{equation}
for a constant $C$ independent of $\phi$. Then by \eqref{eq.CNdual}, (3) of Theorem \ref{ThMPT}, and \eqref{Np2},  we have that 
\[
|\ell(\phi)| \lesssim \|\wt \N(\nabla u_D)\|_p \|\wt \C_1(\delta \nabla \Phi)\|_{p'} + \|\wt \C_1(\delta h)\|_p \|\wt \N(\Phi)\|_{p'} \lesssim \|\wt \C_1(\delta h)\|_p \|\phi\|_{p'}.
\]
Therefore, $\ell$ is a bounded functional on $L^{p'}(\pom)$ and so there exists $g_D \in L^p(\partial \Omega)$ with \begin{equation}\label{eq.gD}
    \norm{g_D}_{L^p(\pom)}\le C\norm {\wt\C_1(\delta h)}_{L^p(\pom)}
\end{equation}
such that 
\(\ell(\phi) = \int_{\partial \Omega} g_D\,\phi \, d\sigma.\)
Moreover, since $\iint_{\Omega} h \, dX = 0$, it is easy to see that $\ell(1)=0$ and thus $\int_{\partial \Omega} g_D \, d\sigma =0$.

We define now $v \in W^{1,2}(\Omega)$ as 
\[v(X) := \int_{\partial \Omega} N(X,y) g_D(y) \, d\sigma(y).\]
On one hand, since $g_D$ has zero average, by exchanging the order of integration and \eqref{eq.RieszNmn'},
we obtain that
\begin{multline*}
    \iint_\Omega A \nabla v \cdot \nabla \Phi \, dX 
=\int_{\pom}g_D(y)\iint_{\om}A(X)\nabla_XN(X,y)\cdot\nabla\Phi(X)dXd\sigma(y)\\
    =\int_{\pom}g_D(y)\br{\phi(y)-\fint_{\pom}\phi\,d\sigma}d\sigma(y)
    = \int_{\partial \Omega} g_D(y) \phi(y) \, d\sigma(y) = \ell(\phi)
\end{multline*}
whenever $\phi \in Lip(\partial \Omega)$ and $\Phi \in W^{1,2}(\Omega) \cap Lip(\Omega)$. 
On the other hand, $u_D-u \in W^{1,2}(\Omega)$ and satisfies
\[\iint_\Omega A \nabla[u_D-u] \cdot \nabla \Phi \, dX = \ell(\phi)- \left( \iint_\Omega A\nabla u\cdot \nabla \Phi\, dX - \iint_\Omega h\, \Phi\, dX \right) = \ell(\phi)\]
since $u$ is the solution to $Lu=h$ with zero Neumann data.
Combining the two estimates, we get that 
\[
\iint_\om A\nabla\br{u_D-u-v}\cdot\nabla\Phi\,dX=0
\]
for all $\Phi\in W^{1,2}(\om)$ by density. 
Taking $\Phi=u_D-u-v$ gives that $\abs{\nabla (u_D-u-v)}=0$ in $\om$, and thus $u_D-u-v$ is constant, which proves \eqref{eq.u-uD=intN}, with $g=g_D$.

We conclude by applying \Np\ to the solution $v$ and using \eqref{eq.gD} as follows.
\[\|\wt \N(\nabla[u-u_D])\|_{L^p(\partial \Omega)} =\|\wt \N(\nabla v)\|_{L^p(\partial \Omega)} \lesssim \|g_D\|_{L^p(\partial \Omega)} \lesssim \|\wt C_1(\delta h)\|_{L^p(\partial \Omega)}\]
as desired.\ep

\appendix

\section{Interpolation between tent spaces.} \label{SAppendix}

The purpose of this appendix is to prove ``by hand'' the real interpolation result given as Theorem \ref{Thinterpolation}.

We shall follow the strategy from \cite[Theorem 1.1]{Auscher07} to prove a weak $L^1$-bound on $\wt \N(\mathcal Z(f))$. The first tool that we require is a Calder\'on-Zygmund decomposition adapted to tent spaces.

\begin{lemma} \label{lemCZ}
Let $\Omega \subset \R^n$ be a 1-sided CAD and $q\in (0,\infty)$. There exists $C>0$ such that for any function $f\in \wt T^q_1(\Omega)$ and any level set $l>0$, we can find a function $g$, a family of balls $\{B_i=B(x_i,r_i)\}_{i\in I}$ centered at the boundary, and a collection of functions $\{b_i\}_{i\in I}$ such that 
\begin{enumerate}[(i)]
    \item $f = g + \sum_{i\in I} b_i$; \smallskip
    \item $\supp b_i \subset B_i$, and $\sum_{i\in I}\1_{B_i}\le C$; \smallskip
    \item $\displaystyle \iint_\Omega \frac{b_i}{\delta} \, dX = 0$; \smallskip
    \item $\sup_{x\in \partial \Omega} \wt \A_1(g)(x) \leq Cl$; \smallskip
    \item $\|b_i\|_{\wt T^q_1} \leq Cl r_i^{(n-1)/q}$; \smallskip
    \item $\displaystyle \sum_{i\in I} r_i^{n-1} \leq l^{-q} \|f\|_{\wt T^q_1(\Omega)}^q$.
\end{enumerate}
\end{lemma}

\bp The proof is essentially the same as the Calder\'on-Zygmund decomposition found in \cite{Huang17}. However, since our setting differs quite a lot from \cite{Huang17} - which  treats $T^q_2(\R^n_+)$ and does not show $(iii)$ - we decided to present again the proof for completeness.

\medskip

\noindent {\bf Construction of $g,b_i$.} Take $f \in \wt T^q_1(\Omega)$ and $l>0$. Define 
\[E_l := \{x\in \partial \Omega, \, \mathcal M(\wt\A_1^*(f)^q)(x)> l^q\},\]
where $\mathcal M$ is the uncentered maximal function on $\partial \Omega$ with the $(n-1)$-Ahlfors-regular - hence doubling - measure $\sigma$ and $\wt \A_1^*$ is defined with cones of large aperture $a^*$. If $E_l = \partial \Omega$, then $\partial \Omega$ and $\Omega$ are bounded (by the weak $(1,1)$ boundedness of $\mathcal M$) and we take $g=f$ and $I = \emptyset$. If $E_l \subsetneq \partial \Omega$, thanks to the Vitali lemma, we can extract from $\{B(x,\dist(x,\partial \Omega \setminus E_l)/20)\}_{x\in E_l}$ a subcollection $\{B_i= B(x_i,r_i)\}_{i\in I}$ such that $\{B_i/2\}_{i\in I}$ covers $E_l$ and $\{B_i/10\}_{i\in I}$ is non-overlapping. We choose $a^*$ large enough such that 
\[\Omega \setminus \bigcup_{x\in \partial \Omega \setminus E_l} \gamma_{a^*}(x) \subset \bigcup_{i\in I} B_i,\]
we assume that $I$ is ordered, and we construct
\[S_0 := \bigcup_{x\in \partial \Omega \setminus E_l} \gamma_{a^*}(x), \qquad D_i := B_i \setminus \left( S_0 \cup \bigcup_{j<i} B_j \right),\]
\[ b_i = f|_{D_i} - \delta \Big[ \frac{f|_{D_i}}{\delta} \Big]_{D_i \cap \Omega} \1_{D_i\cap \Omega} \quad \text{ and } \quad  g := f - \sum_{i\in I} b_i.\]
By construction, we have properties $(i)$ to $(iii)$ of the Lemma. We need to check that the decomposition satisfies $(iv)$ to $(vi)$.

\medskip

\noindent {\bf From $g,b_i$ to $\bar g,\bar b_i$.}
Actually, we just need to prove the properties $(iv)$ to $(vi)$ for
\[\bar b_i := f|_{D_i} \quad \text{ and } \quad \bar g = f|_{S_0} = \sum_{i\in I} \bar b_i.\]
instead of $b_i$ and $g$. Indeed, 
\begin{multline} \label{CZ1}
    \wt \A_1\Big( \delta \left[ \frac{f|_{D_i}}{\delta}\right]_{D_i\cap \Omega} \1_{D_i\cap \Omega} \Big)(x) \lesssim r_i \left( \fiint_{D_i\cap \Omega} \Big| \frac{f(X)}{\delta(X)}\Big|\, dX \right) \1_{9B_i \cap \partial \Omega}(x) \\
    \lesssim r_i^{1-n} \left( \iint_{D_i\cap \Omega} \Big| \frac{f(X)}{\delta(X)}\Big|\, dX \right) \1_{9B_i \cap \partial \Omega}(x)
    \lesssim r_i^{1-n} \|f\1_{D_i}\|_{\wt T^1_1(\Omega)} \1_{9B_i \cap \partial \Omega}(x) \\ \lesssim r_i^{(1-n)/q}  \|\bar b_i\|_{\wt T^q_1(\Omega)} \1_{9B_i \cap \partial \Omega}(x)
\end{multline}
where the first and forth inequalities are because $\wt \A_1(\1_{D_i}$ is supported in $9B_i\cap \partial \Omega$, and the second one is because $|D_i| \geq  |\Omega \cap B_i/10| \gtrsim r_i^n$. Furthermore, we have
\begin{multline*} \Big\|\delta \Big[ \frac{f|_{D_i}}{\delta} \Big]_{D_i \cap \Omega} \1_{D_i \cap \Omega} \Big\|_{\wt T^q_1(\Omega)} \lesssim r_i^{(n-1)/q} \sup_{x\in 9B_i \cap \partial \Omega} \wt \A_1\Big( \delta \left[ \frac{f|_{D_i}}{\delta}\right]_{D_i\cap \Omega} \1_{D_i\cap \Omega} \Big)(x) \lesssim \|\bar b_i\|_{\wt T^q_1(\Omega)}
\end{multline*}
by \eqref{CZ1}. These two last computations show that 
\[\|b_i\|_{\wt T^q_1(\Omega)} \lesssim \|\bar b_i\|_{\wt T^q_1(\Omega)}\]
and, since $\{9B_i \cap \partial \Omega\}_{i\in I}$ is finitely overlapping,
\[\|g - \bar g\|_{\wt T^\infty_1(\Omega)} \lesssim \sup_{i\in I} \|b_i - \bar b_i\|_{\wt T^\infty_1(\Omega)} \lesssim \sup_{i\in I} \|\bar b_i\|_{\wt T^q_1(\Omega)}.\]
That is (iv)-(v) are satisfied for $g,b_i$ as long as they are satisfied for $\bar g,\bar b_i$.

\medskip

\noindent {\bf Proof of $(iv)$ for $\bar g$.} Let $x\in \partial \Omega$. If $x\notin E_\lambda$, then of course
\[\wt \A_1(\bar g) \leq \wt\A_1(f) \leq |\mathcal M(\wt \A_1^*(f)^q)|^{1/q} \leq l\]
by definition of $l$. If $x\in B_i \cap \partial \Omega$, and $X\in \gamma(x)$, then by construction, $\left(\fiint_{B_X} |\bar g(Y)|^2\, dY\right)^\frac12$ is non-zero only when $\delta(X) \gtrsim r_i$. It means that $\wt \A_1(\bar g)(x) \leq \wt \A_1^*(g)(y) \leq \wt \A_1^*(f)(y)$, for $y$ in a small ball $B(x,\epsilon r_i)$ around $x$. Take now $x'\in \partial \Omega \setminus E_l$ such that $|x-x'| = \dist(x,\partial \Omega \setminus E_l) \approx r_i$, then
\[\wt \A_1(\bar g)(x) \leq \fint_{B(x,\epsilon r_i)} \wt \A_1^*(f)(y) \, d\sigma \lesssim \left( \fint_{B(x',200r_i)} |\wt \A_1^*(f)(y)|^q \, d\sigma\right)^\frac1q  \leq \mathcal M[\wt \A_1^*(f)](x') \leq l\]
since $x' \in \partial \Omega \setminus E_l$.

\medskip

\noindent {\bf Proof of $(v)$ for $\bar b_i$.} Take $i\in I$ and $x'_i\in \partial \Omega \setminus E_l$ such that $\dist(B_i,x'_i) \leq 20$, which is possible by construction of $B_i$. Then 
\begin{multline*}
\|\bar b_i\|_{\wt T^q_1(\Omega)} = \|\wt \A_1(f\1_{D_i})\|_{L^q(9B_i \cap \partial \Omega)} \leq \|\wt \A_1^*(f)\|_{L^q(9B_i \cap \partial \Omega)} \leq \|\wt \A_1^*(f)\|_{L^q(B(x'_i,40r_i) \cap \partial \Omega)} \\
\lesssim r^{(n-1)/q} [\mathcal M(|\wt \A_1^*(f)|^q)]^{1/q}(x'_i) \leq r^{(n-1)/q} l
\end{multline*}
since $x'_i \in E_l$.

\medskip

\noindent {\bf Proof of $(vi)$.} That is a simple consequence of the weak $(1,1)$ boundedness of $\mathcal M$. Indeed,
\[\sum_{i\in I} r_i^{n-1} \lesssim \sum_{i\in I} \sigma(B_i/10) \leq \sigma(E_l) \lesssim l^{-q} \|\wt \A^*_1(f)\|_{L^q(\partial \Omega)}^q \lesssim l^{-q} \|f\|_{\wt T^q_1(\Omega)}\]
by Lemma \ref{lemNCA}. The lemma follows.
\ep
\bigskip

We are now ready for the proof of Theorem \ref{Thinterpolation}.
\smallskip

\noindent{\em Proof of Theorem \ref{Thinterpolation}}: Our result is an adaptation pf \cite[Theorem 1.1]{Auscher07}, which is itself inspired from \cite{BK03,DMc99}. It suffices to show that for any $f\in L_c^\infty(\om)$ with $[f/\delta]_{\om}=0$, there exists $C$ such that
\begin{equation}\label{eq.interpo1}
    \sigma\set{\wt\N(\Z(f))>l}\le Cl^{-1}\norm{f}_{\wt T_1^1} \qquad\text{for all }l>0.
\end{equation}
Indeed, we can then interpolate between the $L^{1,\infty}-\wt T^1_1$ and the $L^p-\wt T^p_1$ boundedness of $\wt \N(\mathcal Z(\cdot))$ to get the result, since $[L^{1,\infty},L^p]_{\theta,q} = L^q$ by \cite[Theorem 5.3.1]{BLinterpolation}, and $[T^1_1,T^p_1]_{\theta,q} = T^{q}_1$ by adapting the proof of \cite[Theorem 4']{CMS85}\footnote{We can also see the tent spaces $T^p_1$ as a {\em retracts} - see \cite[Definition 6.4.1]{BLinterpolation} - of $L^p(L^1)$, and use interpolation of $L^p$ on Banach spaces.}. Fix any $l>0$, we use the Calder\'on-Zygmund decomposition (Lemma \ref{lemCZ}) for $f\in \wt T_1^1$ at height $l$, that is, we write 
\(
f=g+\sum_{i\in I} b_i,
\)
where $g$ and $\set{b_i}_{i\in I}$ satisfy properties $(i)-(vi)$ with $q=1$ in Lemma \ref{lemCZ}. 
Then we only need to show that $\sigma\set{\wt\N(\Z(g))>l/2}$ and $\sigma\set{\wt\N(\Z(f-g))>l/2}$ are both bounded above by $Cl^{-1}\norm{f}_{\wt T_1^1}$.

Observe that $\|g\|_{\wt T^1_1} \leq C \|f\|_{\wt T^1_1}$ by (i), (v) and (vi) of the decomposition, and that $\|\wt \A_1(g)\|_\infty \leq Cl$ by (iv). Therefore, we have 
\[
    \norm{\wt\A_1(g)}_p^p\le Cl^{p-1}\norm{g}_{\wt T_1^1}\le Cl^{p-1}\norm{f}_{\wt T_1^1}.
\]
Moreover, $[g/\delta]_\om=0$ because $[f/\delta]_\om=0$ and $[b_i/\delta]_\om=0$, and so the boundedness of $\Z$ from $\wh{\wt T_1^p}$ to $\wt T_\infty^p$ gives that 
\[
\norm{\wt\N(\Z(g))}_p^p\le C\norm{\wt\A_1(g)}_p^p\le C l^{p-1}\norm{f}_{\wt T_1^1}.
\]
Chebyshev's inequality entails now that
\[
\sigma\set{\wt\N(\Z(g))>l/2}\le Cl^{-p}\norm{\wt\N(\Z(g))}_p^p\le Cl^{-1}\norm{f}_{\wt T_1^1},
\]
as desired in \eqref{eq.interpo1}.

Next we estimate the contribution from $\sum_ib_i$. We write
\begin{multline*}
    \sigma\set{\wt\N\br{\Z\left(\sum_{i\in I}b_i\right)}>l/2}\le \sigma\set{\sum_{i\in I}\wt\N(\Z(b_i))>l/2}\\
    \le \sigma\set{\cup_{i\in I}2^7B_i\cap \pom}
    +\sigma\set{x\in\pom\setminus\br{\cup_{i\in I}2^7B_i}: \sum_{i\in I}\wt\N(\Z(b_i))(x)>l/2}=:E+F.
\end{multline*}
The estimate for $E$ is easy: thanks to property $(vi)$ of Lemma \ref{lemCZ} and the fact that $\pom$ is $(n-1)$-Ahlfors regular, we have that 
\[
E\le C\sum_{i\in I}r_i^{n-1}\le l^{-1}\norm{f}_{\wt T_1^1}
\]
as desired. 
To estimate $F$, we first use  Chebyshev's inequality and duality to get the existence of a $\eta \in L^{p'}(\partial \Omega)$ such that $\supp\eta\subset\pom\setminus\br{\cup_{i\in I}2^7B_i}$, $\norm{\eta}_{p'}=1$, and 
\[
F\le Cl^{-p}\int_{\pom\setminus\br{\cup_i 2^7B_i}}\br{\sum_{i\in I}\wt\N(\Z(b_i))}^pd\sigma 
=Cl^{-p}\br{\int\br{\sum_{i\in I}\wt\N(\Z(b_i))}\eta\,d\sigma }^p.
\]
Denote $C_j(B):=2^{j+1}B\setminus 2^j B$. We further write
\begin{multline*}
    \int\br{\sum_{i\in I}\wt\N(\Z(b_i))}\eta\,d\sigma
    \le \sum_{i\in I}\sum_{j=5}^\infty\int  \wt\N\br{\Z(b_i)\1_{C_j(B_i)}}\eta\,d\sigma + \sum_{i\in I}\int \wt\N\br{\Z(b_i)\1_{2^5B_i}}\eta\,d\sigma\\
    =\sum_{i\in I}\sum_{j=5}^\infty\int  \wt\N\br{\Z(b_i)\1_{C_j(B_i)}}\eta\,d\sigma,
\end{multline*}
because $\supp\wt\N(\Z(b_i)\1_{2^5B_i})\subset 2^7 B_i$, and so 
\(\int \wt\N\br{\Z(b_i)\1_{2^5B_i}}\eta\,d\sigma=0\).
By H\"older's inequality, the observation that $\supp\wt\N(\Z(b_i)\1_{C_j(B_i)})\subset 2^{j+3}B_i$, and \eqref{interpolation1}, we get that
\begin{multline*}
    \int  \wt\N\br{\Z(b_i)\1_{C_j(B_i)}}\eta\,d\sigma 
    \le\br{\int \wt\N\br{\Z(b_i)\1_{C_j(B_i)}}^p d\sigma}^{1/p}\br{\int_{2^{j+3}B_i}\abs{\eta}^{p'}d\sigma}^{1/{p'}}\\
    \le C 2^{j(n-1)}g(j)\norm{\wt\A_1(b_i)}_1\mathcal{M}\br{\eta^{p'}}(z)^{1/{p'}} \qquad\text{for any }z\in B_i,
\end{multline*}
where $\mathcal M$ is the Hardy-Littlewood maximal function. Integrating in $z$ over $B_i$ and taking the average, one has that 
\begin{multline*}
    \int  \wt\N\br{\Z(b_i)\1_{C_j(B_i)}}\eta\,d\sigma 
\le C2^{j(n-1)}g(j)\norm{\wt\A_1(b_i)}_1\fint_{B_i\cap\pom}\mathcal{M}\br{\eta^{p'}}(z)^{1/{p'}}d\sigma(z)\\
\le Cl2^{j(n-1)}g(j)r_i^{n-1}\fint_{B_i\cap\pom}\mathcal{M}\br{\eta^{p'}}^{1/{p'}}d\sigma\le  Cl2^{j(n-1)}g(j)\int_{B_i\cap\pom}\mathcal{M}\br{\eta^{p'}}^{1/{p'}}d\sigma
\end{multline*}
as $\norm{b_i}_{\wt T_1^1}\le Clr_i^{n-1}$.
Summing in $i\in I$ and $j\ge5$, since $\sum_{j\ge3}g(j)2^{j(n-1)}<\infty$, we get that 
\[
\int\br{\sum_{i\in I}\wt\N(\Z(b_i))}\eta\,d\sigma \le Cl\int_{\cup_{i\in I}B_i\cap\pom} \mathcal{M}\br{\eta^{p'}}^{1/{p'}}d\sigma,
\]
where we have used the property that the collection $\set{B_i}_{i\in I}$ has finite overlaps.  
Since the maximal function is weak $(1,1)$, we can apply Kolmogorov's lemma (see e.g. \cite[Lemma 5.16]{Duoan}) to get that 
\begin{multline*}
    \int_{\cup_{i\in I}B_i\cap\pom}\mathcal{M}\br{\eta^{p'}}^{1/{p'}}d\sigma \le C\sigma(\cup_{i\in I}B_i\cap\pom)^{1/p}\norm{\abs{\eta}^{p'}}_1^{1/{p'}}\\
    \le
    C\br{\sum_{i\in I}r_i^{n-1}}^{1/p}\le Cl^{-1/p}\norm{f}_{\wt T_1^1}^{1/p},
\end{multline*}
by property $(vi)$ in Lemma \ref{lemCZ}.    

Altogether, we have shown $F\le Cl^{-1}\norm{f}_{\wt T_1^1}$ as desired for \eqref{eq.interpo1}.\ep

\bibliographystyle{alpha}
\bibliography{reference}

\newcommand{\etalchar}[1]{$^{#1}$}
\begin{thebibliography}{HKMP15b}

\bibitem[AA11]{AA11}
Pascal Auscher and Andreas Axelsson.
\newblock Weighted maximal regularity estimates and solvability of non-smooth
  elliptic systems {I}.
\newblock {\em Invent. Math.}, 184(1):47--115, 2011.

\bibitem[AAA{\etalchar{+}}11]{AAAHK11}
M.~Angeles Alfonseca, Pascal Auscher, Andreas Axelsson, Steve Hofmann, and
  Seick Kim.
\newblock Analyticity of layer potentials and {$L^2$} solvability of boundary
  value problems for divergence form elliptic equations with complex
  {$L^\infty$} coefficients.
\newblock {\em Adv. Math.}, 226(5):4533--4606, 2011.

\bibitem[AHM{\etalchar{+}}20]{AHMMT20}
Jonas Azzam, Steve Hofmann, Jos\'{e}~Mar\'{\i}a Martell, Mihalis Mourgoglou,
  and Xavier Tolsa.
\newblock Harmonic measure and quantitative connectivity: geometric
  characterization of the {$L^p$}-solvability of the {D}irichlet problem.
\newblock {\em Invent. Math.}, 222(3):881--993, 2020.

\bibitem[AM14]{AM14}
Pascal Auscher and Mihalis Mourgoglou.
\newblock Boundary layers, {R}ellich estimates and extrapolation of solvability
  for elliptic systems.
\newblock {\em Proc. Lond. Math. Soc. (3)}, 109(2):446--482, 2014.

\bibitem[AR12]{AR12}
Pascal Auscher and Andreas Ros\'{e}n.
\newblock Weighted maximal regularity estimates and solvability of nonsmooth
  elliptic systems, {II}.
\newblock {\em Anal. PDE}, 5(5):983--1061, 2012.

\bibitem[Aus07]{Auscher07}
Pascal Auscher.
\newblock On necessary and sufficient conditions for {$L^p$}-estimates of
  {R}iesz transforms associated to elliptic operators on {$\Bbb R^n$} and
  related estimates.
\newblock {\em Mem. Amer. Math. Soc.}, 186(871):xviii+75, 2007.

\bibitem[BK03]{BK03}
S\"{o}nke Blunck and Peer~Christian Kunstmann.
\newblock Calder\'{o}n-{Z}ygmund theory for non-integral operators and the
  {$H^\infty$} functional calculus.
\newblock {\em Rev. Mat. Iberoamericana}, 19(3):919--942, 2003.

\bibitem[BL76]{BLinterpolation}
J\"{o}ran Bergh and J\"{o}rgen L\"{o}fstr\"{o}m.
\newblock {\em Interpolation spaces. {A}n introduction}, volume No. 223 of {\em
  Grundlehren der Mathematischen Wissenschaften}.
\newblock Springer-Verlag, Berlin-New York, 1976.

\bibitem[BPTT21]{BPTT??}
Simon Bortz, Bruno Poggi, Olli Tapiola, and Xavier Tolsa.
\newblock The {$A_\infty$}-condition, {$\epsilon$}-approximators, and
  varopoulos extensions in uniform domains.
\newblock {\em Preprint}, 2021.

\bibitem[CDMT22]{CDMT22}
Mingming Cao, \'{O}scar Dom\'{\i}nguez, Jos\'{e}~Mar\'{\i}a Martell, and Pedro
  Tradacete.
\newblock On the {$A_\infty$} condition for elliptic operators in 1-sided
  nontangentially accessible domains satisfying the capacity density condition.
\newblock {\em Forum Math. Sigma}, 10:Paper No. e59, 57, 2022.

\bibitem[CHMT20]{CHMT20}
Juan Cavero, Steve Hofmann, Jos\'{e}~Mar\'{\i}a Martell, and Tatiana Toro.
\newblock Perturbations of elliptic operators in 1-sided chord-arc domains.
  {P}art {II}: non-symmetric operators and {C}arleson measure estimates.
\newblock {\em Trans. Amer. Math. Soc.}, 373(11):7901--7935, 2020.

\bibitem[CHPMar]{CHM??}
Mingming Cao, Pablo Hidalgo-Palencia, and Jos{\'e}~Mar{\'i}a Martell.
\newblock Carleson measure estimates, corona decompositions, and perturbation
  of elliptic operators without connectivity.
\newblock {\em Mathematische Annalen}, To appear.

\bibitem[CMS85]{CMS85}
R.~R. Coifman, Y.~Meyer, and E.~M. Stein.
\newblock Some new function spaces and their applications to harmonic analysis.
\newblock {\em J. Funct. Anal.}, 62(2):304--335, 1985.

\bibitem[Dah77]{Dah77}
Bj\"{o}rn E.~J. Dahlberg.
\newblock Estimates of harmonic measure.
\newblock {\em Arch. Rational Mech. Anal.}, 65(3):275--288, 1977.

\bibitem[DDE{\etalchar{+}}24]{DDEMM??}
Guy David, Stephano Decio, Max Engelstein, Svitlana Mayboroda, and Marco
  Michetti.
\newblock Robin harmonic measure on rough domains.
\newblock {\em In preparation}, 2024.

\bibitem[DFM21]{DFMhighcod}
G.~David, J.~Feneuil, and S.~Mayboroda.
\newblock Elliptic theory for sets with higher co-dimensional boundaries.
\newblock {\em Mem. Amer. Math. Soc.}, 274(1346):vi+123, 2021.

\bibitem[DFM23a]{DFMcarl}
Zanbing Dai, Joseph Feneuil, and Svitlana Mayboroda.
\newblock Carleson perturbations for the regularity problem.
\newblock {\em Rev. Mat. Iberoam.}, 39(6):2119--2170, 2023.

\bibitem[DFM23b]{DFMreg}
Zanbing Dai, Joseph Feneuil, and Svitlana Mayboroda.
\newblock The regularity problem in domains with lower dimensional boundaries.
\newblock {\em Journal of Functional Analysis}, 284(11):109903, 2023.

\bibitem[DFM23c]{DFMmixed}
Guy David, Joseph Feneuil, and Svitlana Mayboroda.
\newblock Elliptic theory in domains with boundaries of mixed dimension.
\newblock {\em Ast\'{e}risque}, (442):vi+139, 2023.

\bibitem[DHP23]{DHP23}
Martin Dindo\v{s}, Steve Hofmann, and Jill Pipher.
\newblock Regularity and {N}eumann problems for operators with real
  coefficients satisfying {C}arleson conditions.
\newblock {\em J. Funct. Anal.}, 285(6):Paper No. 110024, 32, 2023.

\bibitem[DJ90]{DJ90}
G.~David and D.~Jerison.
\newblock Lipschitz approximation to hypersurfaces, harmonic measure, and
  singular integrals.
\newblock {\em Indiana Univ. Math. J.}, 39(3):831--845, 1990.

\bibitem[DK87]{DK87}
Bj\"{o}rn E.~J. Dahlberg and Carlos~E. Kenig.
\newblock Hardy spaces and the {N}eumann problem in {$L^p$} for {L}aplace's
  equation in {L}ipschitz domains.
\newblock {\em Ann. of Math. (2)}, 125(3):437--465, 1987.

\bibitem[DK12]{DK12}
Martin Dindo\v{s} and Josef Kirsch.
\newblock The regularity problem for elliptic operators with boundary data in
  {H}ardy-{S}obolev space {$HS^1$}.
\newblock {\em Math. Res. Lett.}, 19(3):699--717, 2012.

\bibitem[DKP11]{DKP11}
Martin Dindos, Carlos Kenig, and Jill Pipher.
\newblock B{MO} solvability and the {$A_\infty$} condition for elliptic
  operators.
\newblock {\em J. Geom. Anal.}, 21(1):78--95, 2011.

\bibitem[DM99]{DMc99}
Xuan~Thinh Duong and Alan MacIntosh.
\newblock Singular integral operators with non-smooth kernels on irregular
  domains.
\newblock {\em Rev. Mat. Iberoamericana}, 15(2):233--265, 1999.

\bibitem[DP19]{DP19}
Martin Dindo\v{s} and Jill Pipher.
\newblock Perturbation theory for solutions to second order elliptic operators
  with complex coefficients and the {$L^p$} {D}irichlet problem.
\newblock {\em Acta Math. Sin. (Engl. Ser.)}, 35(6):749--770, 2019.

\bibitem[DPP07]{DPP07}
Martin Dindos, Stefanie Petermichl, and Jill Pipher.
\newblock The {$L^p$} {D}irichlet problem for second order elliptic operators
  and a {$p$}-adapted square function.
\newblock {\em J. Funct. Anal.}, 249(2):372--392, 2007.

\bibitem[DPR17]{DPR17}
Martin Dindo\v{s}, Jill Pipher, and David Rule.
\newblock Boundary value problems for second-order elliptic operators
  satisfying a {C}arleson condition.
\newblock {\em Comm. Pure Appl. Math.}, 70(7):1316--1365, 2017.

\bibitem[Duo01]{Duoan}
Javier Duoandikoetxea.
\newblock {\em Fourier analysis}, volume~29 of {\em Graduate Studies in
  Mathematics}.
\newblock American Mathematical Society, Providence, RI, 2001.
\newblock Translated and revised from the 1995 Spanish original by David
  Cruz-Uribe.

\bibitem[Fen22]{Fen22}
Joseph Feneuil.
\newblock Absolute continuity of the harmonic measure on low dimensional
  rectifiable sets.
\newblock {\em J. Geom. Anal.}, 32(10):Paper No. 247, 36, 2022.

\bibitem[Fen23]{Fenreg}
Joseph Feneuil.
\newblock An alternative proof of the {$L^p$}-regularity problem for
  dahlberg-kenig-pipher operators on {$\mathbb R^n_+$}.
\newblock {\em Preprint, ArXiv:2310.00645}, 2023.

\bibitem[FP22]{FP23}
Joseph Feneuil and Bruno Poggi.
\newblock Generalized {C}arleson perturbations of elliptic operators and
  applications.
\newblock {\em Trans. Amer. Math. Soc.}, 375(11):7553--7599, 2022.

\bibitem[Gia83]{Giaquinta}
Mariano Giaquinta.
\newblock {\em Multiple integrals in the calculus of variations and nonlinear
  elliptic systems}, volume 105 of {\em Annals of Mathematics Studies}.
\newblock Princeton University Press, Princeton, NJ, 1983.

\bibitem[GMT18]{GMT18}
John Garnett, Mihalis Mourgoglou, and Xavier Tolsa.
\newblock Uniform rectifiability from {C}arleson measure estimates and
  {$\varepsilon$}-approximability of bounded harmonic functions.
\newblock {\em Duke Math. J.}, 167(8):1473--1524, 2018.

\bibitem[GMT23]{GMT??}
Josep Gallegos, Mihalis Mourgoglou, and Xavier Tolsa.
\newblock Extrapolation of solvability of the regularity and the poisson
  regularity problems in rough domains.
\newblock {\em Preprint, Arxiv:2306.06185}, 2023.

\bibitem[GW82]{GW82}
Michael Gr\"{u}ter and Kjell-Ove Widman.
\newblock The {G}reen function for uniformly elliptic equations.
\newblock {\em Manuscripta Math.}, 37(3):303--342, 1982.

\bibitem[Ha96]{Haj96}
Piotr Haj\l~asz.
\newblock Sobolev spaces on an arbitrary metric space.
\newblock {\em Potential Anal.}, 5(4):403--415, 1996.

\bibitem[HaK95]{HK95}
Piotr Haj\l~asz and Pekka Koskela.
\newblock Sobolev meets {P}oincar\'{e}.
\newblock {\em C. R. Acad. Sci. Paris S\'{e}r. I Math.}, 320(10):1211--1215,
  1995.

\bibitem[HKMP15a]{HKMP15b}
Steve Hofmann, Carlos Kenig, Svitlana Mayboroda, and Jill Pipher.
\newblock The regularity problem for second order elliptic operators with
  complex-valued bounded measurable coefficients.
\newblock {\em Math. Ann.}, 361(3-4):863--907, 2015.

\bibitem[HKMP15b]{HKMP15a}
Steve Hofmann, Carlos Kenig, Svitlana Mayboroda, and Jill Pipher.
\newblock Square function/non-tangential maximal function estimates and the
  {D}irichlet problem for non-symmetric elliptic operators.
\newblock {\em J. Amer. Math. Soc.}, 28(2):483--529, 2015.

\bibitem[HL18]{HL18}
Steve Hofmann and Phi Le.
\newblock B{MO} solvability and absolute continuity of harmonic measure.
\newblock {\em J. Geom. Anal.}, 28(4):3278--3299, 2018.

\bibitem[HLM19]{HLM19}
Steve Hofmann, Phi Le, and Andrew~J. Morris.
\newblock Carleson measure estimates and the {D}irichlet problem for degenerate
  elliptic equations.
\newblock {\em Anal. PDE}, 12(8):2095--2146, 2019.

\bibitem[HLMP22]{HLMP22}
Steve Hofmann, Linhan Li, Svitlana Mayboroda, and Jill Pipher.
\newblock The {D}irichlet problem for elliptic operators having a {BMO}
  anti-symmetric part.
\newblock {\em Math. Ann.}, 382(1-2):103--168, 2022.

\bibitem[HM14]{HM14}
Steve Hofmann and Jos\'{e}~Mar\'{\i}a Martell.
\newblock Uniform rectifiability and harmonic measure {I}: {U}niform
  rectifiability implies {P}oisson kernels in {$L^p$}.
\newblock {\em Ann. Sci. \'{E}c. Norm. Sup\'{e}r. (4)}, 47(3):577--654, 2014.

\bibitem[HMM16]{HMM16}
Steve Hofmann, Jos\'{e}~Mar\'{\i}a Martell, and Svitlana Mayboroda.
\newblock Uniform rectifiability, {C}arleson measure estimates, and
  approximation of harmonic functions.
\newblock {\em Duke Math. J.}, 165(12):2331--2389, 2016.

\bibitem[HMM{\etalchar{+}}21]{HMMTZ21}
Steve Hofmann, Jos\'{e}~Mar\'{\i}a Martell, Svitlana Mayboroda, Tatiana Toro,
  and Zihui Zhao.
\newblock Uniform rectifiability and elliptic operators satisfying a {C}arleson
  measure condition.
\newblock {\em Geom. Funct. Anal.}, 31(2):325--401, 2021.

\bibitem[HMTon]{HMTbook14}
Steve Hofmann, Jos{\'e}~Mar{\'i}a Martell, and Tatiana Toro.
\newblock {\em General divergence form elliptic operators on domains with ADR
  boundaries}.
\newblock Book in preparation.

\bibitem[HR13]{HR13}
Tuomas Hyt\"{o}nen and Andreas Ros\'{e}n.
\newblock On the {C}arleson duality.
\newblock {\em Ark. Mat.}, 51(2):293--313, 2013.

\bibitem[HR18]{HR18}
Tuomas Hyt\"{o}nen and Andreas Ros\'{e}n.
\newblock Bounded variation approximation of {$L_p$} dyadic martingales and
  solutions to elliptic equations.
\newblock {\em J. Eur. Math. Soc. (JEMS)}, 20(8):1819--1850, 2018.

\bibitem[HS24]{Steve'student}
Steve Hofmann and Derek Sparrius.
\newblock The neumann function in 1-sided chord arc domains, and extrapolation
  of solvability for the neumann problem.
\newblock {\em In preparation}, 2024.

\bibitem[Hua17]{Huang17}
Yi~Huang.
\newblock Calder\'{o}n-{Z}ygmund decompositions in tent spaces and weak-type
  endpoint bounds for two quadratic functionals of {S}tein and
  {F}efferman-{S}tein.
\newblock {\em Studia Math.}, 239(2):123--132, 2017.

\bibitem[JK81a]{JK81a}
David~S. Jerison and Carlos~E. Kenig.
\newblock The {D}irichlet problem in nonsmooth domains.
\newblock {\em Ann. of Math. (2)}, 113(2):367--382, 1981.

\bibitem[JK81b]{JK81b}
David~S. Jerison and Carlos~E. Kenig.
\newblock The {N}eumann problem on {L}ipschitz domains.
\newblock {\em Bull. Amer. Math. Soc. (N.S.)}, 4(2):203--207, 1981.

\bibitem[JK82]{JK82a}
David~S. Jerison and Carlos~E. Kenig.
\newblock Boundary behavior of harmonic functions in nontangentially accessible
  domains.
\newblock {\em Adv. in Math.}, 46(1):80--147, 1982.

\bibitem[Ken94]{Ken94}
Carlos~E. Kenig.
\newblock {\em Harmonic analysis techniques for second order elliptic boundary
  value problems}, volume~83 of {\em CBMS Regional Conference Series in
  Mathematics}.
\newblock Conference Board of the Mathematical Sciences, Washington, DC; by the
  American Mathematical Society, Providence, RI, 1994.

\bibitem[KKPT00]{KKPT00}
C.~Kenig, H.~Koch, J.~Pipher, and T.~Toro.
\newblock A new approach to absolute continuity of elliptic measure, with
  applications to non-symmetric equations.
\newblock {\em Adv. Math.}, 153(2):231--298, 2000.

\bibitem[KKPT16]{KKPT16}
C.~Kenig, B.~Kirchheim, J.~Pipher, and T.~Toro.
\newblock Square functions and the {$A_\infty$} property of elliptic measures.
\newblock {\em J. Geom. Anal.}, 26(3):2383--2410, 2016.

\bibitem[KP93]{KP93}
Carlos~E. Kenig and Jill Pipher.
\newblock The {N}eumann problem for elliptic equations with nonsmooth
  coefficients.
\newblock {\em Invent. Math.}, 113(3):447--509, 1993.

\bibitem[KP95]{KP95}
Carlos~E. Kenig and Jill Pipher.
\newblock The {N}eumann problem for elliptic equations with nonsmooth
  coefficients. {II}.
\newblock volume~81, pages 227--250. 1995.
\newblock A celebration of John F. Nash, Jr.

\bibitem[KP01]{KP01}
Carlos~E. Kenig and Jill Pipher.
\newblock The {D}irichlet problem for elliptic equations with drift terms.
\newblock {\em Publ. Mat.}, 45(1):199--217, 2001.

\bibitem[KR09]{KR08}
Carlos~E. Kenig and David~J. Rule.
\newblock The regularity and {N}eumann problem for non-symmetric elliptic
  operators.
\newblock {\em Trans. Amer. Math. Soc.}, 361(1):125--160, 2009.

\bibitem[KS08]{KS08}
Aekyoung~Shin Kim and Zhongwei Shen.
\newblock The {N}eumann problem in {$L^p$} on {L}ipschitz and convex domains.
\newblock {\em J. Funct. Anal.}, 255(7):1817--1830, 2008.

\bibitem[KS11]{KS11reg}
Joel Kilty and Zhongwei Shen.
\newblock The {$L^p$} regularity problem on {L}ipschitz domains.
\newblock {\em Trans. Amer. Math. Soc.}, 363(3):1241--1264, 2011.

\bibitem[MPT22]{MPT}
Mihalis Mourgoglou, Bruno Poggi, and Xavier Tolsa.
\newblock Solvability of the {P}oisson-{D}irichlet problem with interior data
  in {$L^{p'}$}-carleson spaces and its applications to the {$L^p$}-regularity
  problem.
\newblock {\em arXiv preprint arXiv:2207.10554}, 2022.

\bibitem[MT03]{MT03}
Marius Mitrea and Michael Taylor.
\newblock Potential theory on {L}ipschitz domains in {R}iemannian manifolds:
  the case of {D}ini metric tensors.
\newblock {\em Trans. Amer. Math. Soc.}, 355(5):1961--1985, 2003.

\bibitem[MTar]{MT??}
Mihalis Mourgoglou and Xavier Tolsa.
\newblock The regularity problem for the laplace equation in rough domains.
\newblock {\em Duke Math. J.}, To appear.

\bibitem[MZ23]{MZ??}
Mihalis Mourgoglou and Thanasis Zacharopoulos.
\newblock Varopoulos extensions in domains with ahlfors-regular boundaries and
  applications to boundary value problems for elliptic systems with
  {$L^\infty$} coefficients.
\newblock {\em Preprint, Arxiv:2303.10717}, 2023.

\bibitem[She07]{Shen07}
Zhongwei Shen.
\newblock A relationship between the {D}irichlet and regularity problems for
  elliptic equations.
\newblock {\em Math. Res. Lett.}, 14(2):205--213, 2007.

\bibitem[Ste93]{Stein93}
Elias~M. Stein.
\newblock {\em Harmonic analysis: real-variable methods, orthogonality, and
  oscillatory integrals}, volume~43 of {\em Princeton Mathematical Series}.
\newblock Princeton University Press, Princeton, NJ, 1993.
\newblock With the assistance of Timothy S. Murphy, Monographs in Harmonic
  Analysis, III.

\end{thebibliography}

\end{document}